\documentclass[a4paper,preprint]{elsarticle}

\usepackage{hyperref}

\journal{}

\usepackage[margin=2.5cm]{geometry}
\usepackage{amsmath,amssymb}
\usepackage{amsfonts}
\usepackage{fancybox}
\usepackage[english,francais]{babel}
\usepackage{subfigure}
\usepackage{float}
\usepackage{comment}

\newtheorem{e-proposition}[theorem]{Proposition}

\newtheorem{e-definition}[theorem]{Definition\rm}
\newtheorem{remark}{\it Remark\/}

\graphicspath{{figures/}}

\newcommand{\intO}{\int_{\Omega}}

\newcommand{\bn}{\boldsymbol{n}}

\newcommand{\bx}{\boldsymbol{x}}

\newcommand{\bs}{\boldsymbol{s}}
\newcommand{\bw}{\boldsymbol{w}}
\newcommand{\bq}{\boldsymbol{q}}

\newcommand{\bp}{\boldsymbol{p}}

\newcommand{\be}{\boldsymbol{e}}
\newcommand{\nab}{\boldsymbol{\nabla}}
\newcommand{\RR}{\mathbb{R}}
\newcommand{\II}{\mathbb{I}}
\newcommand{\Aa}{\mathbb{A}}

\newcommand{\mT}{\mathcal{T}}
\newcommand{\mR}{\mathcal{R}}
\newcommand{\mS}{\mathcal{S}}

\newcommand{\mF}{\mathcal{F}}
\newcommand{\mH}{\mathcal{H}}
\newcommand{\mI}{\mathcal{I}}
\newcommand{\eps}{\varepsilon}

\newcommand{\dis}{\displaystyle}

\newcommand{\vertiii}[1]{{\left\vert\kern-0.25ex\left\vert\kern-0.25ex\left\vert #1 
    \right\vert\kern-0.25ex\right\vert\kern-0.25ex\right\vert}}

\begin{document}

\begin{frontmatter}

\title{Goal-oriented error estimation and adaptivity in MsFEM computations}

\author[1,2]{Ludovic Chamoin\corref{mycorrespondingauthor}}
\cortext[mycorrespondingauthor]{Corresponding author}
\ead{ludovic.chamoin@ens-paris-saclay.fr}

\author[3,2]{Fr\'ed\'eric Legoll}
\ead{frederic.legoll@enpc.fr}

\address[1]{LMT (ENS Paris-Saclay, CNRS, Universit\'e Paris-Saclay), 61 avenue du Pr\'esident Wilson, 94235 Cachan, France}
\address[2]{Inria Paris, MATHERIALS project-team, 2 rue Simone Iff, CS 42112, 75589 Paris Cedex 12, France}
\address[3]{Laboratoire Navier (ENPC, Paris-Est University), 6-8 avenue Blaise Pascal, 77455 Marne-La-Vall\'ee Cedex 2, France}

\begin{abstract}
We introduce a goal-oriented strategy for multiscale computations performed using the Multiscale Finite Element Method (MsFEM). In a previous work, we have shown how to use, in the MsFEM framework, the concept of Constitutive Relation Error (CRE) to obtain a guaranteed and fully computable \textit{a posteriori} error estimate in the energy norm (as well as error indicators on various error sources). Here, the CRE concept is coupled with the solution of an adjoint problem to control the error and drive an adaptive procedure with respect to a given output of interest. Furthermore, a local and non-intrusive enrichment technique is proposed to enhance the accuracy of error bounds. The overall strategy, which is fully automatic and robust, enables to reach an appropriate trade-off between certified reliability and computational cost in the MsFEM context. The performances of the proposed method are investigated on several illustrative numerical test cases. In particular, the error estimation is observed to be very accurate, yielding a very efficient adaptive procedure.
\end{abstract}


\end{frontmatter}


\section{Introduction and objectives}\label{section:intro}

Highly heterogeneous materials are ubiquitous in Material Sciences and Computational Mechanics. Examples include porous or perforated media, composite materials, \dots Developing multiscale numerical methods to simulate the (thermal, mechanical, \dots) response of such materials is an increasingly active research field. As an alternative to solving a full fine-scale problem (the cost of which is usually prohibitively expensive), multiscale modeling aims at linking the different scales by appropriately capturing the effect of the smaller scales on the larger scales. This enables an accurate description of complex physical phenomena and/or the prediction of macroscopic properties at a limited computational cost. In the recent years, significant progress has been witnessed in the design of multiscale methods for elliptic PDEs that take into account scale separation between macroscopic and microscopic features. We can distinguish two families of multiscale methods. The first family includes approaches relying on multilevel finite elements, such as HMM or FE$^2$ methods~\cite{E03,FEY03}. They perform a numerical homogenization of the medium in order to find local effective material properties, which are next inserted in a global macroscopic problem which can be handled with a reduced number of degrees of freedom. Such approaches can be considered as numerical generalizations of classical homogenization approaches (see e.g.~\cite{BEN78}) to non-periodic media. The second family includes approaches based on superposition of scales, such as the GFEM~\cite{STR01}, the MsFEM~\cite{HOU97}, the LOD~\cite{LOD} and the VMS~\cite{HUG98} methods, to name but a few. In this work, we focus on the Multiscale Finite Element Method (MsFEM) which has been the object of numerous research works (see~\cite{EFE09} for an overview). This method is based on the classical FEM framework, defining an approximate solution in a finite dimension space related to a macroscopic mesh, but uses specific basis functions which encode fine-scale details. The basis functions are obtained, in a preliminary offline stage, by solving fine-scale equations over local subdomains. The MsFEM approximation is next computed during an inexpensive online stage, and provides a numerical approximation that describes details of the exact solution at different scales, including the small ones.

\medskip

A recurring concern when using multiscale methods is the control of the accuracy of the numerical solution. In particular, it is critical to understand, quantify and control how errors generated by different sources and at different scales affect the overall error, in order to build relevant adaptive strategies and focus the computational efforts where they are required. Following the large literature on \textit{a posteriori} error estimation in the monoscale framework of the Finite Element Method (FEM) (see e.g.~\cite{VER96,AIN00} for classical textbooks and~\cite{CHA16a} for a recent review), several verification tools have been developed for multiscale methods~\cite{STR06,LAR05,LAR07,NOL08,ABD09,JHU12,ABD13a,HEN14,CHU16}. However, many open questions remain, such as quantitative assessment of error propagation across scales, relevant adaptive strategies, \dots Recently, we have proposed a robust \textit{a posteriori} error estimation tool for MsFEM computations~\cite{CHA16b}. It uses the concept of Constitutive Relation Error (CRE), which has been developed in the Computational Mechanics community for more than thirty years~\cite{LAD04}. It enables both (i) to derive fully computable (i.e. without any unknown multiplicative constant) and guaranteed bounds on the global error measured in the energy norm, which may further be used as stopping criteria in adaptive algorithms, from an extension of the CRE to the multiscale context, and (ii) to assess the various MsFEM error sources and drive an adaptive algorithm from local error maps in order to reach a prescribed error tolerance.

In the present article, we extend the methodology introduced in~\cite{CHA16b} to \textit{a posteriori} goal-oriented error estimation, namely the control of the error on specific scalar-valued outputs of interest for problems discretized with MsFEM. In this specific context, material heterogeneities may indeed require standard adaptive methods to add degrees of freedom in regions with limited influence on the output of interest, or not to add some in regions which significantly affect this output. Consequently, alternative \textit{a posteriori} error estimation and adaptivity methods need to be introduced, taking into account the specific quantities of interest of the problem. For that purpose, we introduce classical duality arguments along with an adjoint problem~\cite{BEC96,PAR97,RAN97,CIR98,PRU99,ODE01,GIL02,CAO03}, which enable to target the error estimation and adaptivity procedure towards the approximation of the considered quantity of interest. We also use strategies to bound the error that are similar to those employed for monoscale applications with CRE~\cite{LAD08a,CHA08,LAD10,PAN10,WAE12,LAD13}. In particular, we introduce a non-intrusive approach for the approximate solution of the adjoint solution (using so-called handbook functions), which enables to obtain both guaranteed and sharp error bounds. We also define local indicators on various error sources (coarse mesh size, size of the local fine meshes, \dots) in order to drive the adaptive strategy, and automatically and iteratively select optimal MsFEM parameters to be used locally for reducing the error in the goal functional. In order to comply with the MsFEM philosophy and the multi-query paradigm, all additional costly computations associated with the proposed procedures are performed in the offline stage. The overall verification tool enables to perform MsFEM computations in which: (i) a given precision on specified outputs of interest is attained at the end of the simulation; (ii) the computational work which is needed is as small as possible. The proposed approach, devoted to the Multiscale Finite Element Method and which fulfills the above two objectives, is an alternative to the goal-oriented adaptive strategy recently proposed in~\cite{CHU16} within the Generalized Multiscale Finite Element framework, which uses residual-based techniques and which mainly focuses on goal-oriented adaptivity alone. We also mention the approach~\cite{LAR05} proposed in the VMS framework, and the approach~\cite{MAI16} proposed in the HMM framework where the overall error is also split between different contributions each related to a discretization parameter.

The performances of our approach are evaluated on several numerical examples. It is interesting to note that it leads in practice to MsFEM computations with guaranteed prescribed accuracy, and in which fine-scale features of the solution are required (and captured using fine discretization parameters) only in very localized parts of the physical domain. A limited accuracy of the fine scale computations is thus sufficient. Consequently, it naturally and automatically shares some resemblance with alternative multiscale computing strategies based on local enrichments or couplings between fine-scale and coarser-scale (upscaled) models~\cite{LIO99,BEN01,BRE01,GEN09,LOZ11}, these strategies being associated to modeling error estimation~\cite{ODE99b,ODE00,VEM01,ODE02,PRU09,AKB15}.

\medskip

This article is organized as follows. In Section~\ref{section:MsFEMbasics}, we introduce the model problem we consider along with basic notions on the MsFEM discretization approach. The CRE concept applied to the MsFEM framework is recalled in Section~\ref{section:CREconcept}. In Section~\ref{section:goal-oriented}, we extend this concept to goal-oriented error estimation within MsFEM computations. This is the core of this work. An associated adaptive strategy, enabling to control the accuracy of outputs of interest, is introduced in Section~\ref{section:adaptive}. We discuss several illustrative numerical experiments in Section~\ref{section:results}. Eventually, conclusions and perspectives are collected in Section~\ref{section:conclusions}.

\section{Model problem and MsFEM approximation}\label{section:MsFEMbasics}

To start with, we describe our model problem in Section~\ref{section:modelpb}. In Section~\ref{section:MsFEMapprox}, we next recall some basic elements on several variants of MsFEM: the initial conforming variant, the oversampling variant (see Section~\ref{section:oversampling}) and the higher-order variant introduced by G. Allaire and R. Brizzi (see Section~\ref{section:highorder}). The reader familiar with MsFEM may easily skip Section~\ref{section:MsFEMapprox} and directly proceed to the {\it a posteriori} error estimation strategy described in Section~\ref{section:CREconcept}.

\subsection{Model problem}\label{section:modelpb}

We consider mathematical models described by partial differential equations (PDEs) and characterized by multiscale material heterogeneities. For the sake of simplicity, we focus on the following elliptic problem posed on some open and bounded domain $\Omega \subset \RR^d$:
\begin{equation}\label{eq:refpb}
- \nab \cdot (\Aa^\eps \nab u^\eps) = f \ \ \text{in $\Omega$}, \qquad 
u^\eps = 0 \ \ \text{on $\Gamma_D$}, \qquad
\Aa^\eps \nab u^\eps \cdot \bn = g \ \ \text{on $\Gamma_N$}.
\end{equation}
The boundary $\Gamma = \partial \Omega$ is split into a subset $\Gamma_D \subset \partial \Omega$ (with measure $|\Gamma_D|>0$) where homogeneous Dirichlet boundary conditions are applied, and a subset $\Gamma_N = \partial \Omega \setminus \Gamma_D \subset \partial \Omega$ where Neumann boundary conditions are applied. We assume that $f \in L^2(\Omega)$ and $g \in L^2(\Gamma_N)$, corresponding to body and traction loadings respectively, are slowly varying functions. The quantity $\bq^\eps = \Aa^\eps \nab u^\eps$ is the flux associated to $u^\eps$ through the constitutive relation.

The second-order diffusion tensor $\Aa^\eps \in (L^\infty(\Omega))^{d\times d}$ is a heterogeneous object, with rapid oscillations and possibly high contrast. The parameter $\eps$ refers to the typical (presumably small) length scale of the heterogeneities. We assume that $\Aa^\eps$ is a symmetric matrix that is uniformly elliptic and bounded, in the sense that there exist $\alpha >0$ and $\beta >0$ such that
$$
\forall \eps \ge 0, \ \ \forall \boldsymbol{\xi}\in\RR^d, \quad \alpha |\boldsymbol{\xi}|^2 \le \Aa^\eps(\bx) \boldsymbol{\xi} \cdot \boldsymbol{\xi} \le \beta |\boldsymbol{\xi}|^2 \quad \text{a.e. in $\Omega$.}
$$
Introducing the functional space $V=\{v \in \mH^1(\Omega), \ \ \text{$v=0$ on $\Gamma_D$} \}$, the weak formulation of~\eqref{eq:refpb} consists in finding $u^\eps \in V$ such that
\begin{equation}\label{eq:refpbweak}
\forall v \in V, \quad B^\eps(u^\eps,v)=F(v)
\end{equation}
with
$$
B^\eps(u,v) = \intO \Aa^\eps \nab u \cdot \nab v, \qquad F(v)=\intO f \, v + \int_{\Gamma_N} g \, v.
$$
The bilinear form $B^\eps$ is symmetric, continuous and coercive on $V$. It thus defines an inner product and induces the energy norm $\vertiii{v} = \sqrt{B^\eps(v,v)}$ on $V$. The $\mH^1$-norm is denoted $\| \cdot \|_1$, and it is equivalent to the norm $\vertiii{\, \cdot \, }$ on $V$.

\medskip

The model problem~\eqref{eq:refpbweak} may be numerically solved using classical and direct numerical techniques, such as the finite element method (FEM). In this method, a conforming partition $\mT_H$ of $\Omega$ is introduced ($H$ denotes the characteristic size of the elements of $\mT_H$) and a finite dimensional space $V_H^0 \subset V$ is constructed from Lagrange basis functions $\phi_i^0$ associated with each node $i$ of $\mT_H$. The solution to~\eqref{eq:refpbweak} is thus approximated by $\dis u_H=\sum_i u_i \, \phi_i^0 \in V_H^0$, solution to
$$
\forall v \in V^0_H, \quad B^\eps(u_H,v)=F(v).
$$
A very fine mesh is necessary to represent the multiscale nature of the solution $u^\eps$. This leads to choosing $H < \eps$, and thus a prohibitively large number of degrees of freedom. The resulting method has a too high computational complexity.

As an alternative, multiscale model reduction techniques based on a scale separation assumption are typically used to solve~\eqref{eq:refpbweak} accurately and efficiently. On the one hand, homogenization strategies can be employed when the structure of $\Aa^\eps$ is periodic (or locally periodic, or random and stationary, \dots; see~\cite{BEN78}), or when one is only interested in macroscopic/upscaled properties of the model~\eqref{eq:refpb} (we refer e.g. to~\cite{E03,FEY03}). On the other hand, multiscale numerical approaches have the capability to approximate, efficiently and for general structures of $\Aa^\eps$, the details of $u^\eps$ at all scales. In the following, we focus on the Multiscale Finite Element Method (MsFEM)~\cite{HOU97,EFE09,LEB14}, which is one of these latter approaches. 

\subsection{MsFEM approximation}\label{section:MsFEMapprox}

\subsubsection{General strategy}

We define two partitions of $\Omega$: (i) a coarse grid partition $\mT_H$, referred to as the macro mesh; (ii) a fine grid partition $\mT_h$, subordinate to $\mT_H$, referred to as the micro mesh (the mesh size $h$ is chosen smaller, and often much smaller, than $H$). The mesh $\mT_h$ is in practice obtained by refining each coarse element of $\mT_H$. It should be chosen so that oscillations of the data can be accurately captured on $\mT_h$.

The main idea in the MsFEM approach is to construct a set $\{\phi_i^\eps\}$ of multiscale basis functions, associated with each node $i$ of the macro mesh $\mT_H$, which encode fine scale information and local heterogeneities of the solution space. These are constructed locally, in an offline phase, as the solutions, over each element $K$ of the coarse mesh (or over a local region $S_K \supset K$), to problems of the form
\begin{equation}\label{eq:localpbMsFEM}
\nab \cdot (\Aa^\eps \nab \phi_i^\eps)=0 \quad \text{in $K$}, 
\end{equation}
complemented with boundary conditions discussed below. We notice that the local problems~\eqref{eq:localpbMsFEM} are decoupled one from each other. They are thus well-suited to parallel computations. Furthermore, they do not depend on the loads $f$ and $g$, which is a main advantage in a multi-query context. These local problems are in practice solved numerically using the micro mesh $\mT_h$. 

\medskip

\begin{remark}
A variant in the construction of multiscale basis functions, referred to as the Generalized Multiscale Finite Element Method (GMsFEM), was recently introduced in~\cite{EFE13}. It consists in defining a larger number of basis functions over each element $K$ of $\mT_H$, which are obtained using local spectral problems.
\end{remark}

\medskip

Once the multiscale basis functions are computed, a classical Galerkin approximation of~\eqref{eq:refpbweak} is performed on the finite dimensional space $V^\eps_H = \text{Span}\{\phi_i^\eps\}$. The solution to~\eqref{eq:refpbweak} is thus approximated by some $\dis u^\eps_H=\sum_i u_i \, \phi_i^\eps \in V_H^\eps$, defined as the solution (computed during the online phase) to the following global problem:
\begin{equation}\label{eq:MsFEMpb}
\forall v \in V^\eps_H, \quad B^\eps(u^\eps_H,v) = F(v).
\end{equation}
We note that $V_H^\eps$ and $V_H^0$ share the same dimension. The above problem hence involves a limited number of degrees of freedom. The small scale information incorporated in the basis functions $\phi_i^\eps$ is thus brought to the large scales through couplings in the global stiffness matrix. The assembly of this matrix is inexpensive since it reuses local matrices computed and stored in the offline stage.

The MsFEM approach leads to a considerable decrease in memory requirements and CPU cost as an accurate approximation may be obtained for $H \gg \eps$. Nevertheless, \textit{a posteriori} error estimates and adaptive techniques are required to ensure a sufficiently high-fidelity approximation at a limited cost. These are addressed in Sections~\ref{section:CREconcept} and~\ref{section:goal-oriented}.

\subsubsection{Conforming MsFEM}

The choice of the boundary conditions for the local problems~\eqref{eq:localpbMsFEM} (when constructing the multiscale basis functions) is of critical importance, since it can significantly affect the accuracy of the MsFEM solution. A first choice, leading to a conforming numerical discretization, is to impose an affine evolution of $\phi_i^\eps$ along $\partial K$. Introducing the classical first-order and piecewise affine FE basis functions $\phi_i^0$, we thus consider the local problems
\begin{equation}\label{eq:localpbMsFEM2}
\nab \cdot (\Aa^\eps \nab \phi_i^\eps)=0 \quad \text{in $K$}, \qquad \phi_i^\eps=\phi_i^0 \quad \text{on $\partial K$}.
\end{equation}
In the regime of interest $H \ge \eps$, assuming that $\Aa^\eps$ has a periodic structure (i.e. $\Aa^\eps = \Aa^{\rm per}(\cdot/\eps)$ for some fixed matrix $\Aa^{\rm per}$) and that the local problems~\eqref{eq:localpbMsFEM2} are solved exactly (and under some technical regularity assumptions), the following \textit{a priori} error estimate holds~\cite{HOU99}:
\begin{equation}\label{eq:aprioriestimate1}
\|u^\eps-u^\eps_H\|_1 \le C\left(H+\sqrt{\frac{\eps}{H}}+\sqrt{\eps}\right),
\end{equation}
where the constant $C$ is independent of $\eps$ and $H$. This estimate is derived using the periodic homogenization theory that provides (when $\eps \ll 1$) a detailed description of $u^\eps$ and $\phi_i^\eps$ in terms of a two-scale expansion~\cite{BEN78}. It is the sum of three terms:
\begin{itemize}
\item the term in $O(\sqrt{\eps})$ comes from the homogenization theory, and quantifies the difference between $u^\eps$ and its two-scale expansion (this error is essentially located close to $\partial \Omega$; recall in particular that $u^\eps$ satisfies homogeneous Dirichlet boundary conditions on $\Gamma_D$, whereas its two-scale expansion does not);
\item the term in $O(\sqrt{\eps/H})$ again comes from the homogenization theory, and quantifies the difference between any basis function $\phi_i^\eps$ and its two-scale expansion (this error is essentially located close to $\partial K$);
\item the term in $O(H)$ is a classical \textit{a priori} error estimate contribution. It appears when approximating the homogenized solution (associated to $u^\eps$) by means of the piecewise affine basis functions $\phi^0_i$.
\end{itemize} 

\subsubsection{Oversampling technique}\label{section:oversampling}

The basis functions defined by~\eqref{eq:localpbMsFEM2} do not oscillate on the coarse mesh edges, which leads to a poor approximation of $u^\eps$ close to these edges. To improve on that, a possible strategy is to impose affine Dirichlet boundary conditions on the boundary of a domain $S_K \supset K$ which is slightly larger than $K$ (and often homothetic to $K$, for practical reasons), and to only use the interior information, on $K$, to construct the basis functions. The motivation stems from the fact that the information inside $K$ depends only mildly on the (incorrect) boundary conditions imposed on $\partial S_K$. This so-called oversampling technique thus consists in first solving the local problems
$$
\nab \cdot (\Aa^\eps \nab \psi_i^\eps)=0 \quad \text{in $S_K$}, \qquad \psi_i^\eps \, \text{is affine on $\partial S_K$}, \qquad \psi_i^\eps(\bs_j)=\delta_{ij},
$$
where $\bs_j$ are the coordinates of the vertices of $S_K$, before defining the MsFEM basis functions in $K$ as $\dis \phi_i^\eps=\psi_{i|K}^\eps$. In general, the obtained basis functions $\phi_i^\eps$ are discontinuous at the boundaries of the macro elements of $\mT_H$. This thus leads to a non-conforming MsFEM approximation: $V^\eps_H = \text{Span } \{ \phi_i^\eps \} \not\subset V$. In this framework, the bilinear form in~\eqref{eq:MsFEMpb} is replaced by $\dis B^\eps_H(u,v) = \sum_{K \in \mT_H}\int_K \Aa^\eps \nab u \cdot \nab v$, and the approximation $u^\eps_H \in V^\eps_H$ is defined as the solution to
$$
\forall v \in V^\eps_H, \quad B^\eps_H(u^\eps_H,v) = F(v).
$$

\medskip

\begin{remark}
It is also possible to define $\phi_i^\eps$ in $K$ as a linear combination of the functions $\psi^\eps_{j|K}$. We refer e.g. to~\cite[Eq.~(2.8)]{LEB14b} for more details.
\end{remark}

\medskip

An \textit{a priori} error estimate is given in~\cite{EFE00} for the oversampling variant of MsFEM, under the same assumptions as for~\eqref{eq:aprioriestimate1} (in particular, $\Aa^\eps$ is again assumed to be periodic). It reads
\begin{equation}\label{eq:aprioriestimate2}
\|u^\eps-u^\eps_H\|_{1,H} \le C\left(H+\frac{\eps}{H}+\sqrt{\eps}\right),
\end{equation}
where the constant $C$ is independent of $\eps$ and $H$ and where $\| \cdot \|_{1,H}$ is the broken $\mH^1$-norm. Furthermore, a Petrov-Galerkin formulation is proposed in~\cite{HOU04}, with the aim to further decrease the cell resonance error (which occurs when $H \approx \eps$). The formulation, in which the above non-conforming multiscale basis functions are used for trial functions while conforming linear Lagrange test functions are considered, reads
$$
\text{Find $u_H^\eps\in V_H^\eps$ such that, for any $v\in V^0_H$, \quad $B^\eps_H(u_H^\eps,v)=F(v)$.}
$$

\subsubsection{Higher-order MsFEM}\label{section:highorder}

In the original version of MsFEM, the multiscale basis functions $\phi_i^\eps$ are constructed from linear coarse-scale functions $\phi_i^0$. The construction of higher-order MsFEM basis functions was introduced in~\cite{ALL06} using local harmonic coordinates. For each element $K \in \mT_H$, these coordinates $w_j^{\eps,K}$ ($j=1,\dots,d$) are defined as the solutions to
\begin{equation}\label{eq:localpbho}
\nab \cdot(\Aa^\eps\nab w_j^{\eps,K}) = 0 \; \text{in $K$}, \qquad w_j^{\eps,K}=\bx \cdot \be_j \; \text{on $\partial K$}.
\end{equation}
We next define $w_j^\eps \in \mH^1(\Omega)$ by $w_{j|K}^\eps = w_j^{\eps,K}$, and set $\bw^\eps=(w_1^\eps,\dots,w_d^\eps)^T$.  Introducing the set $\{\phi_i^0\}$ of the Lagrange basis functions of order $k$, the higher-order MsFEM basis $\{\phi_i^\eps\}$ is constructed from the composition
$$
\phi_i^\eps = \phi_i^0 \circ \bw^\eps.
$$
Choosing Lagrange functions $\{\phi_i^0\}$ of order $1$ leads to the classical MsFEM basis defined by~\eqref{eq:localpbMsFEM2}. The approach is easy to implement since the computation of the oscillating functions $w_j^\eps$ is independent of the order $k$ of the coarse scale basis. This approach is advantageously used in the following sections to both compute \textit{a posteriori} error estimates and drive the adaptive strategy.

In the case when $\Aa^\eps$ is periodic, an \textit{a priori} error estimate is given in~\cite{ALL06} when considering MsFEM basis functions of order $k$. It reads
$$
\|u^\eps-u^\eps_H\|_1 \le C\left(H^k+\sqrt{\frac{\eps}{H}}+\sqrt{\eps}\right),
$$
where the constant $C$ is independent of $\eps$ and $H$. 

\medskip

\begin{remark}
As in the oversampling technique described in Section~\ref{section:oversampling}, the local problems~\eqref{eq:localpbho} can be posed on a larger domain $S_K \supset K$ before using the restriction to $K$ of the solutions. In this case, the support of $\phi_i^\eps$ may be different from that of $\phi_i^0$, and a non-conforming MsFEM approximation is obtained. 
\end{remark}

\medskip

\begin{remark}
An alternative higher-order MsFEM approach was recently proposed in~\cite{HES14}. 
\end{remark}

\section{The CRE concept applied to MsFEM}\label{section:CREconcept}

In this section, we review the method we have introduced in~\cite{CHA16b} to obtain a rigorous and fully computable upper bound on the MsFEM error measured in the energy norm, i.e. an upper bound on $\vertiii{u^\eps-u^\eps_H}$. It uses the concept of Constitutive Relation Error (CRE) based on duality~\cite{FRA72,ODE74}. A detailed overview on CRE (in a single scale framework) can be read in~\cite{LAD04}.

\subsection{Definition of an \textit{a posteriori} error estimate}\label{section:globalestimate}

Introduce the space of equilibrated fluxes
$$
\mS = \left \{\bp \in \mH(div,\Omega); \quad \forall v \in V, \quad \intO \bp \cdot \nab v = \intO fv + \int_{\Gamma_N} gv \right\}
$$
and define, for any pair $(\widehat{u},\widehat{\bp}) \in V \times \mS$, the CRE functional $E_{CRE}$ as
$$
E_{CRE}(\widehat{u},\widehat{\bp}) = \left( \intO (\Aa^\eps)^{-1}(\widehat{\bp}-\Aa^\eps \nab \widehat{u})\cdot (\widehat{\bp}-\Aa^\eps \nab \widehat{u}) \right)^{1/2} = \vertiii{\widehat{\bp}-\Aa^\eps \nab \widehat{u}}_{\mF},
$$
where $\vertiii{\, \cdot \, }_{\mF}$ denotes the energy norm on flux fields. The equality
\begin{equation}\label{eq:prager}
\forall (\widehat{u},\widehat{\bp}) \in V \times \mS, \qquad
\left( E_{CRE}(\widehat{u},\widehat{\bp}) \right)^2 = \vertiii{\bq^\eps-\widehat{\bp}}^2_{\mF} + \vertiii{u^\eps- \widehat{u}}^2,
\end{equation}
referred to as the Prager-Synge equality, is a direct consequence of the property $\dis \intO(\bq^\eps-\widehat{\bp})\cdot \nab(u^\eps- \widehat{u})=0$ (we recall that $\bq^\eps = \Aa^\eps \nab u^\eps$). Another equality, derived from~\eqref{eq:prager}, reads
\begin{equation}\label{eq:prager2}
E_{CRE}(\widehat{u},\widehat{\bp}) = 2 \vertiii{\bq^\eps-\widehat{\bp}^\ast}_{\mF} \quad \text{where} \quad \widehat{\bp}^\ast=\frac{1}{2}(\widehat{\bp} + \Aa^\eps \nab \widehat{u}).
\end{equation}
With a given $\widehat{\bp} \in \mS$ at hand, a guaranteed and computable bound on the global MsFEM error $\vertiii{u^\eps-u^\eps_H}$ can be derived from~\eqref{eq:prager}, as we now explain. When using a conforming MsFEM, we have $V^\eps_H \subset V$. We can hence choose $\widehat{u} = u^\eps_H$ in~\eqref{eq:prager}, which directly yields
\begin{equation} \label{eq:giens1}
\vertiii{u^\eps- u^\eps_H} \le E_{CRE}(u^\eps_H,\widehat{\bp}).
\end{equation}
When considering a nonconforming version of MsFEM, we may recover a field $\widehat{u}^\eps_H \in V$ from $u^\eps_H$ by a local averaging on the edges of the partition $\mT_H$. Using the triangle inequality, we get
\begin{equation} \label{eq:giens2}
\vertiii{u^\eps- u^\eps_H}_H \le E_{CRE}(\widehat{u}^\eps_H,\widehat{\bp}) + \vertiii{\widehat{u}^\eps_H- u^\eps_H}_H,
\end{equation}
where $\vertiii{\, \cdot \, }_H$ is the broken energy norm. The right-hand sides of~\eqref{eq:giens1} and~\eqref{eq:giens2} are fully computable.

Of course, the accuracy of the CRE estimate depends on the choice of $\widehat{\bp}$. In view of~\eqref{eq:prager}, we see that the difference $\vertiii{\bq^\eps-\widehat{\bp}}_{\mF}$ should be small in comparison with $\vertiii{u^\eps- \widehat{u}}$.

\subsection{Computation of an equilibrated flux field}\label{section:admissible_fields}

\begin{figure}[h]
\begin{center}
\includegraphics[width=40mm]{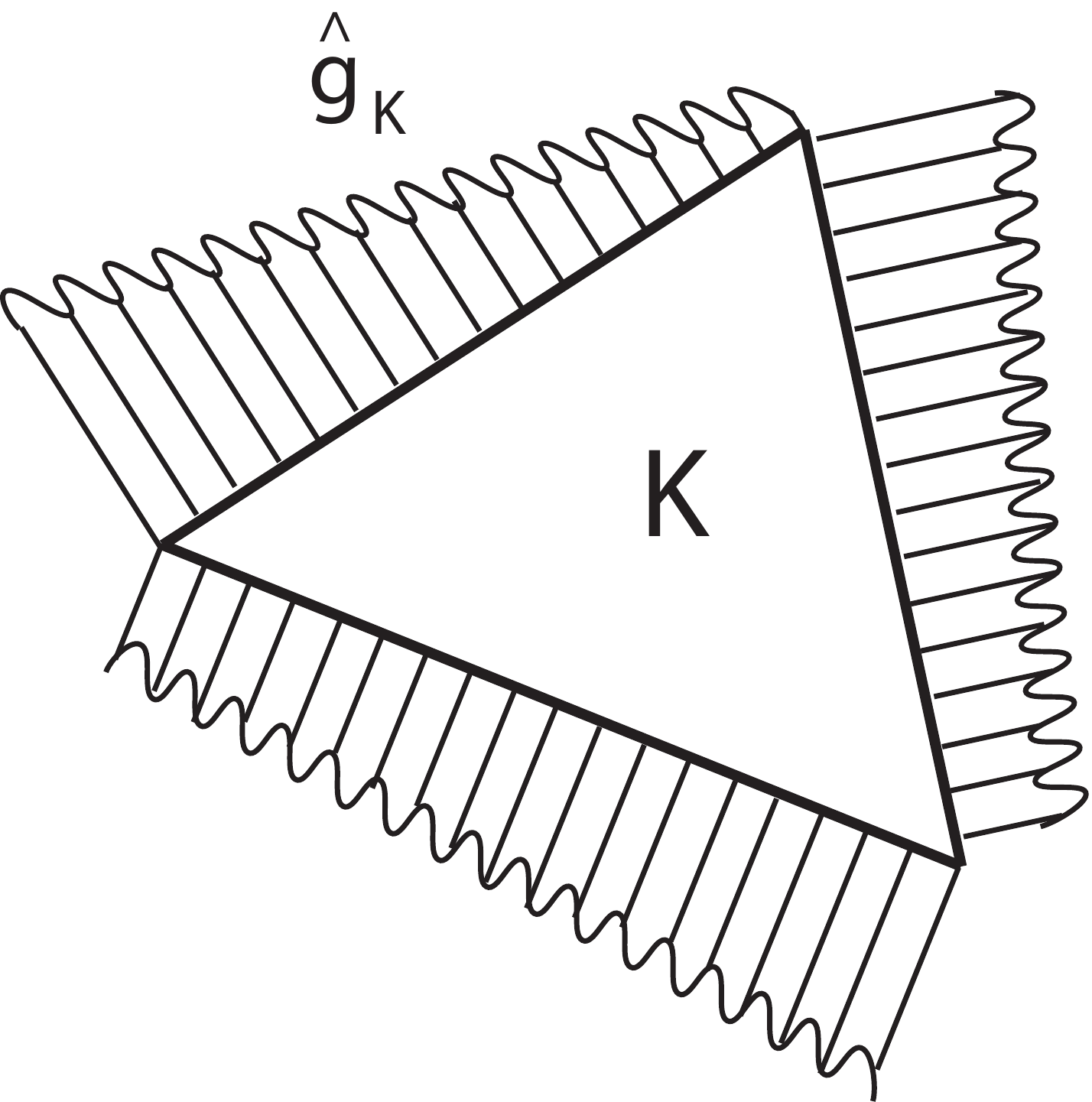}
\end{center}
\caption{Local problem~\eqref{eq:neumann} on the element $K$ with prescribed tractions $\widehat{g}_K$.}\label{fig:pblocalK}
\end{figure}

The technical point in the CRE concept is the construction of a relevant flux field $\widehat{\bp} \in \mS$. A convenient strategy, that we have introduced for all variants of MsFEM in~\cite{CHA16b} as an extension to that usually employed in the monoscale framework (see e.g.~\cite{LAD83,LAD96,PLE11}), consists in post-processing the approximate MsFEM flux $\bq^\eps_H=\Aa^\eps \nab u^\eps_H$ at hand in order to build a flux $\widehat{\bq}^\eps_H$ which belongs to $\mS$. The strategy consists in two steps:
\begin{itemize}
\item in Step 1, some tractions $\widehat{g}_K(\bx)$ are defined along the edges of each element $K \in \mT_H$, with $\widehat{g}_K=g$ on $\Gamma_N$ (see Fig.~\ref{fig:pblocalK}). They satisfy an equilibrium property at the element level,
\begin{equation}
\label{eq:equil}
\forall K, \quad \int_K f + \int_{\partial K}\widehat{g}_K = 0,
\end{equation}
as well as continuity between neighboring elements. Such tractions $\widehat{g}_K$ are constructed by solving local low-dimensional linear systems associated to each vertex in $\mT_H$, and using the partition of unity property verified by the test functions $\phi_i^\ast$ used in the Galerkin MsFEM formulation, i.e. $\dis \sum_i \phi_i^\ast=1$ ($\phi_i^\ast=\phi_i^\eps$ in Galerkin formulations, while $\phi_i^\ast=\phi_i^0$ in Petrov-Galerkin formulations). We refer to~\cite{CHA16b} for details. The tractions are represented as affine combinations of the test functions $\phi_i^\ast$:
$$
\widehat{g}_{K|\Gamma} = \frac{1}{2} \left( \bq^\eps_{H|K} + \bq^\eps_{H|K'} \right) \cdot \bn_K + \sum_i \alpha_i \, \phi_{i|\Gamma}^\ast,
$$
where $K$ and $K'$ are the elements such that $\Gamma = \overline{K} \cap \overline{K'}$;
\item in Step 2, the flux $\widehat{\bq}^\eps_H$ is constructed over each element $K$ as a solution to the local problem
$$
\forall v \in \mH^1(K), \qquad \int_K \widehat{\bq}^\eps_{H|K} \cdot \nab v = \int_K fv + \int_{\partial K} \widehat{g}_Kv.
$$
In practice, this problem is solved using a dual formulation. We hence search for some $w^\eps \in \mH^1(K)$ satisfying the Neumann problem
\begin{equation}
\label{eq:neumann}
\forall v \in \mH^1(K), \quad \int_K \Aa^\eps \nab w^\eps \cdot \nab v = \int_K fv + \int_{\partial K} \widehat{g}_Kv
\end{equation}
and next set $\widehat{\bq}^\eps_{H|K} = \Aa^\eps \nab w^\eps$. The Neumann problem~\eqref{eq:neumann} is well-posed (up to the addition of a constant) in view of~\eqref{eq:equil}. In practice, the problem~\eqref{eq:neumann} is numerically approximated by a Galerkin procedure on $\text{Span } \{ \phi_i^0 \circ \bw^\eps \}$ (see Section~\ref{section:highorder}), where $\phi_i^0$ are Lagrange function of high-order (in practice, of order 4). Since $\widehat{g}_{K|\Gamma}$ is a linear combination of the functions $\phi_{i|\Gamma}^\ast$ and $\Aa^\eps \nab \phi_{i|\Gamma}^\eps \cdot \bn_K$, we can use the superposition principle and obtain $w^\eps$ as a linear combination of pre-computed elementary solutions.
\end{itemize}

All technicalities in the construction of $\widehat{g}_K$ and $\widehat{\bq}^\eps_{H|K}$ are given in~\cite{CHA16b}. It is fruitful to notice that the overall strategy for constructing $\widehat{\bq}^\eps_H \in \mS$ is particularly suited to the MsFEM framework. Indeed, it is based on fine-scale computations at the macro element level only. These can be performed in the offline stage of MsFEM while constructing the multiscale basis functions.

\section{Goal-oriented error estimation}\label{section:goal-oriented}

\textit{A posteriori} error estimation on some outputs of the computation (which are typically needed for decision-making) is often very interesting for engineering applications. In many practical problems, one may indeed be interested in approximating some functional $Q(u^\eps)$, referred to as the quantity of interest, rather than the solution $u^\eps$ itself. In this section, we address goal-oriented error estimation in the MsFEM framework, again using the CRE concept introduced in Section~\ref{section:CREconcept}. 

\subsection{Adjoint problem and classical error indicator}\label{section:adjoint}

We consider a linear and continuous functional $Q: V \to \RR$, and assume that our quantity of interest is $Q(u^\eps)$. Typical examples include the average of the multiscale solution $u^\eps$, or of the exact flux $\bq^\eps = \Aa^\eps \nab u^\eps$, in some region of interest $\omega_Q \subset \Omega$. We then introduce the following adjoint problem~\cite{BEC01,GIL02}:
\begin{equation}\label{eq:adjointpbweak}
\text{Find $\widetilde{u}^\eps \in V$ such that, for any $v\in V$,}\quad B^{\eps\ast}(\widetilde{u}^\eps,v) = Q(v),
\end{equation}
where $B^{\eps\ast}(u,v)$ is the formal adjoint of the primal form $B^{\eps}(u,v)$, satisfying $B^{\eps\ast}(u,v)=B^\eps(v,u)$. In the current symmetric case, the forms $B^\eps$ and $B^{\eps\ast}$ are identical. We also define the adjoint flux $\widetilde{\bq}^\eps=\Aa^\eps\nab \widetilde{u}^\eps$ and the space of equilibrated fluxes (with regard to the adjoint problem):
$$
\widetilde{\mS} = \left \{\widetilde{\bp} \in \mH(div,\Omega); \quad \forall v \in V, \quad \intO \widetilde{\bp} \cdot \nab v = Q(v) \right\}.
$$

\medskip

\begin{remark}
In practice, the functional $Q$ is often defined by
\begin{equation} \label{eq:rep_Q}
Q(v) = \intO \left( \widetilde{f}_Q \, v + \widetilde{\bq}_Q \cdot \nab v \right) + \int_{\Gamma_N}\widetilde{g}_Q \, v,
\end{equation}
where $\left(\widetilde{f}_Q,\widetilde{\bq}_Q,\widetilde{g}_Q\right)$ is a set of extraction functions which can be mechanically interpreted in the adjoint problem as body, pre-flux and traction loadings, respectively. 
\end{remark}

\medskip

\begin{remark}
To handle some particular cases (such as quantities of interest related to the normal flux on $\Gamma_D$, which correspond to a functional $Q$ defined on a strict subset of $V$), a more general formulation of the adjoint problem can be found in~\cite{GIL02}.
\end{remark}

\medskip

\begin{remark}
Pointwise quantities of interest, such as $u^\eps(\bx_0)$ for some $\bx_0 \in \Omega$, cannot be directly considered in this framework since the solution $u^\eps \in \mH^1(\Omega)$ may not be continuous (for $d \ge 2$). To circumvent this issue, a regularization technique was introduced in~\cite{PRU99} for (monoscale) FEM. We do not consider such pointwise quantities of interest here.
\end{remark}

\medskip

The primal-dual equivalence follows for $u^\eps$ solution to the primal problem~\eqref{eq:refpbweak} and $\widetilde{u}^\eps$ solution to the adjoint problem~\eqref{eq:adjointpbweak}:
$$
Q(u^\eps) = B^{\eps\ast}(\widetilde{u}^\eps,u^\eps) = B^\eps(u^\eps,\widetilde{u}^\eps) = F(\widetilde{u}^\eps).
$$
Let $u^\eps_H$ be the MsFEM approximation of $u^\eps$, and let $\widehat{u}^\eps_H$ be a function in $V$ close to $u^\eps_H$ (if we work with a conforming MsFEM approach, then $\widehat{u}^\eps_H = u^\eps_H$). From the adjoint solution $\widetilde{u}^\eps$, an exact representation of the MsFEM error on the quantity $Q(u^\eps)$ is available:
\begin{equation}\label{eq:represQ}
Q(u^\eps)-Q(u^\eps_H) = Q(u^\eps-\widehat{u}^\eps_H) + \Delta Q = B^\eps(u^\eps-\widehat{u}^\eps_H,\widetilde{u}^\eps) + \Delta Q = R^\eps(\widetilde{u}^\eps) + \Delta Q,
\end{equation}
with $\Delta Q=Q(\widehat{u}^\eps_H-u^\eps_H)$, and where $R^\eps(v)=F(v)-B^\eps(\widehat{u}^\eps_H,v)$ is the residual functional associated with $\widehat{u}^\eps_H \in V$. The representation~\eqref{eq:represQ} shows that the adjoint solution $\widetilde{u}^\eps$, often termed generalized Green's function, provides for the sensitivity of the error on $Q$ to the local sources of the MsFEM error in the whole domain $\Omega$.

A convenient error indicator may be derived from~\eqref{eq:represQ} using the dual weighted residual (DWR) method proposed in~\cite{BEC96} for FEM and extended to multiscale computations in recent works (see for instance~\cite{ABD13a,CHU16}). The DWR method is based on hierarchical \textit{a posteriori} error analysis. We compute an approximate adjoint solution $\widetilde{u}^\eps_+$ from a hierarchically refined MsFEM space $V^\eps_+ \supset V^\eps_H$ (higher-order basis functions, refined mesh) and we assume that the error between $\widetilde{u}^\eps$ and $\widetilde{u}^\eps_+$ is small. We then write
\begin{equation}\label{eq:DWRindicator}
Q(u^\eps)-Q(u^\eps_H) - \Delta Q = R^\eps(\widetilde{u}^\eps_+) + R^\eps(\widetilde{u}^\eps-\widetilde{u}^\eps_+) \approx R^\eps(\widetilde{u}^\eps_+),
\end{equation}
where the right-hand side and $\Delta Q$ are fully computable.

\medskip

\begin{remark}
Due to the Galerkin orthogonality, choosing the same MsFEM space $V^\eps_H$ to compute $u^\eps_H$ and $\widetilde{u}^\eps_+$ would lead to a meaningless indicator (in that case, we would have $R^\eps(\widetilde{u}^\eps_+) = 0$). Consequently, a convenient approximation of the adjoint solution should involve a subspace in $(V^\eps_H)^\perp$.
\end{remark}

\medskip

The above argument can be somewhat quantified. Let $u^\eps_+ \in V^\eps_+$ be the approximation (in the refined space) of the primal solution. We assume that the saturation assumption $|Q(u^\eps)-Q(u^\eps_+)| \le \beta |Q(u^\eps)-Q(\widehat{u}^\eps_H)|$ is satisfied for some $0 \le \beta < 1$. Using the fact that $u^\eps_+-\widehat{u}^\eps_H \in V^\eps_+$ and next the fact that $\widetilde{u}^\eps_+ \in V^\eps_+$, we obtain
$$
Q(u^\eps_+)-Q(\widehat{u}^\eps_H)=B^\eps(u^\eps_+-\widehat{u}^\eps_H,\widetilde{u}^\eps_+)=B^\eps(u^\eps-\widehat{u}^\eps_H,\widetilde{u}^\eps_+)=R^\eps(\widetilde{u}^\eps_+).
$$
Using the triangle inequality, we get
$$
\frac{|R^\eps(\widetilde{u}^\eps_+)|}{1+\beta} \le |Q(u^\eps)-Q(u^\eps_H)-\Delta Q| \le \frac{|R^\eps(\widetilde{u}^\eps_+)|}{1-\beta},
$$
which quantifies the approximation~\eqref{eq:DWRindicator}. The indicator defined in~\eqref{eq:DWRindicator} may be used to drive adaptive algorithms, but it does not provide for {\em guaranteed} error bounds on $Q$. Alternative techniques are proposed in the next sections.

\subsection{Guaranteed error bounds using CRE}\label{section:bounding}

In the FEM context, classical approaches to bound the error on a quantity of interest $Q$ consist in solving the primal and adjoint problems with the same approximation space, before using \textit{a posteriori} bounds available for the error in the energy norm~\cite{AIN00}. These approaches can be easily extended to MsFEM with general, conforming or nonconforming, formulations.

Denoting by $\widetilde{u}^\eps_H$ the approximate adjoint solution obtained using the space $V^\eps_H$, and potentially post-processing it to define $\widehat{\widetilde{u}}^\eps_H\in V$ in case $V^\eps_H \not\subset V$ (see Section~\ref{section:globalestimate}), the equality~\eqref{eq:represQ} yields that
\begin{equation}\label{eq:ineqgo1_part1}
Q(u^\eps)-Q(u^\eps_H)-\Delta Q-R^\eps\left(\widehat{\widetilde{u}}^\eps_H\right) = B^\eps \left(u^\eps-\widehat{u}^\eps_H,\widetilde{u}^\eps-\widehat{\widetilde{u}}^\eps_H \right).
\end{equation}
The Cauchy-Schwarz inequality then yields
\begin{equation}\label{eq:ineqgo1_part2}
\left| Q(u^\eps)-Q(u^\eps_H)-\Delta Q-R^\eps\left(\widehat{\widetilde{u}}^\eps_H\right) \right| \le \vertiii{u^\eps-\widehat{u}^\eps_H} \ \vertiii{\widetilde{u}^\eps-\widehat{\widetilde{u}}^\eps_H}.
\end{equation}

\medskip

\begin{remark}
In the conforming formulation of MsFEM (i.e. when $V^\eps_H \subset V$), we can take $\widehat{u}^\eps_H=u^\eps_H$ and $\widehat{\widetilde{u}}^\eps_H=\widetilde{u}^\eps_H$ so that the Galerkin orthogonality $\dis R^\eps\left(\widehat{\widetilde{u}}^\eps_H\right)=0$ holds. The result~\eqref{eq:ineqgo1_part2} can be recast as
$$
|Q(u^\eps)-Q(u^\eps_H)| \le \vertiii{u^\eps-u^\eps_H} \ \vertiii{\widetilde{u}^\eps-\widetilde{u}^\eps_H}.
$$
This is a classical result in the monoscale conforming FEM context.
\end{remark}

\medskip

Consequently, guaranteed error bounds on $Q$ can be obtained from any guaranteed \textit{a posteriori} global error estimate on $\vertiii{u^\eps-\widehat{u}^\eps_H}$ and $\vertiii{\widetilde{u}^\eps-\widehat{\widetilde{u}}^\eps_H}$. Using the CRE concept introduced in Section~\ref{section:CREconcept}, \eqref{eq:ineqgo1_part2} yields the \textit{a posteriori} estimate
\begin{equation}\label{eq:CREgobound1}
\left| Q(u^\eps)-Q(u^\eps_H)-\Delta Q-R^\eps\left(\widehat{\widetilde{u}}^\eps_H\right) \right| \le E_{CRE}(\widehat{u}^\eps_H,\widehat{\bq}^\eps_H) \, E_{CRE}\left(\widehat{\widetilde{u}}^\eps_H,\widehat{\widetilde{\bq}}^\eps_H\right),
\end{equation}
where $\widehat{\bq}^\eps_H \in \mS$, resp. $\widehat{\widetilde{\bq}}^\eps_H \in\widetilde{\mS}$, is constructed as a post-processing of $\bq^\eps_H$, resp. $\widetilde{\bq}^\eps_H$, as detailed in Section~\ref{section:admissible_fields}. However, the estimate~\eqref{eq:CREgobound1}, based on the Cauchy-Schwarz inequality over the whole domain $\Omega$, does not exploit the locality of $Q$. It thus leads to crude upper bounds when MsFEM errors for primal and adjoint problems have disjoint supports. A better estimate would consist in using the Cauchy-Schwarz inequality at the element level, under the form
$$
\left| Q(u^\eps)-Q(u^\eps_H)-\Delta Q-R^\eps \left(\widehat{\widetilde{u}}^\eps_H \right) \right| \le \sum_K \vertiii{u^\eps-\widehat{u}^\eps_H}_{|K} \ \vertiii{\widetilde{u}^\eps-\widehat{\widetilde{u}}^\eps_H}_{|K}.
$$
However, computable upper bounds on the error in the energy norm restricted to each element $K$ are hardly available, even though recent works such as~\cite{LAD13} have proposed first ideas in that direction.

\medskip

An alternative and more accurate estimate than~\eqref{eq:CREgobound1} can be constructed by taking advantage of the property~\eqref{eq:prager2} of the CRE concept. Indeed, we infer from~\eqref{eq:ineqgo1_part1} that
\begin{align}
Q(u^\eps)-Q(u^\eps_H)-\Delta Q-R^\eps\left(\widehat{\widetilde{u}}^\eps_H\right) 
&=
B^\eps\left(u^\eps-\widehat{u}^\eps_H,\widetilde{u}^\eps-\widehat{\widetilde{u}}^\eps_H\right)
\nonumber
\\
&= \intO \nab(u^\eps-\widehat{u}^\eps_H) \cdot \left(\widetilde{\bq}^\eps-\Aa^\eps \nab \widehat{\widetilde{u}}^\eps_H \right)
\nonumber
\\
&= \intO \nab(u^\eps-\widehat{u}^\eps_H) \cdot \left(\widehat{\widetilde{\bq}}_H^\eps-\Aa^\eps \nab \widehat{\widetilde{u}}^\eps_H \right)
\nonumber
\\
&= \intO (\Aa^\eps)^{-1} \left( \bq^\eps-\widehat{\bq}^{\eps\ast}_H \right) \cdot \left(\widehat{\widetilde{\bq}}_H^\eps-\Aa^\eps \nab \widehat{\widetilde{u}}^\eps_H \right) + C^\eps_H,
\label{eq:relationsgo}
\end{align}
with $\widehat{\bq}^{\eps\ast}_H = (\widehat{\bq}^\eps_H + \Aa^\eps \nab \widehat{u}^\eps_H)/2$ and $C^\eps_H$ a computable term defined by
$$
C^\eps_H =  \intO (\Aa^\eps)^{-1} \left( \widehat{\bq}^{\eps\ast}_H - \Aa^\eps \nab \widehat{u}^\eps_H \right) \cdot \left( \widehat{\widetilde{\bq}}_H^\eps-\Aa^\eps \nab \widehat{\widetilde{u}}^\eps_H \right) = \frac{1}{2} \intO (\Aa^\eps)^{-1} \left( \widehat{\bq}^\eps_H - \Aa^\eps \nab \widehat{u}^\eps_H \right) \cdot \left(\widehat{\widetilde{\bq}}_H^\eps-\Aa^\eps \nab \widehat{\widetilde{u}}^\eps_H\right).
$$
We next get from~\eqref{eq:relationsgo} that
$$
Q(u^\eps)-Q(u^\eps_H)-\Delta Q-\overline{C}^\eps_H = \intO (\Aa^\eps)^{-1} \left(\bq^\eps-\widehat{\bq}^{\eps\ast}_H \right) \cdot \left(\widehat{\widetilde{\bq}}_H^\eps-\Aa^\eps \nab \widehat{\widetilde{u}}^\eps_H \right)
$$
with
$$
\overline{C}^\eps_H = C^\eps_H + R^\eps\left(\widehat{\widetilde{u}}^\eps_H\right) = \frac{1}{2}\intO (\Aa^\eps)^{-1} \left(\widehat{\bq}^\eps_H - \Aa^\eps \nab \widehat{u}^\eps_H \right) \cdot \left(\widehat{\widetilde{\bq}}_H^\eps+\Aa^\eps \nab \widehat{\widetilde{u}}^\eps_H \right).
$$
Using the Cauchy-Schwarz inequality and next~\eqref{eq:prager2}, we obtain that
\begin{equation}\label{eq:CREgobound2}
\left| Q(u^\eps)-Q(u^\eps_H)-\Delta Q-\overline{C}^\eps_H \right|
\le
\vertiii{\bq^\eps-\widehat{\bq}^{\eps\ast}_H}_{\mF} \ \vertiii{\widehat{\widetilde{\bq}}_H^\eps-\Aa^\eps \nab \widehat{\widetilde{u}}^\eps_H}_{\mF}
=
\frac{1}{2} E_{CRE}\left(\widehat{u}^\eps_H,\widehat{\bq}^\eps_H\right) \, E_{CRE}\left(\widehat{\widetilde{u}}^\eps_H,\widehat{\widetilde{\bq}}^\eps_H\right).
\end{equation}
Note that, thanks to the factor $1/2$ on the right-hand side, the estimate~\eqref{eq:CREgobound2} yields a more accurate estimation of $Q(u^\eps)-Q(u^\eps_H)$ than~\eqref{eq:CREgobound1}.

The fully computable quantity $Q(u^\eps_H)+\Delta Q+\overline{C}^\eps_H$ can be interpreted as a corrected approximation of $Q(u^\eps)$. Furthermore, the estimate~\eqref{eq:CREgobound2} can also be written as
$$
Q(u^\eps)-Q(u^\eps_H) = \Delta Q+\overline{C}^\eps_H + \zeta \quad \text{for some $\zeta \in \left[ -\frac{1}{2} E_{CRE} \, \widetilde{E}_{CRE}, \frac{1}{2}E_{CRE} \, \widetilde{E}_{CRE} \right]$},
$$
hence
\begin{equation}\label{eq:goestimate}
|Q(u^\eps)-Q(u^\eps_H)| \le \eta^Q := \max_{\theta=\pm 1} \left| \Delta Q + \overline{C}^\eps_H + \frac{\theta}{2} E_{CRE} \, \widetilde{E}_{CRE} \right|.
\end{equation}
In contrast to~\eqref{eq:CREgobound1}, the quantity $\eta^Q$ partially takes into account error cancellations.

\medskip

It is fruitful to notice that the bounds~\eqref{eq:CREgobound1} and~\eqref{eq:CREgobound2}, constructed from the CRE concept, hold for independent (MsFEM) discretizations of the primal and adjoint problems as they are merely based on equilibrium properties. In the following section, we show how a convenient local enrichment of the adjoint approximate solution can thus be advantageously used to get sharp error bounds from~\eqref{eq:CREgobound2}. 

\subsection{Non-intrusive local enrichment of the adjoint solution}\label{section:nonintrusive}

We extend to the MsFEM context an enrichment technique initially introduced in~\cite{GRA03,CHA08,LAD10} for goal-oriented error estimation in FEM models. The general idea is to improve the accuracy of the approximate adjoint solution $\widetilde{u}^\eps_H$ (or $\widehat{\widetilde{u}}^\eps_H$) by adding \textit{a priori} information on the problem, i.e. given functions whose purpose is to reproduce localized high-gradient components of $\widetilde{u}^\eps$ (which are difficult to numerically capture) in the vicinity of the adjoint loading region. These functions, referred to as \textit{handbook functions} and denoted $\widetilde{u}^\eps_{hand}$ (with associated flux $\widetilde{\bq}^\eps_{hand}=\Aa^\eps \nab \widetilde{u}^\eps_{hand}$), correspond to generalized Green's functions. In analogy to~\eqref{eq:adjointpbweak}, they satisfy
\begin{equation} \label{eq:def_u_hand}
\text{Find $\widetilde{u}^\eps_{hand} \in \mH^1_0(\Omega_{hand})$ such that, for any $v\in \mH^1_0(\Omega_{hand})$,}\quad B^{\eps\ast}(\widetilde{u}^\eps_{hand},v) = Q_{hand}(v),
\end{equation}
where $Q_{hand}$ is an auxiliary quantity of interest (the eventual quantity of interest $Q$ in~\eqref{eq:adjointpbweak} may be equal to $Q_{hand}$, or a linear combination of various $Q_{hand}$, or various translations of $Q_{hand}$, \ldots). The domain $\Omega_{hand} \subset \Omega$ is chosen to be sufficiently large, so that the solution to the above problem with $\mH^1_0(\Omega_{hand})$ replaced by $V$ would indeed be close to zero on $\partial \Omega_{hand}$. In practice, $\widetilde{u}^\eps_{hand}$ is numerically pre-computed, using an accurate fine scale discretization. Examples of such enrichment functions, for a 2D multiscale medium described by $\Aa^\eps = A^\eps \, \II_2$ with $A^\eps(x_1,x_2)=3+\cos(2\pi x_1/\eps)+\cos(2\pi x_2/\eps)$ and $\eps=1/20$ (see Fig.~\ref{fig:mediumhandbook}), are given below. 

\begin{figure}[h]
\begin{center}
\includegraphics[width=75mm]{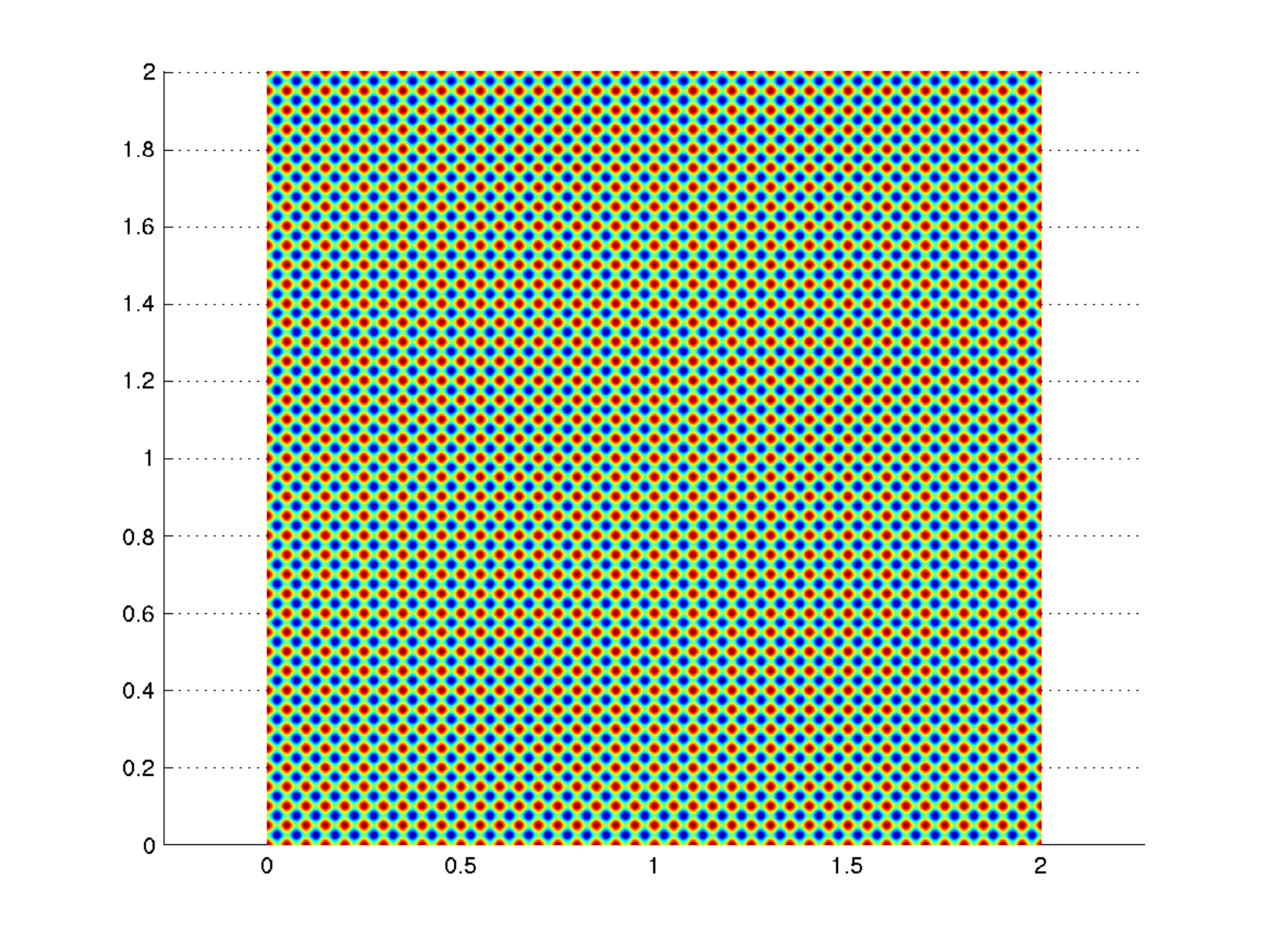}
\end{center}
\caption{Evolution of the material parameter.}
\label{fig:mediumhandbook}
\end{figure}

\begin{figure}[h]
\begin{center}
\includegraphics[width=75mm]{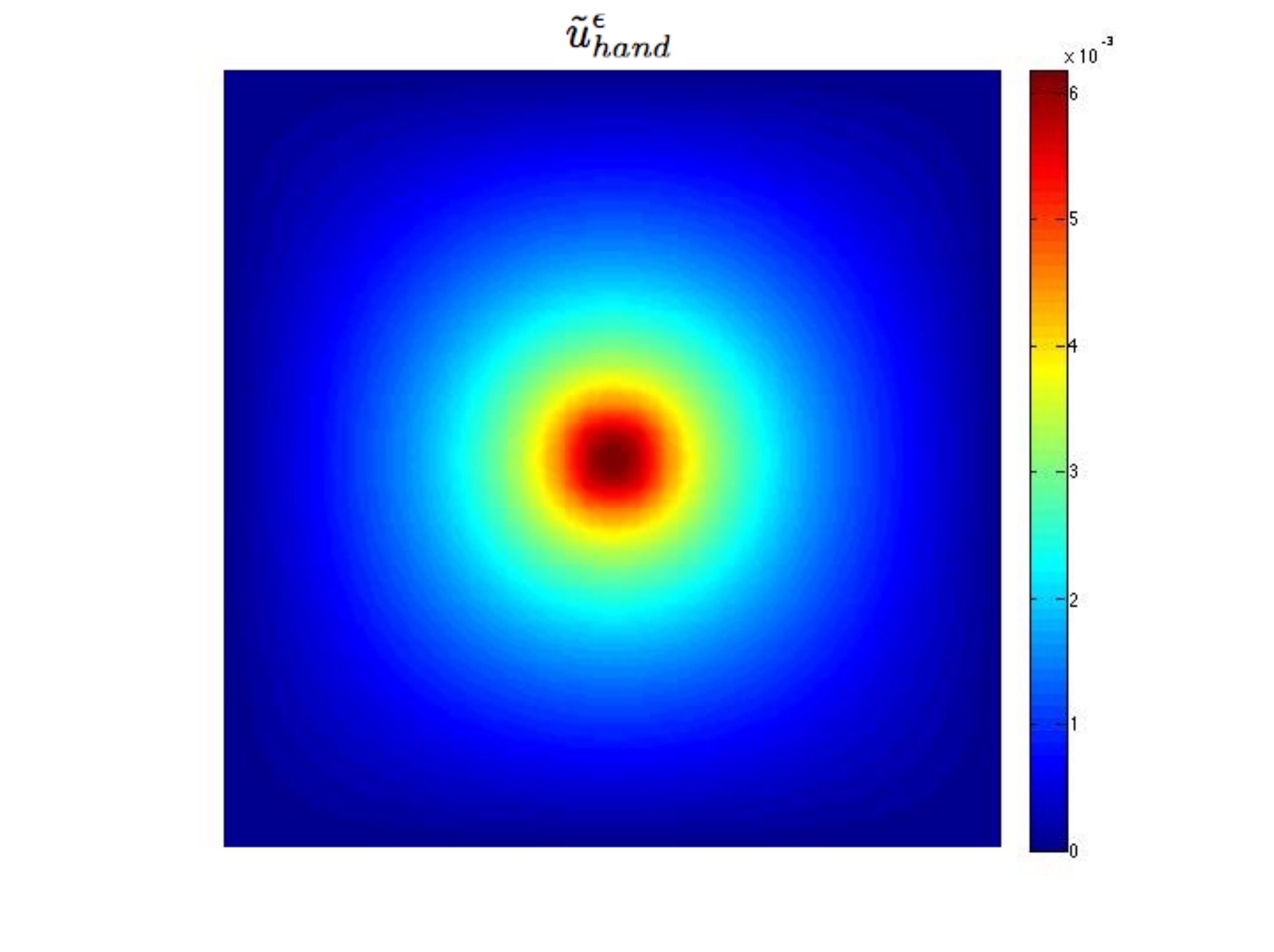}
\includegraphics[width=75mm]{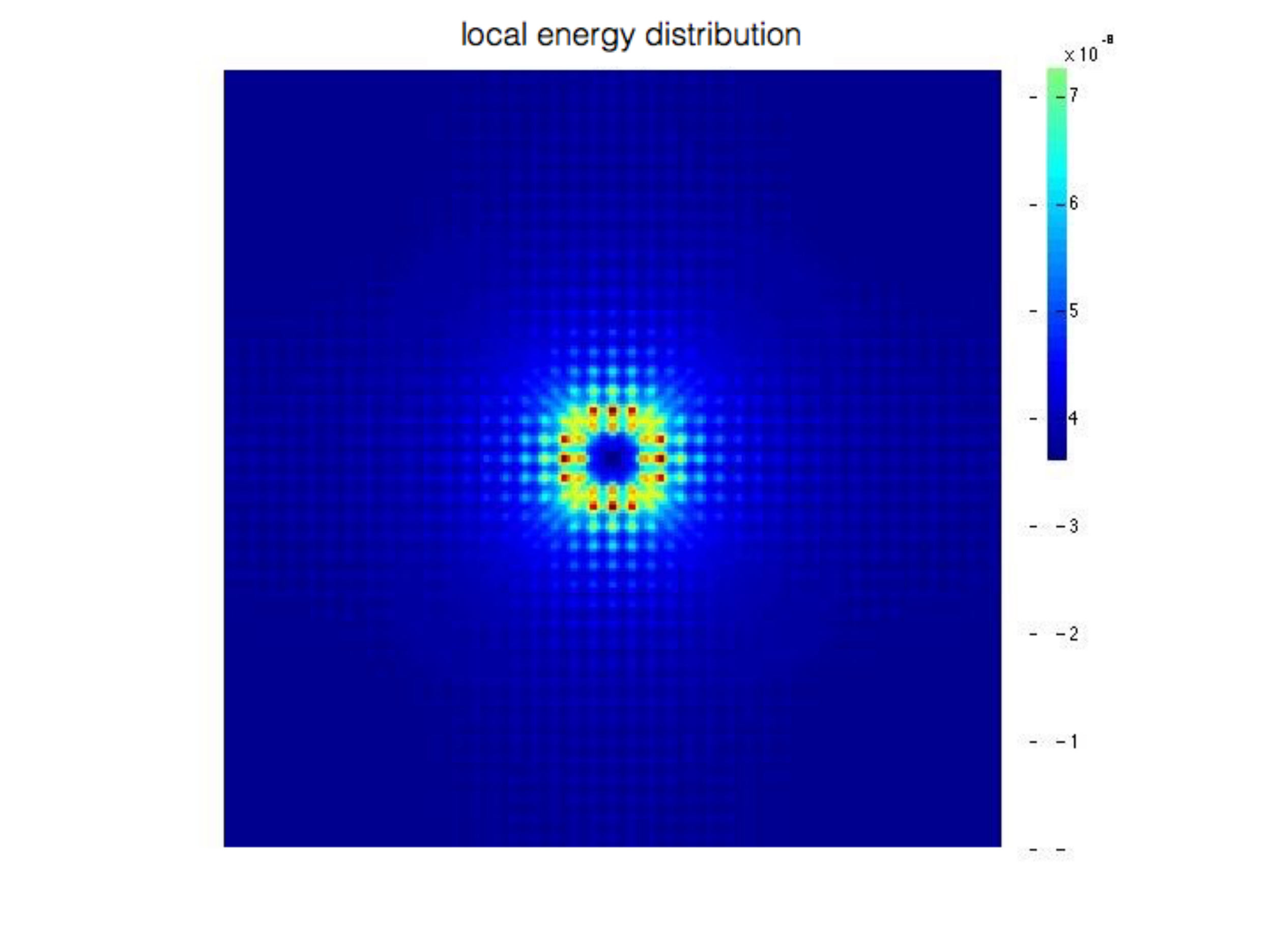}\\
\includegraphics[width=75mm]{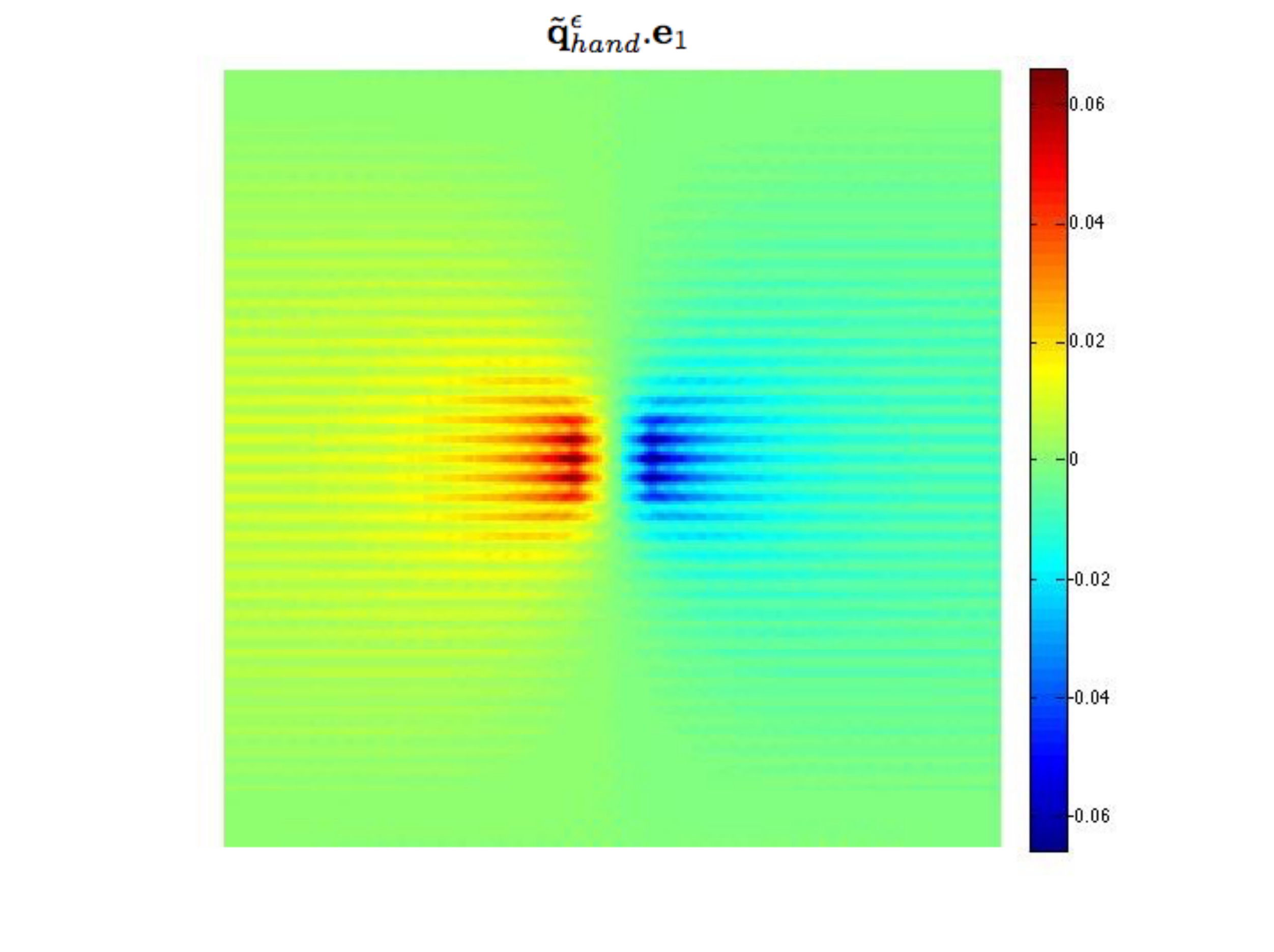}
\includegraphics[width=75mm]{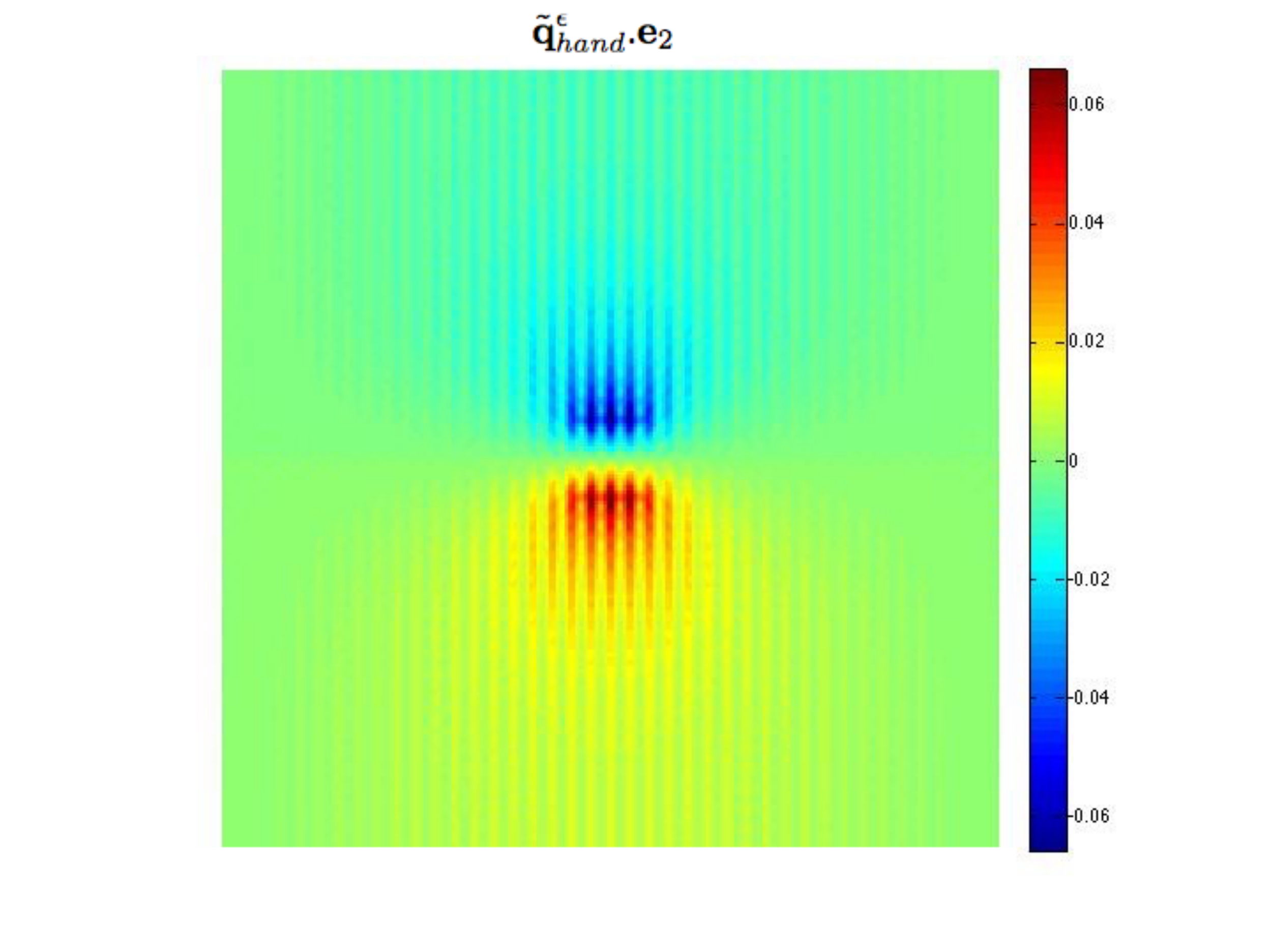}
\end{center}
\caption{Handbook function in a sufficiently large domain subjected to a constant body force in a local square region $\omega_Q$: $\widetilde{u}^\eps_{hand}$ (top left), energy distribution (top right), $\widetilde{\bq}^\eps_{hand} \cdot \be_1$ (bottom left) and $\widetilde{\bq}^\eps_{hand} \cdot \be_2$ (bottom right).}
\label{fig:handbook1}
\end{figure}

A first example corresponds to a linear functional $Q$ given by~\eqref{eq:rep_Q} with $\widetilde{\bq}_Q = 0$, $\widetilde{g}_Q = 0$ and a localized body force $\widetilde{f}_Q = 0$ except in the square domain $\omega_Q$, centered at the center of $\Omega$ and of size $4 \eps \times 4 \eps$, where $\widetilde{f}_Q = 1$. Such a quantity of interest corresponds to the average of $u^\eps$ in $\omega_Q$. The corresponding handbook function is shown in Fig.~\ref{fig:handbook1}.

\medskip

\begin{remark} \label{rem:pas_de_gQ}
Here and in all what follows, the quantity of interest is given by~\eqref{eq:rep_Q} for some $\widetilde{f}_Q$ and $\widetilde{\bq}_Q$, with $\widetilde{g}_Q = 0$. We take for $\Omega_{hand}$ a domain 10 times as large (in each direction) as the support of $\widetilde{f}_Q$ and $\widetilde{\bq}_Q$. For the example discussed above, we thus choose for $\Omega_{hand}$ a square domain centered at the center of $\Omega$ and of size $40 \eps \times 40 \eps$.
\end{remark}

\medskip

\begin{figure}[h]
\begin{center}
\includegraphics[width=75mm]{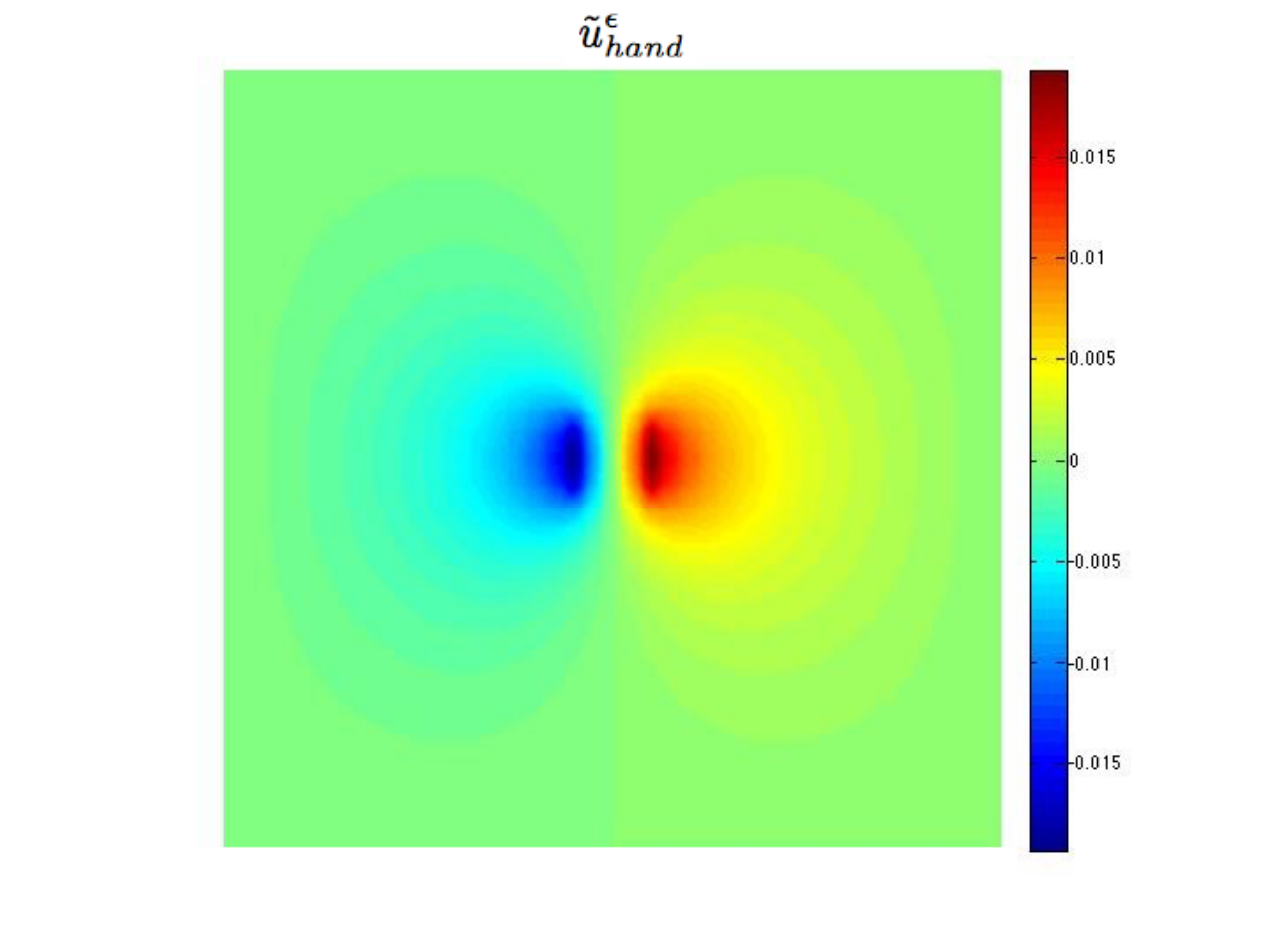}
\includegraphics[width=75mm]{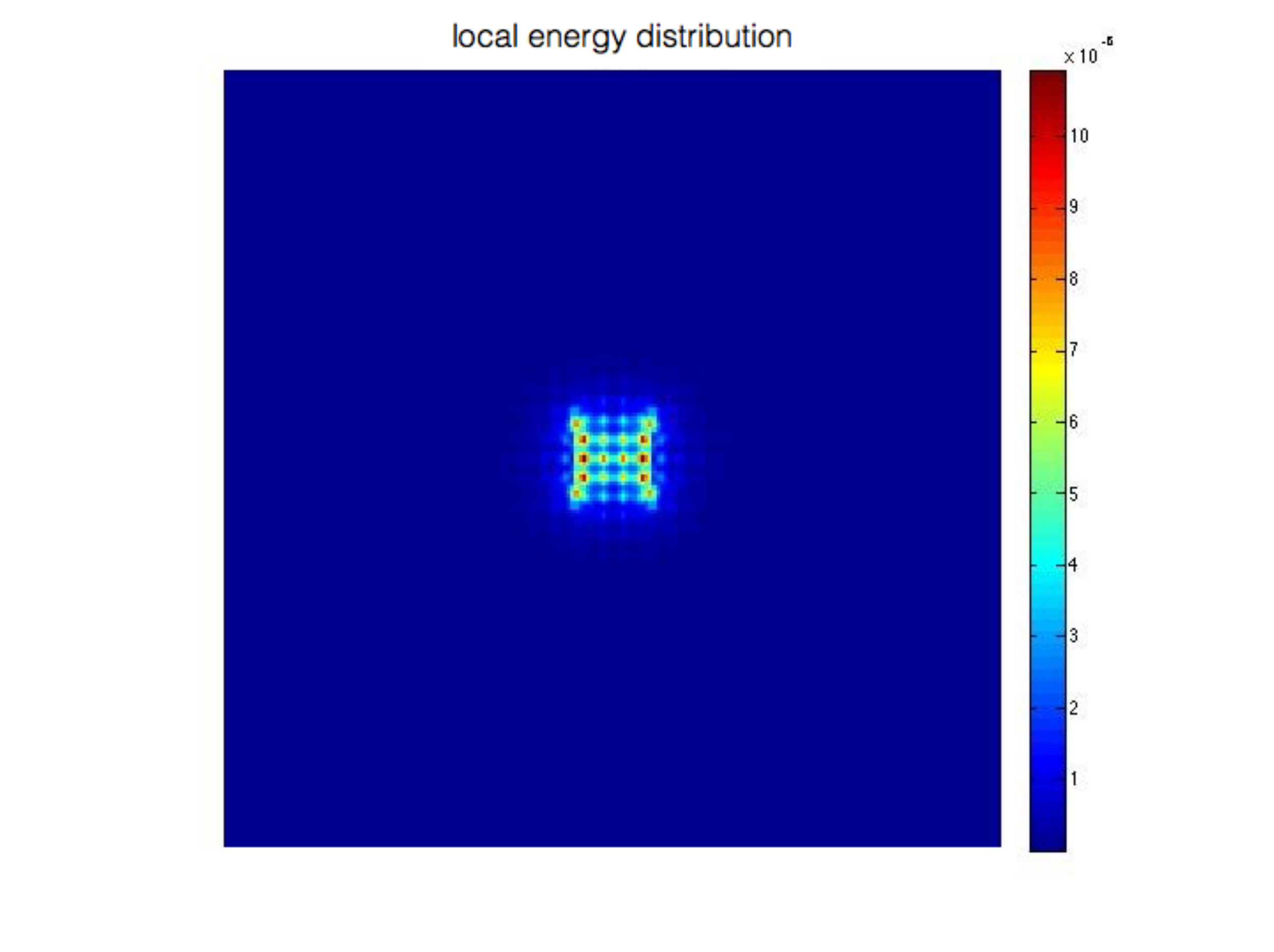}\\
\includegraphics[width=75mm]{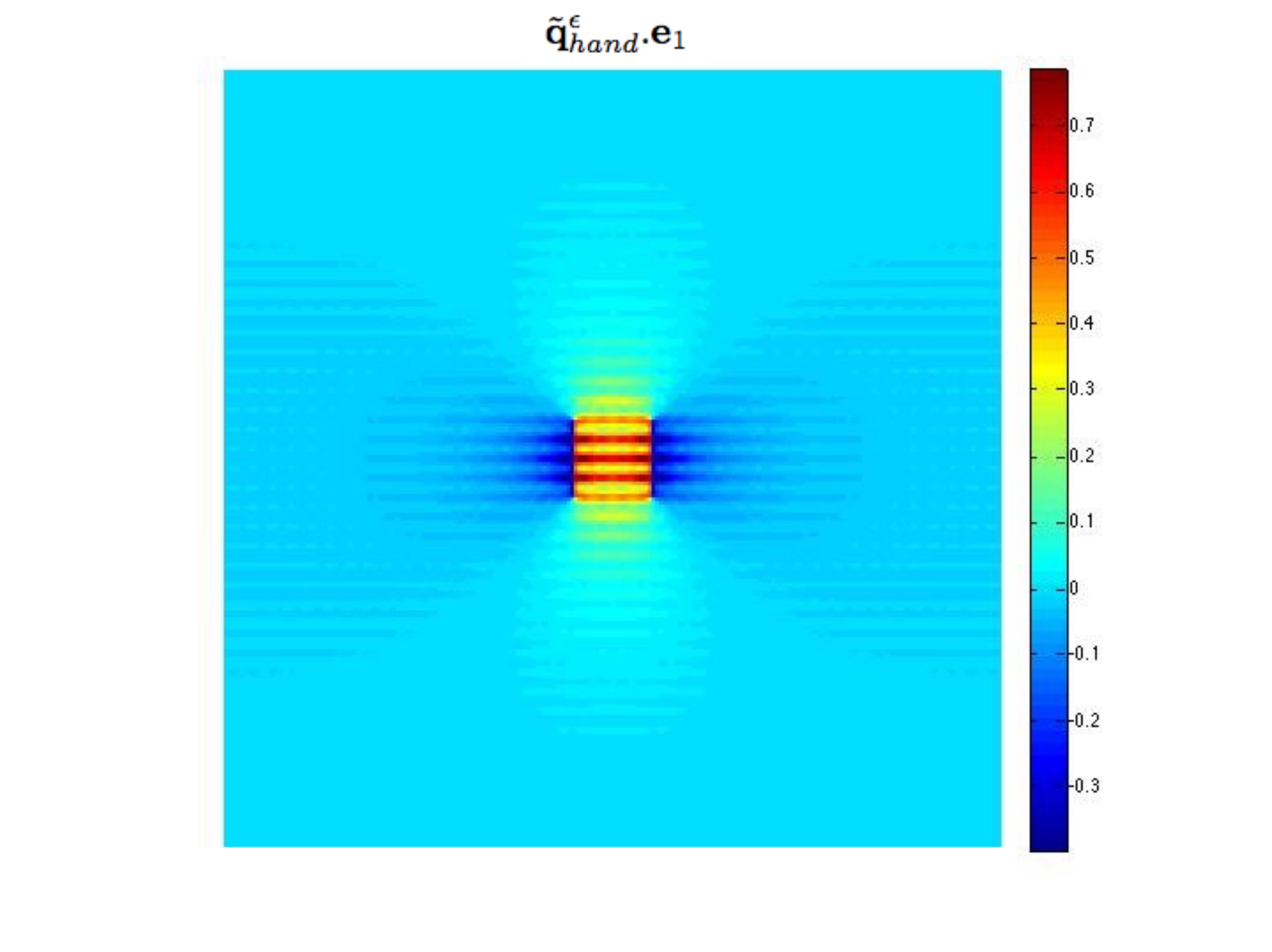}
\includegraphics[width=75mm]{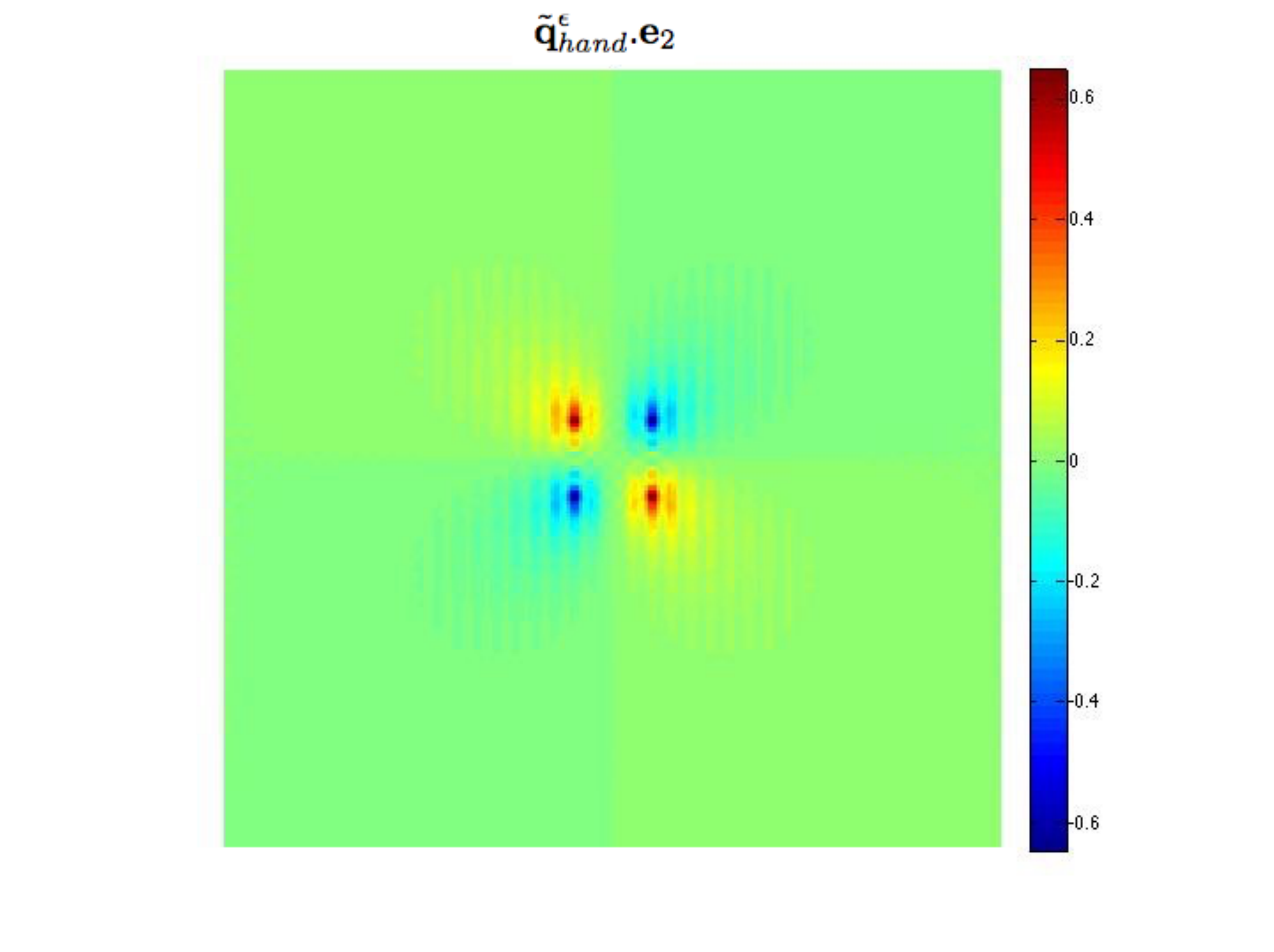}
\end{center}
\caption{Handbook function in a sufficiently large domain subjected to a constant pre-flux in a local square region $\omega_Q$: $\widetilde{u}^\eps_{hand}$ (top left), energy distribution (top right), $\widetilde{\bq}^\eps_{hand} \cdot \be_1$ (bottom left) and $\widetilde{\bq}^\eps_{hand} \cdot \be_2$ (bottom right).}
\label{fig:handbook2}
\end{figure}

A second example corresponds to a linear functional $Q$ given by~\eqref{eq:rep_Q} with $\widetilde{f}_Q = 0$, $\widetilde{g}_Q = 0$ and a localized pre-flux $\widetilde{\bq}_Q = (1,0)^T$ in the same domain $\omega_Q$ as above, and $\widetilde{\bq}_Q = 0$ elsewhere in $\Omega$. Such a quantity of interest corresponds to the average of $\be_1 \cdot \nabla u^\eps$ in $\omega_Q$. The corresponding handbook function is shown in Fig.~\ref{fig:handbook2}. 

\begin{figure}[h]
\begin{center}
\includegraphics[height=40mm]{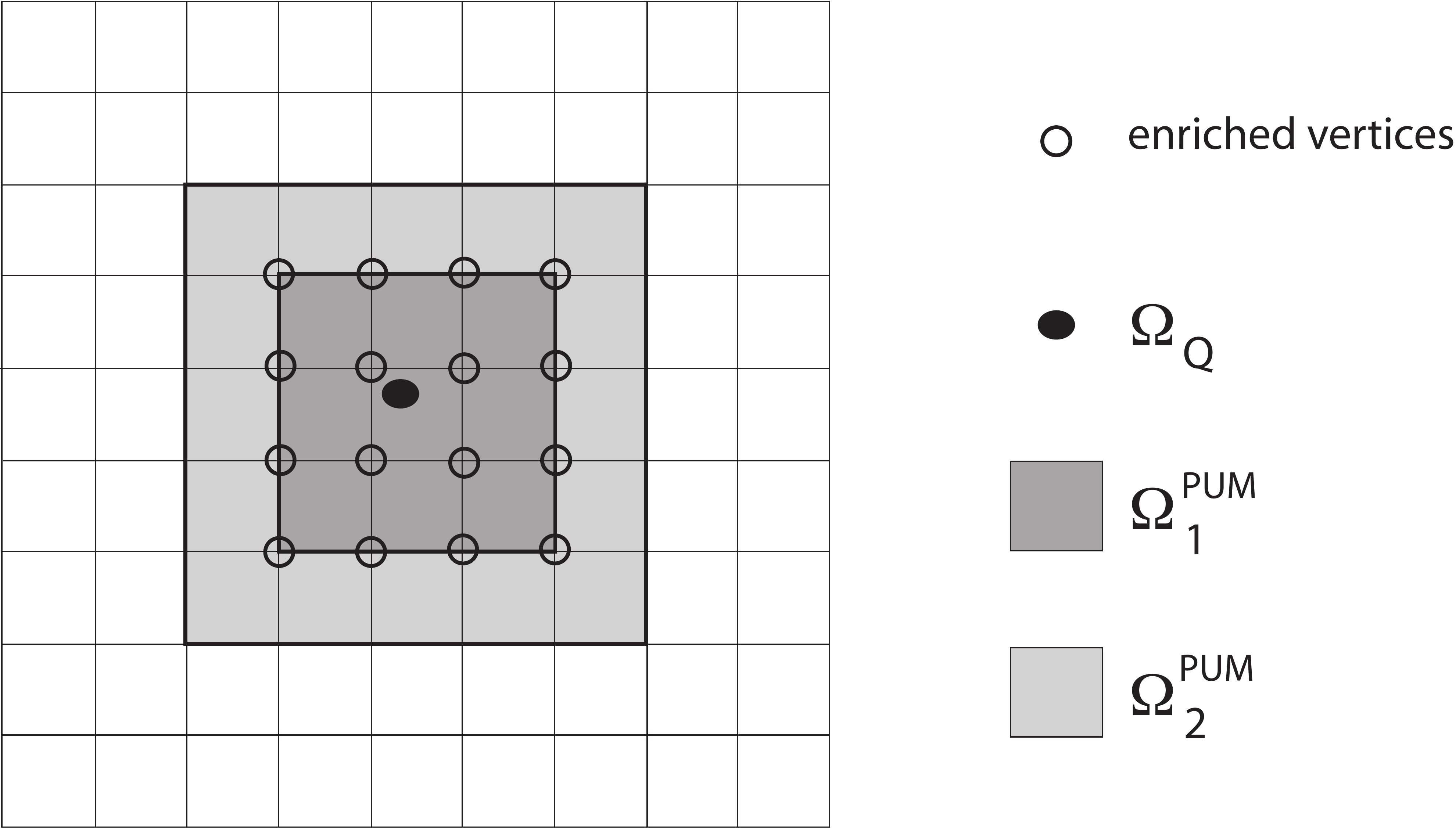}
\end{center}
\caption{Definition of $\Omega_1^{PUM}$ and $\Omega_2^{PUM}$, for a given $\Omega_Q$.}\label{fig:regions}
\end{figure}

Handbook functions are then locally and seamlessly inserted in the overall procedure, as we explain now. Let $\Omega_Q$ be the area of interest over which the quantity $Q$ is defined, i.e. the region which supports the extraction functions $\widetilde{f}_Q$ and $\widetilde{\bq}_Q$ (recall~\eqref{eq:rep_Q} and Remark~\ref{rem:pas_de_gQ}). We next introduce some domain $\Omega_1^{PUM} \subset \Omega_{hand}$, which is the union of some coarse elements of $\mT_H$, and which contains $\Omega_Q$ (see Fig.~\ref{fig:regions}). Its boundary is denoted $\partial \Omega_1^{PUM}$. Let $\mI^{PUM}$ be the subset of vertices in $\mT_H$ which are located in $\overline{\Omega_1^{PUM}}$, so that
$$
\Omega_1^{PUM} = \left\{ \bx \in \Omega; \quad \sum_{i\in \mI^{PUM}} \phi^0_i(\bx) = 1 \right\}.
$$
We next introduce
$$
\Omega_2^{PUM} = \left\{ \bx \in \Omega; \quad \sum_{i\in \mI^{PUM}} \phi^0_i(\bx) \in (0,1) \right\}
$$
and see that $\dis \Omega_1^{PUM} \cup \Omega_2^{PUM} = \text{Supp} \left( \sum_{i\in \mI^{PUM}} \phi^0_i\right) \subset \Omega$.

\begin{remark}
In the numerical examples of Section~\ref{section:defect}, $\Omega_Q$ is a square of size $4 \eps \times 4 \eps$, and we take for $\Omega_1^{PUM}$ the coarse element containing $\Omega_Q$. In the numerical example of Section~\ref{section:flow}, $\Omega_Q$ is a square of size $0.1 \times 0.1$, and the initial coarse meshsize is $H=0.1$. We then take for $\Omega_1^{PUM}$ the union of the coarse element equal to $\Omega_Q$ and of the 8 coarse elements surrounding it.
\end{remark}

The enrichment approach consists in approximating the adjoint solution $\widetilde{u}^\eps$ by 
\begin{equation}\label{eq:PUMdecomp}
\widetilde{u}^\eps_+ = \widetilde{u}^\eps_{hand} \sum_{i\in \mI^{PUM}} \phi^0_i + \widetilde{u}^\eps_{res,H},
\end{equation}
where $\widetilde{u}^\eps_{res,H}$ is an approximate residual solution which is computed using the MsFEM approximation space $V^\eps_H$. We emphasize that
\begin{itemize}
\item no additional degree of freedom is associated with the enrichment part, in contrast to other similar enrichment methods such as XFEM or GFEM (where the approximation would read $\dis \sum_{i\in \mI^{PUM}} \alpha_i \, \phi^0_i \, \widetilde{u}^\eps_{hand} + \widetilde{u}^\eps_{res,H}$ for some additional unknown coefficients $\alpha_i$ and with $\widetilde{u}^\eps_{res,H} \in V^\eps_H$). This is due to the fact that $\widetilde{u}^\eps_{hand}$ aims at representing the local high-gradients due to the singular adjoint loading, which can be fully assessed \textit{a priori};
\item the strategy uses the same MsFEM discretization space $V^\eps_H$ to approximate both the primal and the adjoint problems, even though a higher level of approximation is obtained for the adjoint solution due to the enrichment term. Consequently, all the numerical technicalities (mesh connectivities, stiffness matrix, \ldots) implemented to solve the primal problem can be directly reused to solve the adjoint problem. Only the loadings need to be changed. For this reason, the strategy is called non-intrusive. It fully makes sense in the MsFEM framework where the discretization space is obtained from costly computations performed only once in the offline stage.
\end{itemize}

\medskip

\begin{remark}
The above enrichment methodology can also be used to address goal-oriented error estimation on pointwise quantities of interest (when these are defined), as shown in~\cite{CHA08,LAD10}. In this case, handbook functions correspond to classical Green's functions. They can be analytically computed for homogeneous and isotropic media.
\end{remark}

\medskip

\begin{remark}
We have chosen the notations $\Omega_1^{PUM}$, $\Omega_2^{PUM}$ and $\mI^{PUM}$ in analogy to those introduced in the XFEM and the GFEM approaches. However, here, no partition of unity property is used.
\end{remark}

\medskip

In view of~\eqref{eq:adjointpbweak} and~\eqref{eq:PUMdecomp}, and assuming that $Q_{hand} \equiv Q$ in~\eqref{eq:def_u_hand}, the problem to be solved numerically consists in finding $\widetilde{u}^\eps_{res,H} \in V^\eps_H$ such that
$$
\begin{aligned}
  \forall v \in V^\eps_H, \quad B^{\eps\ast} (\widetilde{u}^\eps_{res,H},v)
  &=
  Q(v) - B^{\eps\ast}\left(\widetilde{u}^\eps_{hand}\sum_{i\in \mI^{PUM}} \phi^0_i,v \right)
  \\
  &= Q(v) - \int_{\Omega_1^{PUM}} \Aa^\eps \nab \widetilde{u}^\eps_{hand} \cdot \nab v - \int_{\Omega_2^{PUM}}\Aa^\eps \nab \left( \widetilde{u}^\eps_{hand} \sum_{i\in \mI^{PUM}} \phi^0_i \right) \cdot \nab v
  \\
&= - \int_{\partial \Omega_1^{PUM}} \left( \Aa^\eps \nab \widetilde{u}^\eps_{hand} \cdot \bn_{12} \right) v - \int_{\Omega_2^{PUM}} \Aa^\eps \nab \left( \widetilde{u}^\eps_{hand} \sum_{i\in \mI^{PUM}} \phi^0_i \right) \cdot \nab v,
\end{aligned}
$$
with $\bn_{12}$ is the outgoing normal vector to $\Omega_1^{PUM}$. In the last line, we have used the fact (assuming e.g. that $\widetilde{g}_Q = 0$ and $\widetilde{\bq}_Q=0$ in~\eqref{eq:rep_Q}) that $- \nab \cdot (\Aa^\eps \nab \widetilde{u}^\eps_{hand}) = \widetilde{f}_Q$ in $\Omega_{hand} \supset \Omega_1^{PUM}$ (see~\eqref{eq:def_u_hand}), and that $\dis \int_{\Omega_1^{PUM}} \widetilde{f}_Q \, v = \int_{\Omega} \widetilde{f}_Q \, v$ because $\widetilde{f}_Q$ is supported in $\Omega_Q \subset \Omega_1^{PUM}$. 

Post-processing the numerical flux $\widetilde{\bq}^\eps_{res,H}=\Aa^\eps \nab \widetilde{u}^\eps_{res,H}$ as explained in Section~\ref{section:admissible_fields}, we compute a flux field $\widehat{\widetilde{\bq}}^\eps_{res,H} \in \widetilde{\mS}_{res}$, with
\begin{multline*}
  \widetilde{\mS}_{res} = \left\{ \widetilde{\bp} \in \mH(div,\Omega); \quad \forall v \in V, \quad \intO \widetilde{\bp} \cdot \nab v = - \int_{\partial \Omega_1^{PUM}} \left( \Aa^\eps \nab \widetilde{u}^\eps_{hand} \cdot \bn_{12} \right) v \right.
  \\
  \left. - \int_{\Omega_2^{PUM}} \Aa^\eps \nab \left( \widetilde{u}^\eps_{hand} \sum_{i\in \mI^{PUM}} \phi^0_i \right) \cdot \nab v \right\}.
\end{multline*}
We next define
$$
\widehat{\widetilde{\bq}}^\eps_{+}= \Aa^\eps\nab \left( \widetilde{u}^\eps_{hand} \sum_{i\in \mI^{PUM}} \phi^0_i \right) + \widehat{\widetilde{\bq}}^\eps_{res,H} \in \widetilde{\mS}.
$$
We also introduce
$$
\widehat{\widetilde{u}}^\eps_+ = \widetilde{u}^\eps_{hand} \sum_{i\in \mI^{PUM}} \phi^0_i + \widehat{\widetilde{u}}^\eps_{res,H},
$$
where $\widehat{\widetilde{u}}^\eps_{res,H} \in V$ is equal to $\widetilde{u}^\eps_{res,H}$ when $V_H^\eps \subset V$, and is obtained after a post-processing of $\widetilde{u}^\eps_{res,H}$ otherwise (see Section~\ref{section:globalestimate}). 

Using the estimate~\eqref{eq:CREgobound2}, we eventually obtain
\begin{multline*}
\left| Q(u^\eps)-Q(u^\eps_H)-\Delta Q-\overline{C}^\eps_{H+} \right|
\le
\frac{1}{2} E_{CRE}(\widehat{u}^\eps_H,\widehat{\bq}^\eps_H) \, E_{CRE}\left(\widehat{\widetilde{u}}^\eps_+,\widehat{\widetilde{\bq}}^\eps_+\right)
\\
=
\frac{1}{2} E_{CRE}(\widehat{u}^\eps_H,\widehat{\bq}^\eps_H) \, E_{CRE}\left(\widehat{\widetilde{u}}^\eps_{res,H},\widehat{\widetilde{\bq}}^\eps_{res,H}\right),
\end{multline*}
with
$$
\overline{C}^\eps_{H+} = \frac{1}{2}\intO (\Aa^\eps)^{-1} \left( \widehat{\bq}^\eps_H - \Aa^\eps \nab \widehat{u}^\eps_H \right) \cdot \left( \widehat{\widetilde{\bq}}_+^\eps+\Aa^\eps \nab \widehat{\widetilde{u}}^\eps_+ \right).
$$
The fully computable quantity $Q(u^\eps_H)+\Delta Q+\overline{C}^\eps_{H+}$ is a corrected approximation of $Q(u^\eps)$. As in~\eqref{eq:goestimate}, the quantity
\begin{equation}\label{eq:defetaQ}
\eta^Q = \max_{\theta=\pm 1} \left| \Delta Q+\overline{C}^\eps_{H+} + \frac{\theta}{2}E_{CRE}(\widehat{u}^\eps_H,\widehat{\bq}^\eps_H) \, E_{CRE}\left(\widehat{\widetilde{u}}^\eps_{res,H},\widehat{\widetilde{\bq}}^\eps_{res,H}\right) \right|
\end{equation}
is a computable upper bound of $|Q(u^\eps)-Q(u^\eps_H)|$. This bound is thought to be sharp as the enrichment of the adjoint solution is supposed to make the term $\dis E_{CRE}\left(\widehat{\widetilde{u}}^\eps_{res,H},\widehat{\widetilde{\bq}}^\eps_{res,H}\right)$ small. When $\eta^Q$ is higher than a given error tolerance $\eta_{tol}^Q$, an adaptive strategy can be introduced. This is the subject of the next section.

\section{Goal-oriented adaptive strategy}\label{section:adaptive}

The objective here is to design an adaptive algorithm that optimally modifies the parameters of the MsFEM discretization in order to reduce the error on $Q$ and reach a given error tolerance $\eta_{tol}^Q$ while keeping the computational cost as small as possible. This is performed by introducing local error indicators associated to each error source and defined from (hierarchical) auxiliary reference problems and CRE properties. 

\subsection{Recast of MsFEM}\label{section:recastMsFEM}

We first introduce a hierarchy of MsFEM-type problems, which is identical to that introduced in~\cite{CHA16b} (we also refer to~\cite{HEN14}, where a similar hierarchy is introduced). We briefly present them below, and refer to~\cite[Section~5.1]{CHA16b} for more details.

We first introduce a projection operator $I_H:V \to V^0_H \subset V$ chosen such that its image is $V^0_H$, where $V^0_H$ is spanned by the Lagrange basis functions associated with each node of $\mT_H$ (see Section~\ref{section:modelpb}). A possible choice for $I_H$ is the Cl\'ement interpolant. We next define the fine-scale space
$$
V^f=\big\{ \phi \in V, \ \ I_H(\phi)=0 \big\}
$$
and observe that $V = V^0_H \oplus V^f$.

For any function $\varphi_H \in V^0_H$, we define the corrector $\Theta^\eps(\varphi_H ) \in V^f$ as the solution to the problem
\begin{equation}\label{eq:correctorpb2}
\forall v \in V^f, \quad B^\eps(\varphi_H+\Theta^\eps(\varphi_H),v) = F(v).
\end{equation}
The associated operator $\Theta^\eps:V^0_H \to V^f$ is an affine operator:
$$
\Theta^\eps(\varphi_H) = \Theta^\eps_0 + \Theta_{\ell}^\eps(\varphi_H) \quad \text{with} \quad 
\left\{
\begin{array}{l}
\Theta^\eps_0 \in V^f \quad \text{s.t.} \quad \forall v \in V^f, \quad B^\eps(\Theta^\eps_0,v) = F(v), \\
\Theta_{\ell}^\eps(\varphi_H) \in V^f \quad \text{s.t.} \quad \forall v \in V^f, \quad B^\eps(\varphi_H+\Theta_{\ell}^\eps(\varphi_H),v) = 0.
\end{array}
\right.
$$
We then define a reconstruction operator $\Lambda^\eps : V^0_H \to V$ as $\Lambda^\eps = I+\Theta^\eps$, as well as the multiscale functional space $\dis \overline{V}^\eps_H =\{\Lambda^\eps(\varphi_H), \ \varphi_H \in V^0_H\}=\Theta^\eps_0 + \text{Span} \big\{\phi^0_i+\Theta_{\ell}^\eps(\phi^0_i)\big\}$, where $\{\phi^0_i\}$ is the set of Lagrange basis functions of $V^0_H$.

\medskip

Consider the following MsFEM-like problem: 
\begin{equation}\label{eq:formulMsFEM02}
\text{Find $\overline{u}^\eps_H \in \overline{V}^\eps_H$ such that, for any $v \in \overline{V}^\eps_H$,}\quad B^\eps(\overline{u}^\eps_H,v) = F(v),
\end{equation}
which is equivalent to the problem
\begin{equation}\label{eq:formulMsFEM0}
\text{Find $\overline{u}^\eps_H \in \overline{V}^\eps_H$ such that, for any $v \in V^0_H$,}\quad B^\eps(\overline{u}^\eps_H,v) = F(v).
\end{equation}
It is easy to show that~\eqref{eq:formulMsFEM02} and~\eqref{eq:formulMsFEM0} are well-posed, and that their solution $\overline{u}^\eps_H$ satisfies the variational formulation~\eqref{eq:refpbweak} of the reference problem. We hence have $\overline{u}^\eps_H = u^\eps$. Note in particular that the fine scale problems~\eqref{eq:correctorpb2} are solved exactly, over the whole domain $\Omega$, and considering the reference loading $F$ on the right-hand side.

\medskip

Starting from the exact formulations~\eqref{eq:formulMsFEM02} or~\eqref{eq:formulMsFEM0}, we now construct a fully discrete MsFEM formulation, on the MsFEM space $V^\eps_H$, by making the following simplifying changes (denoted \textbf{S\#}):
\begin{itemize}
\item \textbf{S1}: replace the right-hand side term with $0$ in the fine-scale problem~\eqref{eq:correctorpb2}, which makes the operator $\Theta^\eps$ linear (i.e. $\Theta^\eps_0=0$) and independent of the loadings. This leads to the low-dimensional vector space $\text{Span} \{\Lambda^\eps(\phi^0_i)\} = \text{Span} \big\{\phi^0_i+\Theta_{\ell}^\eps(\phi^0_i)\big\}$;
\item \textbf{S2}: restrict to $\Omega_i$ the correction $\Theta^\eps(\phi^0_i)$, where $\Omega_i$ is the support of $\phi^0_i$. This leads to the vector space
\begin{equation} \label{eq:def_overlineV1}
  \overline{V}^\eps_{H,1} = \text{Span} \big\{\phi^0_i + \Theta_{\ell}^\eps(\phi^0_i)_{|\Omega_i} \big\}.
\end{equation}
\item \textbf{S3}: localize the fine-scale problem~\eqref{eq:correctorpb2} by replacing $\Omega$ by smaller domains. In practice, these are subdomains $S_K$ that include the element $K$, on the boundary $\partial S_K$ of which we prescribe homogeneous Dirichlet boundary conditions;
\item \textbf{S4}: replace $V^f$ by a discrete $P_1$-Lagrange FE space $V^f_h$ associated to a fine grid $\mT_h$ obtained from a regular refinement of $\mT_H$. The fine mesh size $h$ should be chosen so that oscillations of the data can be accurately captured.
\end{itemize}

\medskip

Applying the simplifications~\textbf{S1}-\textbf{S4} on~\eqref{eq:formulMsFEM02} or~\eqref{eq:formulMsFEM0}, we obtain a fully discrete MsFEM formulation. For each $K \in \mT_H$ and each basis function $\phi^0_i \in V^0_H$, the local correction $\Theta^\eps_h(\phi^0_{i|K})\in V^f_h(S_K)$ is defined as the solution to
\begin{equation*}
\forall v \in V^f_h(S_K), \quad \int_{S_K} \Aa^\eps \nab \big( \pi(\phi^0_{i|K}) + \Theta^\eps_h(\phi^0_{i|K})\big)\cdot \nab v =0,
\end{equation*}
where $\pi(\phi^0_{i|K})$ is the affine function on $S_K$ which is equal to $\phi^0_{i|K}$ on $K$, and $V^f_h(S_K)$ is a discrete subspace of $V^f(S_K) = \mH^1_0(S_K)$. The local correctors $\Theta^\eps_h(\phi^0_{i|K})$ are then connected together in the global corrector $\dis \Theta^\eps_h(\phi^0_i) := \sum_{K \in \mT_H} \Theta^\eps_h(\phi^0_{i|K})_{|K}$. The MsFEM functional space eventually reads $V^\eps_H = \text{Span} \{ \phi^0_i+\Theta^\eps_h(\phi^0_i) \}$.

\medskip

\begin{remark}
Applying the simplifications~\textbf{S1}-\textbf{S4} on~\eqref{eq:formulMsFEM02} with the choice $S_K=K$ yields the MsFEM approach without oversampling given by~\eqref{eq:MsFEMpb}--\eqref{eq:localpbMsFEM2}.
\end{remark}

\subsection{Splitting of error sources and error indicators}\label{section:errorsources}

The MsFEM formulation given in Section~\ref{section:recastMsFEM} shows that 3 parameters can be tuned in an adaptive procedure to improve the accuracy of $Q(u^\eps_H)$:
\begin{itemize}
\item[$\bullet$] the size $H_K$ of the elements of the coarse mesh $\mT_H$. Simplification~\textbf{S1} amounts to ignoring the influence of $f$ and $g$ in micro-scale equations. This assumption is valid provided $H_K$ is comparable with the characteristic length of spatial variations of the external loading and macroscale solution. As detailed in~\cite[Remark~16]{CHA16b}, the smaller $H_K$ is, the better Simplification~\textbf{S1} is justified. The validity of Simplification~\textbf{S2} is also related to $H_K$. If the set $\{H_K\}$ is not chosen correctly, the coarse-scale discretization error can dominate the overall approximation error;
\item[$\bullet$] the size $d_K$ (minimum patch radius) of the computational domains $S_K$ used to solve the micro-scale equations obtained after Simplification~\textbf{S3}. If the set $\{d_K\}$ is not chosen correctly, the artificial boundary conditions set
on $\partial S_K$ may strongly dominate the overall approximation error;
\item[$\bullet$] the local size $h_K$ of the fine mesh used to solve the micro-scale equations (in each $S_K$) after applying Simplification~\textbf{S4}. If the set $\{h_K\}$ is not chosen correctly, a fine-scale discretization error can dominate the overall approximation error.
\end{itemize}

\medskip

We wish to define a robust adaptive algorithm that is able: (i) to detect regions where the MsFEM approximation is not sufficiently accurate for the quantity of interest $Q(u^\eps)$; (ii) to adapt locally and optimally the relevant MsFEM parameters among the three ones mentioned above.

The item (i) can be directly addressed with a spatial decomposition, over the macro mesh $\mT_H$, of the \textit{a posteriori} error estimate $\eta_Q$ defined by~\eqref{eq:defetaQ}:
\begin{equation}\label{eq:splitindicator_a}
  \eta^Q \leq \sum_{K \in \mT_H} \eta^Q_K
\end{equation}
with
\begin{multline}\label{eq:splitindicator_b}
\eta^Q_K = \left| B^\eps_{|K}\left(\widehat{u}^\eps_H-u^\eps_H,\widehat{\widetilde{u}}^\eps_+\right) + \frac{1}{2} \int_K (\Aa^\eps)^{-1}(\widehat{\bq}^\eps_H - \Aa^\eps \nab \widehat{u}^\eps_H) \cdot \left( \widehat{\widetilde{\bq}}_+^\eps + \Aa^\eps \nab \widehat{\widetilde{u}}^\eps_+ \right) \right.
\\
\left. + \frac{\theta_{\rm max}}{2} \, E_{CRE|K}(\widehat{u}^\eps_H,\widehat{\bq}^\eps_H) \, E_{CRE|K}\left(\widehat{\widetilde{u}}^\eps_{res,H},\widehat{\widetilde{\bq}}^\eps_{res,H}\right) \right|
\end{multline}
where $\theta_{\rm max}$ is the maximizer in~\eqref{eq:defetaQ}.
The contribution $\eta^Q_K$, which is fully computable, is considered as a local error indicator in the coarse element $K$.

\medskip

\begin{figure}[h]
\begin{center}
\includegraphics[width=95mm]{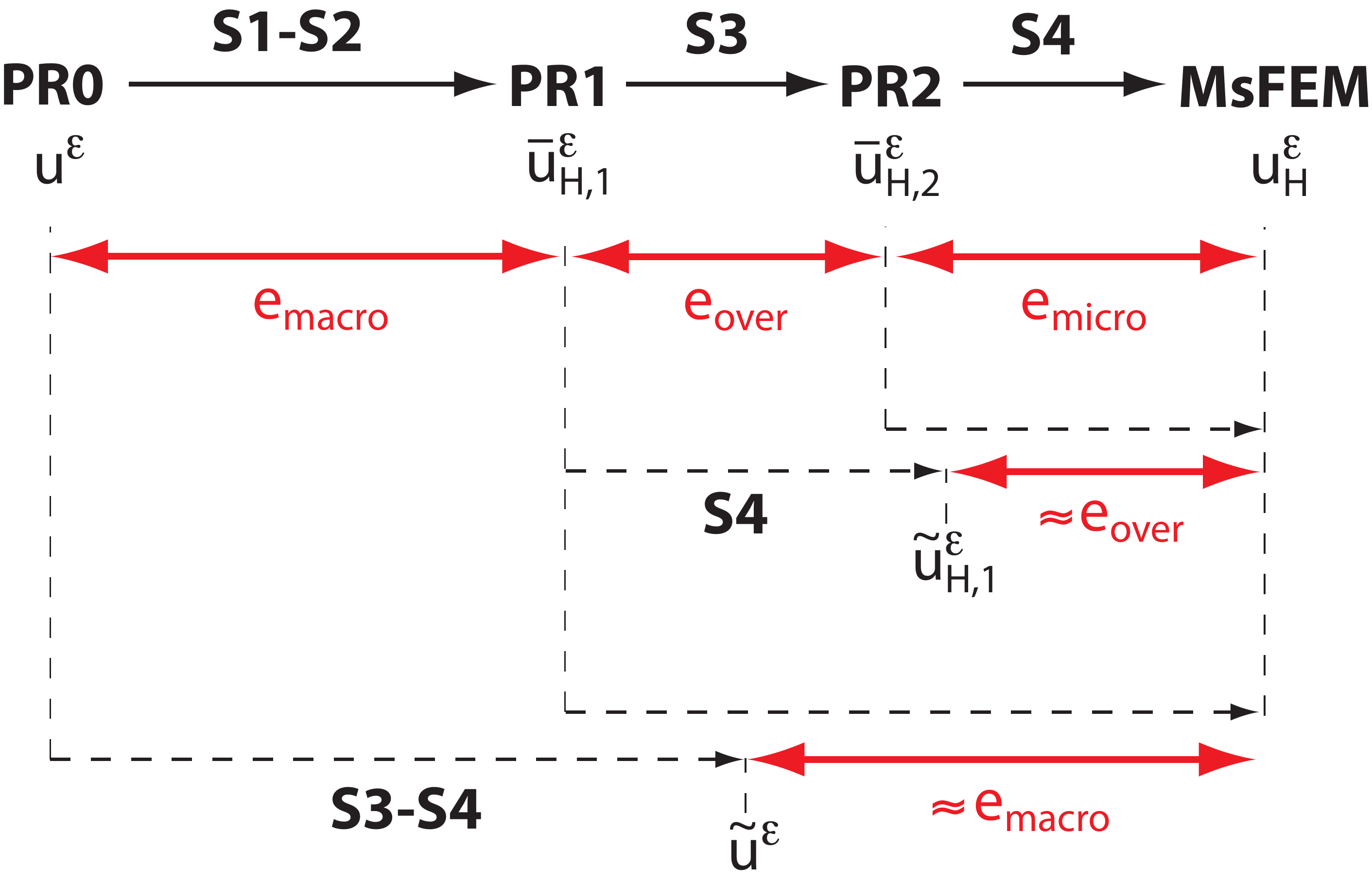}
\end{center}
\caption{Illustration of the strategy used to define indicators on each error source.}
\label{fig:errorsources}
\end{figure}

The above item (ii) requires the definition of specific error indicators associated to each individual error source. In order to set up such error indicators, we follow our approach described in details in~\cite[Section~5.2]{CHA16b} and introduce two intermediate reference problems:
\begin{itemize}
\item the first one, denoted \textbf{PR1} in the following, is obtained from the initial reference problem~\eqref{eq:refpbweak} (recast as~\eqref{eq:formulMsFEM02} or~\eqref{eq:formulMsFEM0} and denoted \textbf{PR0}) by applying Simplifications~\textbf{S1} and \textbf{S2}. The functional space associated with \textbf{PR1} is the space $\overline{V}^\eps_{H,1}$ defined by~\eqref{eq:def_overlineV1}, and the solution to \textbf{PR1} is denoted $\overline{u}_{H,1}^\eps$. Note that, in the problem \textbf{PR1}, fine-scale computations are still defined over the whole domain $\Omega$;
\item the second one, denoted \textbf{PR2} in the following, is obtained from the initial reference problem \textbf{PR0} by applying Simplifications \textbf{S1}, \textbf{S2} and \textbf{S3}. The functional space associated with \textbf{PR2} is denoted $\overline{V}^\eps_{H,2}$, and the solution to \textbf{PR2} is denoted $\overline{u}_{H,2}^\eps$. In the problem \textbf{PR2}, the fine-scale computations are defined over the subdomains $S_K$ and are performed exactly (no discretization of $V^f(S_K) = \mH^1_0(S_K)$ is yet introduced).
\end{itemize}
Note that $\overline{u}_{H,1}^\eps$ and $\overline{u}_{H,2}^\eps$ may be nonconforming. We next write
\begin{equation*}
Q(u^\eps-u^\eps_H) = Q(u^\eps-\overline{u}_{H,1}^\eps)+Q(\overline{u}_{H,1}^\eps-\overline{u}_{H,2}^\eps)+Q(\overline{u}_{H,2}^\eps-u^\eps_H), \quad \text{hence} \quad |Q(u^\eps-u^\eps_H)| \le e^Q_{macro}+e^Q_{over}+e^Q_{micro},
\end{equation*}
where
\begin{itemize}
\item[$\bullet$] $e^Q_{macro}=|Q(u^\eps-\overline{u}_{H,1}^\eps)|$ is the part of the error on $Q$ due to the coarse-scale discretization;
\item[$\bullet$] $e^Q_{over}=|Q(\overline{u}_{H,1}^\eps-\overline{u}_{H,2}^\eps)|$ is the part of the error on $Q$ due to the localization on $S_K$ of the problems~\eqref{eq:correctorpb2};
\item[$\bullet$] $e^Q_{micro}=|Q(\overline{u}_{H,2}^\eps-u^\eps_H)|$ is the part of the error on $Q$ due to the fine-scale discretization.
\end{itemize}
The strategy to assess $e^Q_{macro}$, $e^Q_{over}$ and $e^Q_{micro}$ uses error indicators constructed from the CRE concept, with admissibility in the weaker sense of the associated reference problems (both for the primal and the adjoint problems). This is illustrated in Fig.~\ref{fig:errorsources}. We have used the same ideas in~\cite[Section~5.2]{CHA16b} to define specific error indicators (on the whole field $u^\eps$) associated to each individual error source. Error indicators $\eta^Q_{macro}$, $\eta^Q_{over}$ and $\eta^Q_{micro}$ associated to each error source can then be computed. We emphasize that the expensive computations (involving the resolution of highly oscillatory problems) can all be performed in the \textit{offline} stage of our approach. These indicators are used to drive a goal-oriented greedy adaptive algorithm detailed in the next section.

\subsection{Adaptive algorithm}\label{section:algorithm}

\begin{figure}[H]
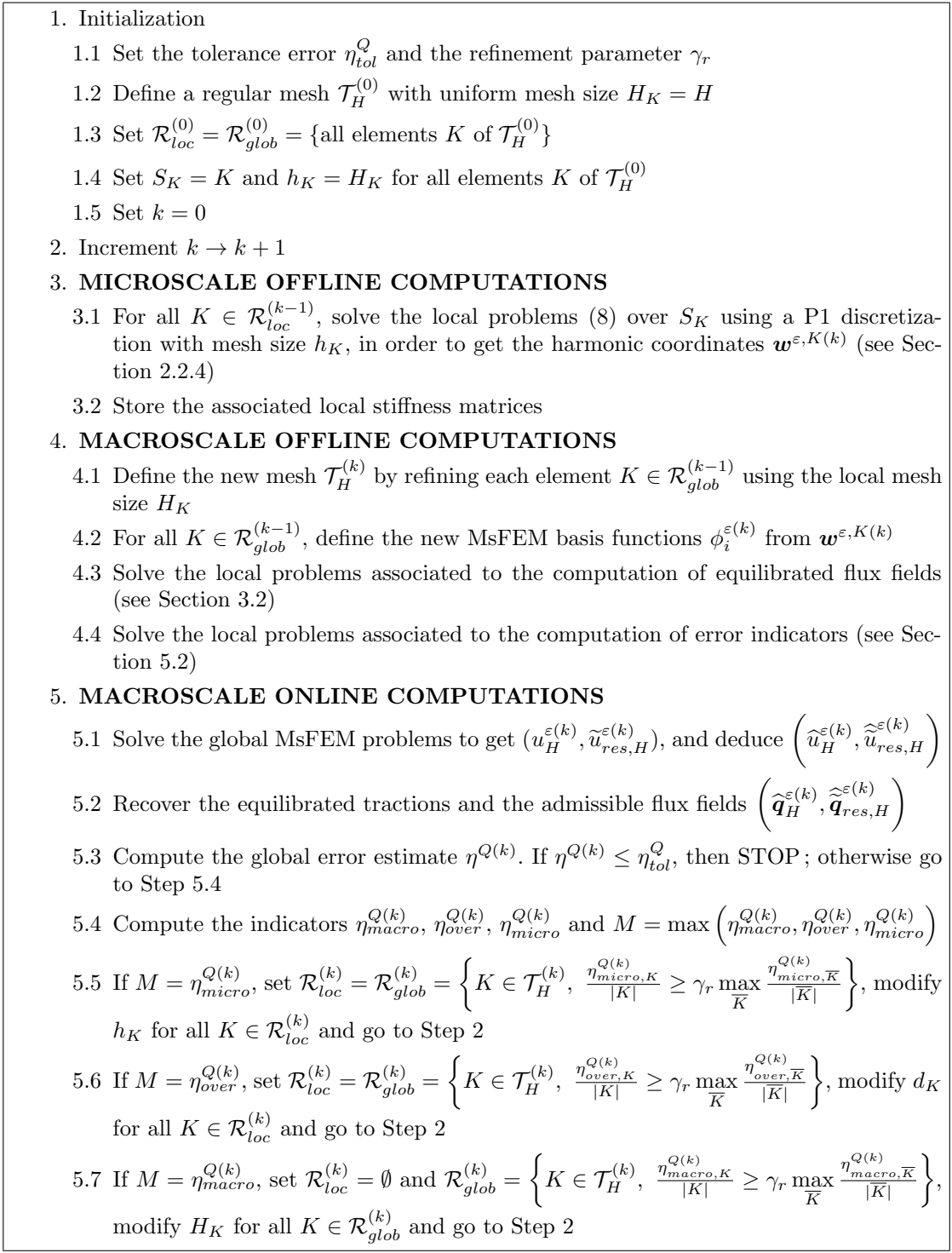

\begin{center}
\boxed{
\begin{minipage}{0.9\textwidth}
\begin{enumerate}
\item[1.] Initialization 
\begin{enumerate}
\item[1.1] Set the tolerance error $\eta^Q_{tol}$ and the refinement parameter $\gamma_r$
\item[1.2] Define a regular mesh $\mT^{(0)}_H$ with uniform mesh size $H_K=H$
\item[1.3] Set $\mR^{(0)}_{loc}=\mR^{(0)}_{glob}=\{\text{all elements $K$ of}\; \mT^{(0)}_H\}$
\item[1.4] Set $S_K=K$ and $h_K=H_K$ for all elements $K$ of $\mT^{(0)}_H$
\item[1.5] Set $k=0$
\end{enumerate}
\item[2.] 
Increment $k \rightarrow k+1$
\item[3.] \textbf{MICROSCALE OFFLINE COMPUTATIONS}
\begin{enumerate}
\item[3.1] For all $K \in \mR^{(k-1)}_{loc}$, solve the local problems~\eqref{eq:localpbho} over $S_K$ using a P1 discretization with mesh size $h_K$, in order to get the harmonic coordinates $\bw^{\eps,K(k)}$ (see Section~\ref{section:highorder})
\item[3.2] Store the associated local stiffness matrices
\end{enumerate}
\item[4.] \textbf{MACROSCALE OFFLINE COMPUTATIONS}
\begin{enumerate}
\item[4.1] Define the new mesh $\mT_H^{(k)}$ by refining each element $K \in \mR^{(k-1)}_{glob}$ using the local mesh size $H_K$
\item[4.2] For all $K \in \mR^{(k-1)}_{glob}$, define the new MsFEM basis functions $\phi_i^{\eps(k)}$ from $\bw^{\eps,K(k)}$
\item[4.3] Solve the local problems associated to the computation of equilibrated flux fields (see Section~\ref{section:admissible_fields})
\item[4.4] Solve the local problems associated to the computation of error indicators (see Section~\ref{section:errorsources})
\end{enumerate}
\item[5.] \textbf{MACROSCALE ONLINE COMPUTATIONS}
\begin{enumerate}
\item[5.1] Solve the global MsFEM problems to get $(u^{\eps(k)}_{H},\widetilde{u}^{\eps(k)}_{res,H})$, and deduce $\left(\widehat{u}^{\eps(k)}_H,\widehat{\widetilde{u}}^{\eps(k)}_{res,H}\right)$ 
\item[5.2] Recover the equilibrated tractions and the admissible flux fields $\left(\widehat{\bq}^{\eps(k)}_H,\widehat{\widetilde{\bq}}^{\eps(k)}_{res,H}\right)$
\item[5.3] Compute the global error estimate $\eta^{Q(k)}$. If $\eta^{Q(k)} \le \eta^Q_{tol}$, then STOP; otherwise go to Step 5.4
\item[5.4] Compute the indicators $\eta^{Q(k)}_{macro}$, $\eta^{Q(k)}_{over}$, $\eta^{Q(k)}_{micro}$ and $M=\max \left(\eta^{Q(k)}_{macro},\eta^{Q(k)}_{over},\eta^{Q(k)}_{micro} \right)$
\item[5.5] If $M=\eta^{Q(k)}_{micro}$, set $\mR^{(k)}_{loc}=\mR^{(k)}_{glob}=\left\{ K \in \mT_H^{(k)}, \ \frac{\eta^{Q(k)}_{micro,K}}{|K|} \ge \gamma_r \, \underset{\overline{K}}{\max} \, \frac{\eta^{Q(k)}_{micro,\overline{K}}}{|\overline{K}|} \right\}$, modify $h_K$ for all $K\in\mR^{(k)}_{loc}$ and go to Step 2
\item[5.6] If $M=\eta^{Q(k)}_{over}$, set $\mR^{(k)}_{loc}=\mR^{(k)}_{glob}=\left\{K \in \mT_H^{(k)}, \ \frac{\eta^{Q(k)}_{over,K}}{|K|}\ge \gamma_r \, \underset{\overline{K}}{\max} \, \frac{\eta^{Q(k)}_{over,\overline{K}}}{|\overline{K}|} \right\}$, modify $d_K$ for all $K\in\mR^{(k)}_{loc}$ and go to Step 2
\item[5.7] If $M=\eta^{Q(k)}_{macro}$, set $\mR^{(k)}_{loc}=\emptyset$ and $\mR^{(k)}_{glob}=\left\{ K \in \mT_H^{(k)}, \ \frac{\eta^{Q(k)}_{macro,K}}{|K|}\ge \gamma_r \, \underset{\overline{K}}{\max} \, \frac{\eta^{Q(k)}_{macro,\overline{K}}}{|\overline{K}|} \right\}$, modify $H_K$ for all $K\in\mR^{(k)}_{glob}$ and go to Step 2
\end{enumerate}
\end{enumerate}
\end{minipage}
}
\caption{Goal-oriented adaptive algorithm to drive the MsFEM computations.}
\label{fig:adaptivescheme}
\end{center}
\end{figure}

From the error estimate $\eta^Q$ defined by~\eqref{eq:defetaQ} and the error indicators ($\eta^Q_{macro}$, $\eta^Q_{over}$, $\eta^Q_{micro}$) defined in Section~\ref{section:errorsources}, we design a goal-oriented greedy adaptive algorithm that yields an adaptive MsFEM approach. The overall algorithm is presented in Fig.~\ref{fig:adaptivescheme} and is similar to the algorithm designed in~\cite[Section~5.3]{CHA16b} for the control of the MsFEM error measured in the energy norm.

The algorithm is based on a splitting~\eqref{eq:splitindicator_a} of the goal-oriented error estimator $\eta^Q$ into local contributions $\{\eta^Q_{macro,K},\eta^Q_{over,K},\eta^Q_{micro,K}\}$ in each element $K$ of the mesh $\mT_H$. These are used to mark coarse elements of the mesh $\mT_H$ for further enrichment, and thus to drive the adaptive process by selecting optimal MsFEM parameters $(H_K,d_K,h_K)$ in each coarse region $K$ to compute both primal and adjoint approximate solutions.

Each iteration in the algorithm is made of three groups of computations: (i) microscale (possibly costly) computations that can be performed in the \textit{offline} phase of the iteration; (ii) macroscale computations that can be performed in this \textit{offline} phase; (iii) macroscale computations that need to be performed in the \textit{online} phase of the iteration, but which are inexpensive.

The adaptive procedure is associated with the following considerations:
\begin{itemize} 
\item[$\bullet$] the algorithm is initialized with the coarsest MsFEM configuration that can be employed (regular coarse mesh, no oversampling, fine-scale problems solved with a single macro element). In this configuration, the diffusion tensor that locally appears to compute the MsFEM shape functions is the average of $\Aa^\eps$ on $K$, as discussed in~\cite[Remark~20]{CHA16b};
\item[$\bullet$] when modifying the parameters $h_K$ from their initial value $h_K=H_K$ (in Step 5.5), two values are considered in the adaptive process: $h_K=\eps/5$ or $h_K=\eps/20$ (finest mesh size at the microscale);
\item[$\bullet$] when modifying the parameters $d_K$ from their initial value $d_K=H_K$ i.e. $S_K=K$ (in Step 5.6), the oversampling size is determined by adding layers of progressive thickness $\eps$, $2\eps$, $3\eps$, \dots, around the element $K$. Such an adaptation is not performed when $h_K=H_K$;
\item[$\bullet$] refining the mesh $\mT_H$, i.e. modifying the parameters $H_K$ (in Step 5.7), is performed using a quadtree (or octree in 3D) method with nested elements. This requires to handle hanging nodes;
\item[$\bullet$] the adapted parameters $h_K$, $d_K$ and $H_K$ are defined from classical adaptive strategies (see~\cite{LAD04}), based on the convergence rates given by the \textit{a priori} estimates~\eqref{eq:aprioriestimate1} or~\eqref{eq:aprioriestimate2}, and with an error target $\eta^Q_{macro/over/micro}=\eta^Q_{tol}/3$;
\item[$\bullet$] the technique proposed in~\cite{ALL06} and based on composition rules (see Section~\ref{section:highorder}) is beneficially used throughout the adaptive process. It enables independent computations without coming back to the fine scale \textit{offline} computations (no additional costly computations). The technique is used in its $p$-refinement version in Steps 4.3 and 4.4, while it is used in its $H$-refinement version in Step 4.2 (when the macro mesh has to be refined). In short, this technique allows us to refine the coarse mesh without the need to solve new fine-scale problems of the type~\eqref{eq:localpbho}.
\end{itemize}

\medskip

The adaptive algorithm allows us to adjust automatically (and hopefully optimally) the calculation parameters during the simulation. We expect that it provides MsFEM parameters which ensure a correct balance between computational cost and accuracy on the output of interest. In particular, it should lead to fine-scale solutions only in the parts of $\Omega$ which have an impact on the value of the quantity of interest. As will be shown below, and as expected, it also leads (for the same tolerance on the error) to less expensive computations than the adaptive algorithm introduced in~\cite{CHA16b} to drive the MsFEM computations according to the error measured in the energy norm.

\section{Numerical results}\label{section:results}

We consider two test-cases. The first one corresponds to a structure with a local defect, for which we consider two possible loadings (see Section~\ref{section:defect}). The second one is a flow problem in fractured porous media (see Section~\ref{section:flow}), which has already been considered in~\cite{CHU16}. For all these test-cases, the conclusions of our numerical investigations are the same:
\begin{itemize}
\item It is worth working with a goal-oriented adaptive strategy rather than an adaptive strategy based on the global error: for the same error tolerance, the computational cost is much lower, in the sense that the obtained mesh (satisfying the prescribed error tolerance) is much coarser. 
\item As the mesh is iteratively adapted, the error on the quantity of interest monotically decreases until reaching the prescribed tolerance.
\item The error estimate $\eta^Q$ is always an upper bound on the exact error, in view of~\eqref{eq:defetaQ}. We see below that it only slightly overestimates the error (roughly by 10\% on the considered numerical examples), showing thus its accuracy. 
\end{itemize}

\subsection{Structure with a local defect}\label{section:defect}

To start with, we consider a highly oscillatory problem with defects, as introduced and analyzed in details in~\cite{BLA12,BLA15a,BLA15b} (see also the recent works~\cite{defauts-BLL1,defauts-BLL2}). We assume that the oscillating coefficient is a local perturbation (which models the defect) of a function with simple geometric properties, for instance a periodic function. We thus consider problem~\eqref{eq:refpb} with a (scalar, to simplify) coefficient $A^\eps(\bx) = a(\bx/\eps)$ with $a=a_{per}+b$ where $b \in L^r(\RR^d)$ is the local perturbation. The homogenization of such a problem is addressed in~\cite{BLA12} for the Hilbertian case $r=2$ and in~\cite{BLA15a,BLA15b} for the cases $1 <r <+\infty$ (see also~\cite{defauts-BLL1,defauts-BLL2}). In general, the defect does not modify the macroscopic behavior (i.e. the homogenized solution) but it affects the solution locally, and at the small scale. 

\begin{figure}[h]
\begin{center}
\includegraphics[width=75mm]{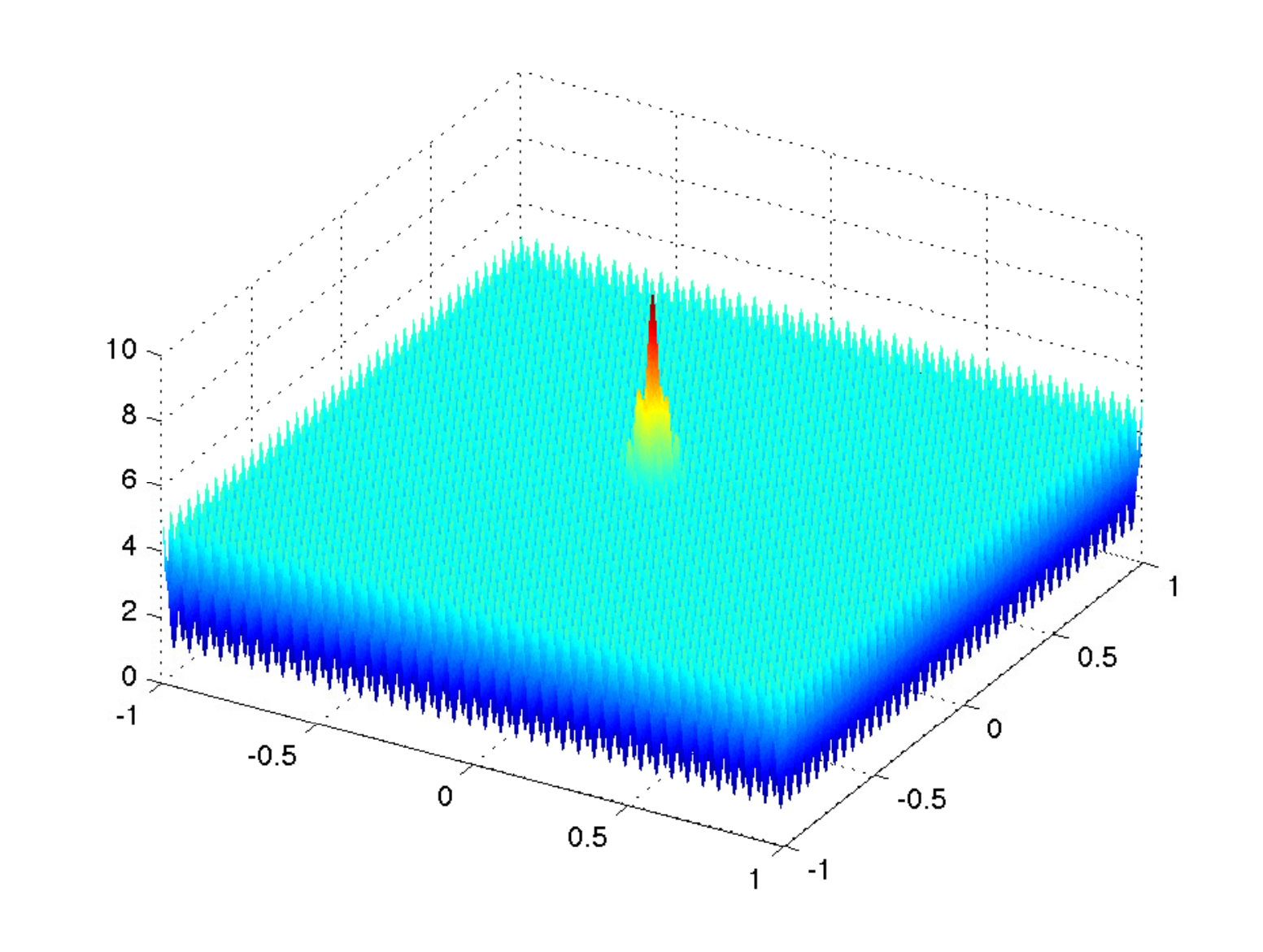}
\includegraphics[width=75mm]{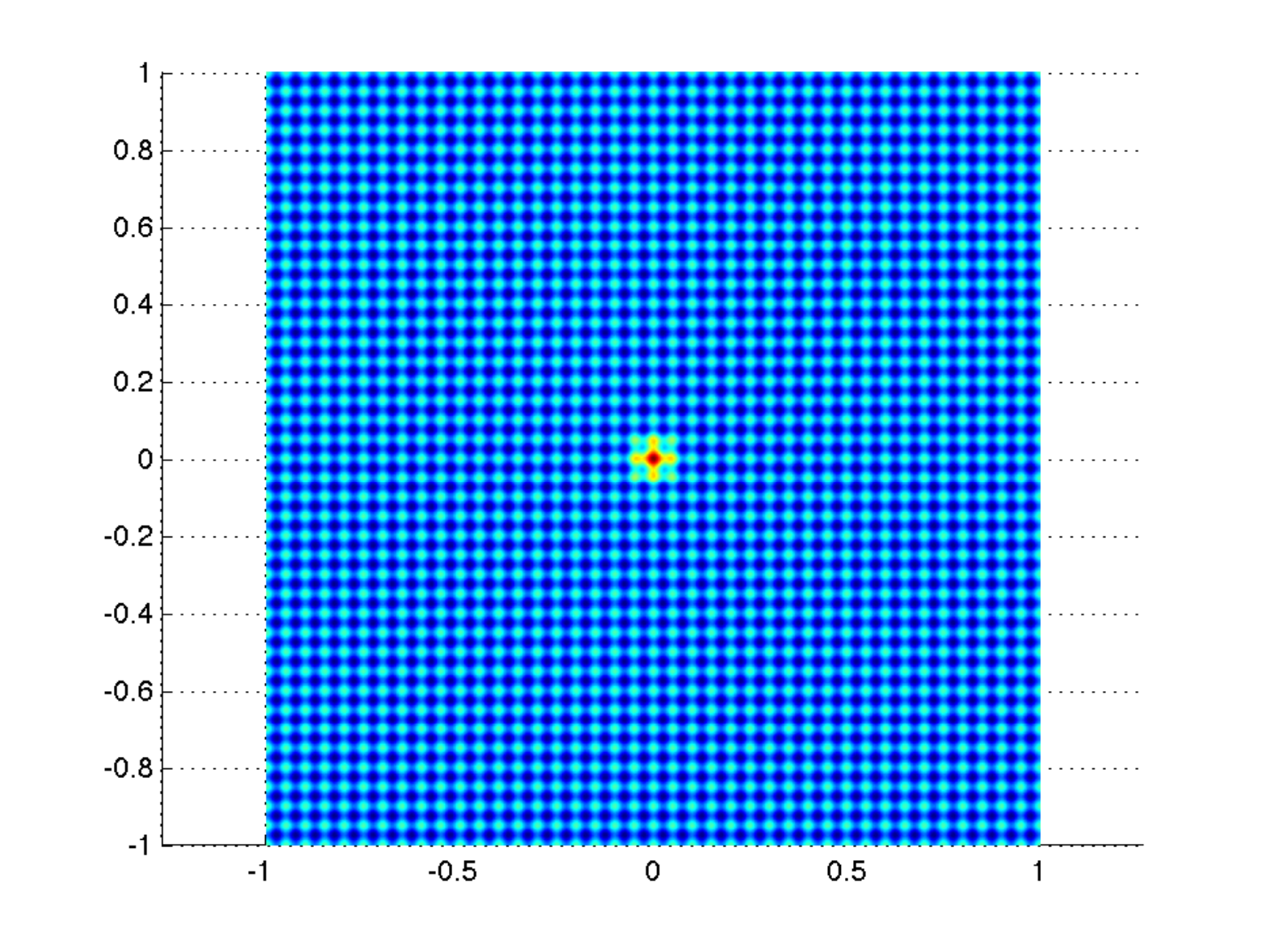}
\end{center}
\caption{Defect problem: evolution of $A^\eps$, modeling a defect, in the domain $\Omega$ (we set $\eps = 1/20$).}
\label{fig:evolAE2Ddefect}
\end{figure}

Inspired by these works, we consider here the following example. We set $\Omega=(-1,1)^2$, choose $\eps=1/N$ for some integer $N$, and consider problem~\eqref{eq:refpb} with $\Aa^\eps(\bx) = A^\eps(\bx) \II$, where $A^\eps(\bx) = a_{per}(\bx/\eps) + b(\bx/\eps)$ with
$$
a_{per}(x_1,x_2)=3+\cos(2\pi x_1)+\cos(2\pi x_2), \qquad b(x_1,x_2) = 5 \exp(-(x_1^2+x_2^2)).
$$
In the following, we choose $\eps=1/20$. The evolution of $A^\eps$ is shown on Fig.~\ref{fig:evolAE2Ddefect}.

\subsubsection{Sinusoidal loading}
\label{sec:sinus}

\begin{figure}[H]
\begin{center}
\includegraphics[width=75mm]{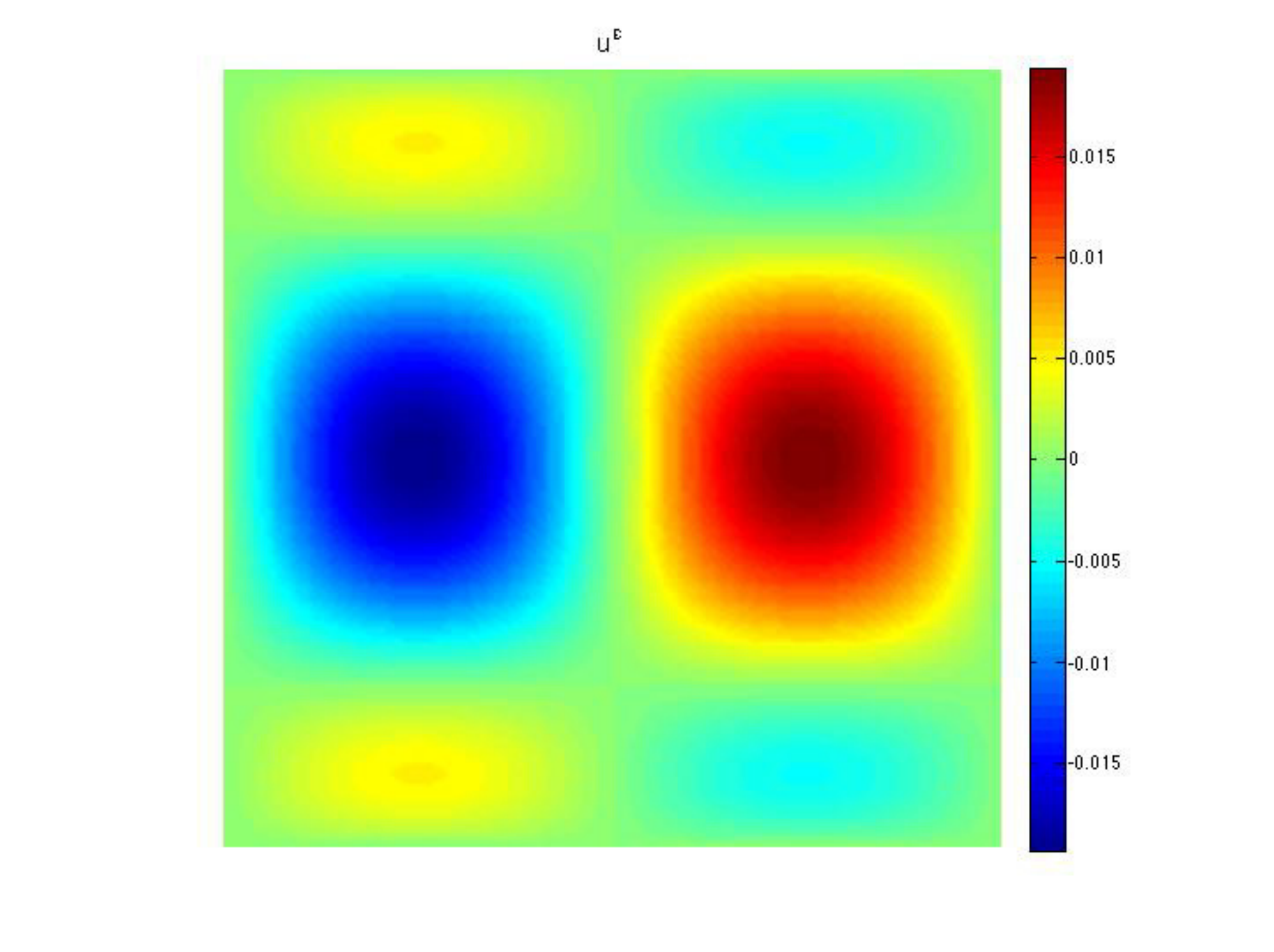}
\includegraphics[width=75mm]{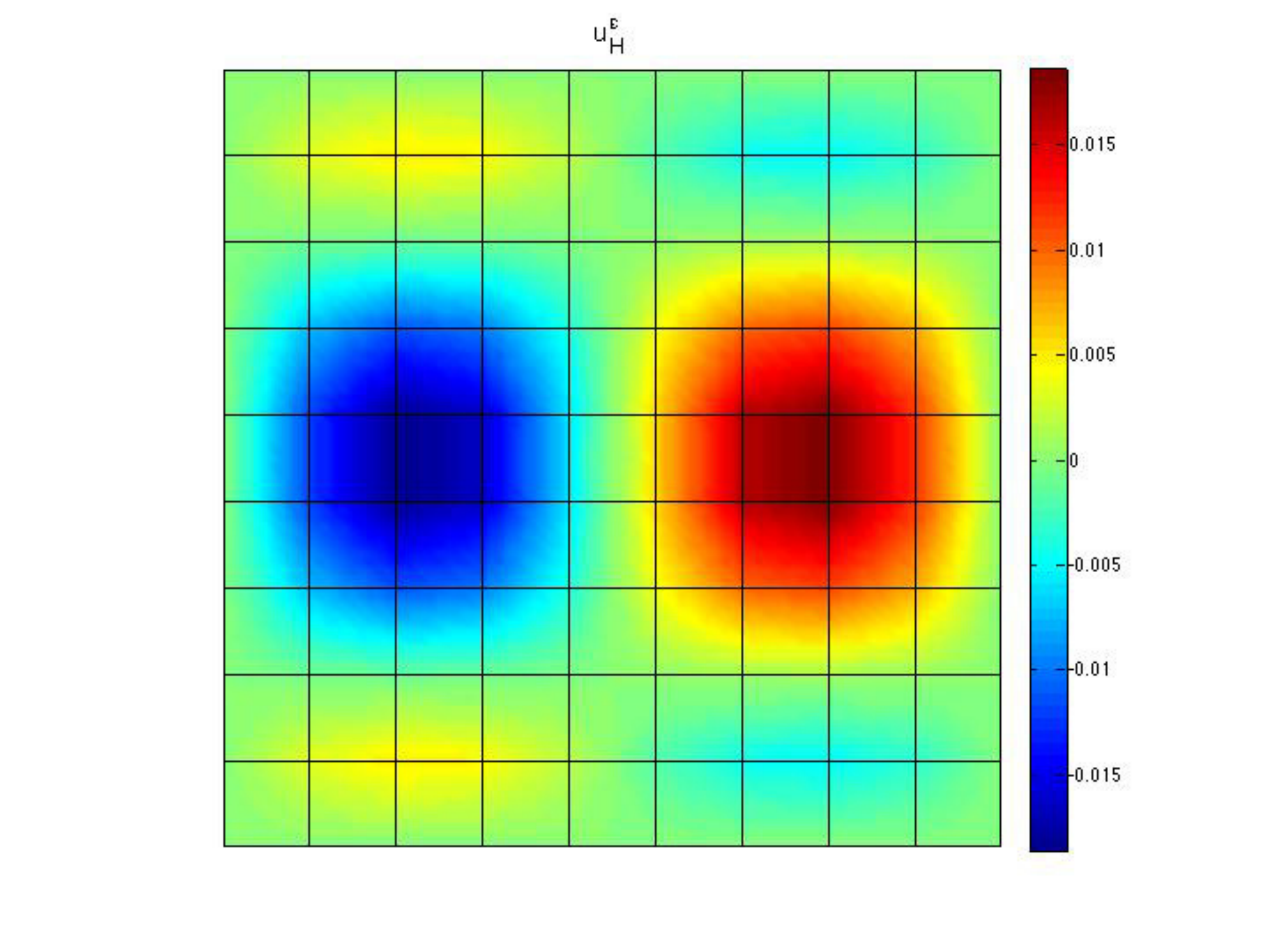}
\end{center}
\caption{Defect problem, sinusoidal loading: exact solution $u^\eps$ (left) and MsFEM solution $u^\eps_H$ (right).}
\label{fig:solMsFEM2Ddefectu}
\end{figure}

\begin{figure}[H]
\begin{center}
\includegraphics[width=75mm]{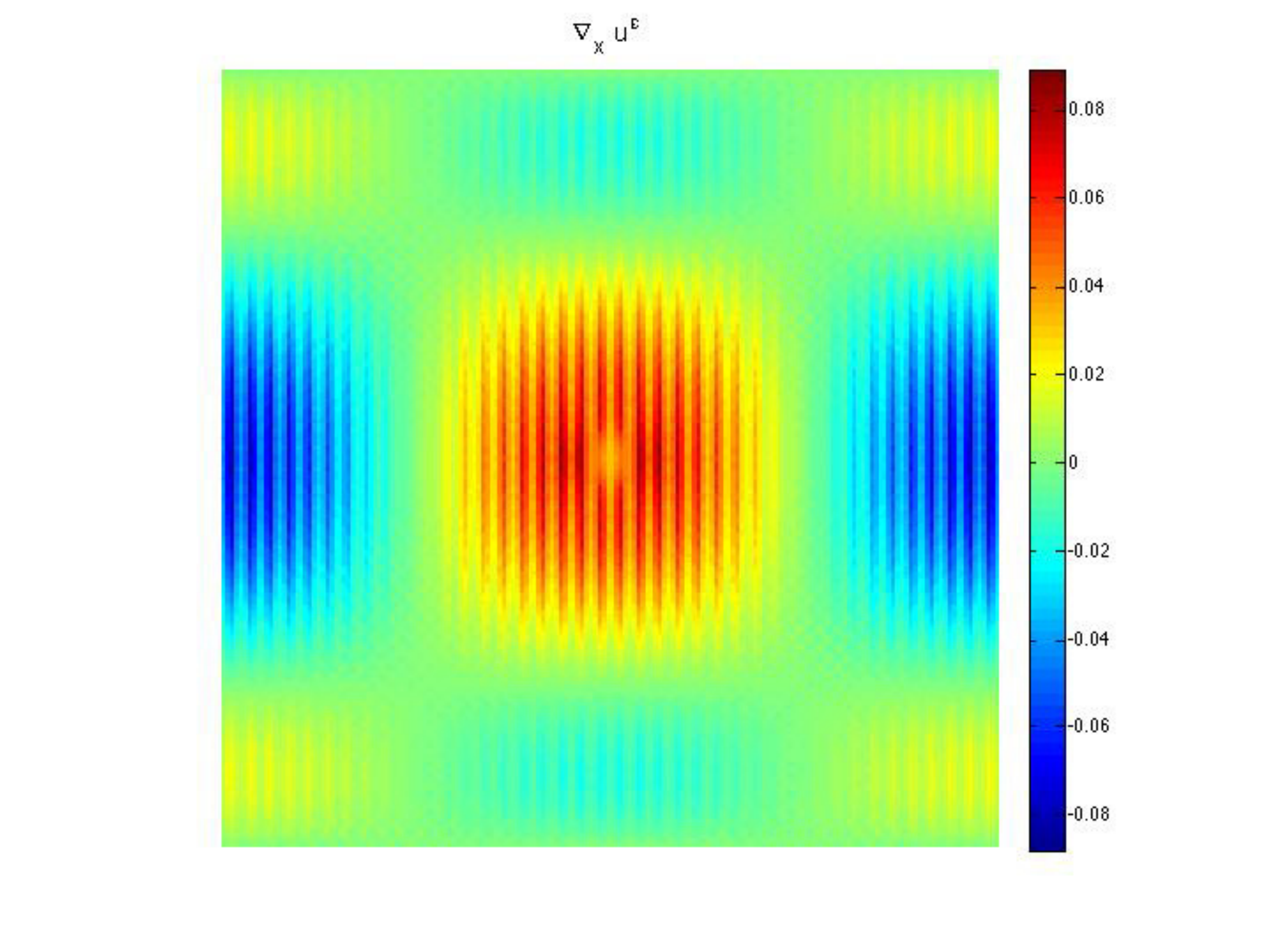}
\includegraphics[width=75mm]{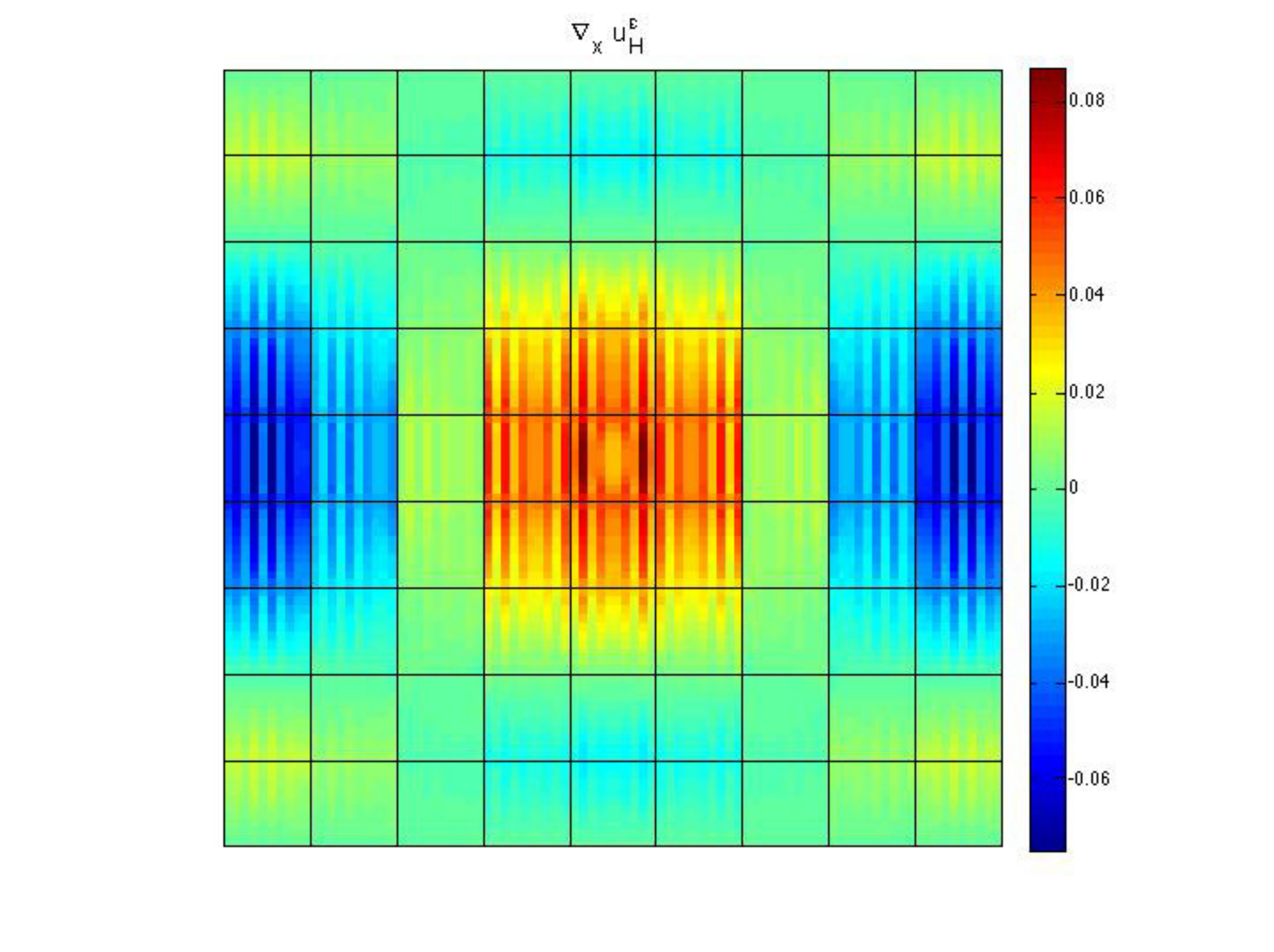}\\
\includegraphics[width=75mm]{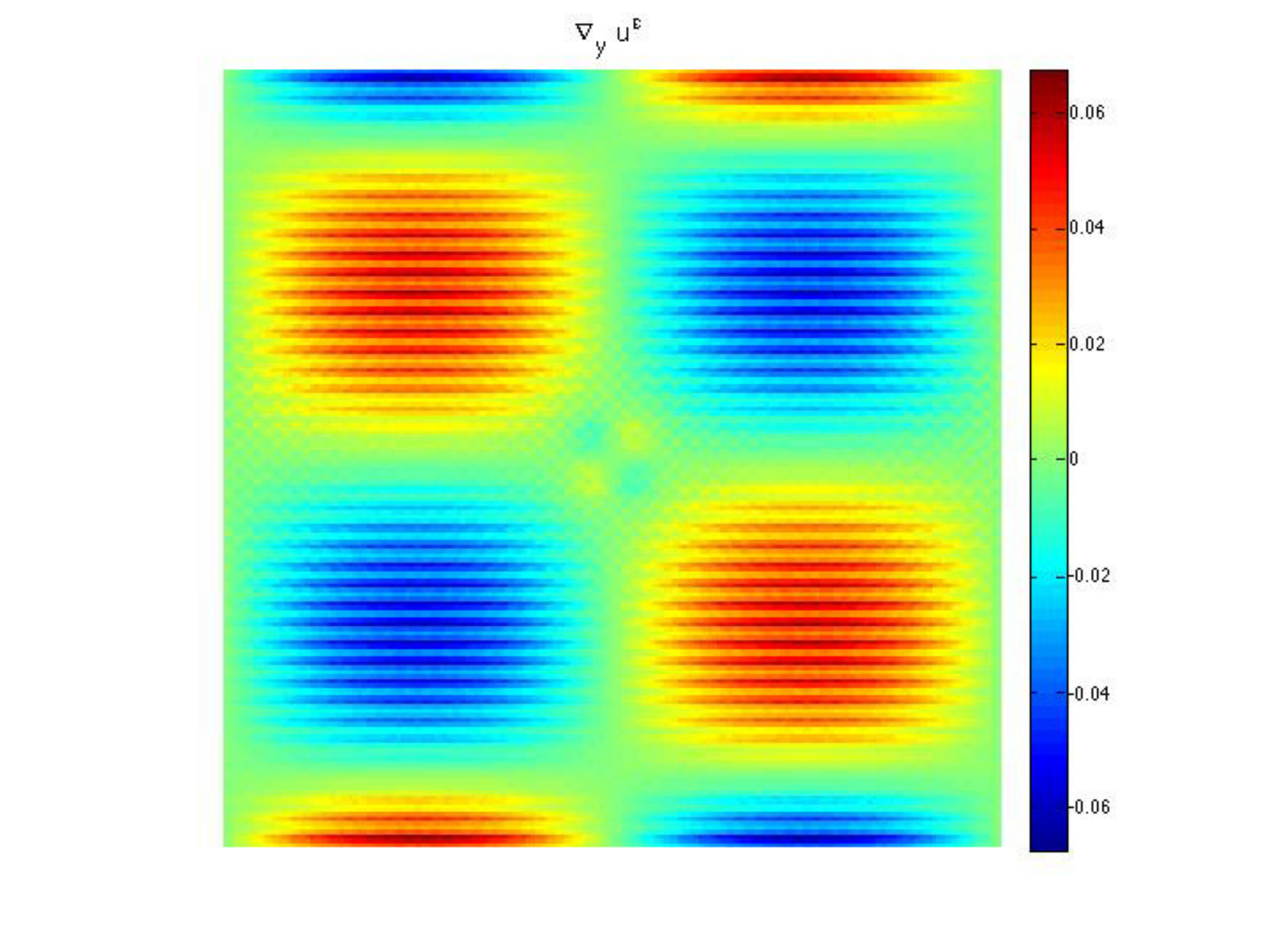}
\includegraphics[width=75mm]{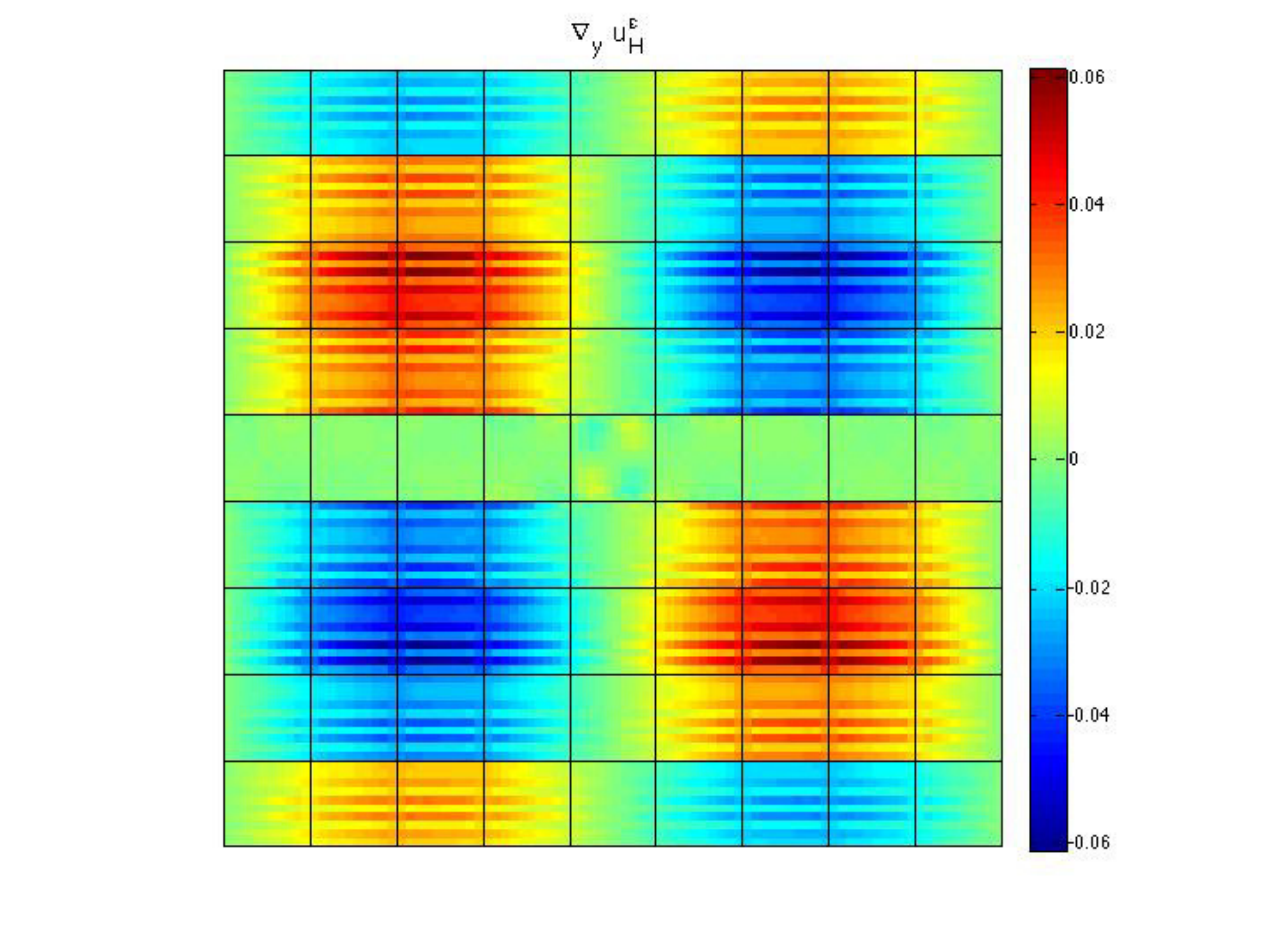}
\end{center}
\caption{Defect problem, sinusoidal loading: exact gradient $\nab u^\eps$ (left) and MsFEM gradient $\nab u^\eps_H$ (right) (top row: components $\nab u^\eps \cdot \be_1$ and $\nab u_H^\eps \cdot \be_1$; bottom row: components $\nab u^\eps \cdot \be_2$ and $\nab u_H^\eps \cdot \be_2$).}
\label{fig:solMsFEM2Ddefectgrad}
\end{figure}

\begin{figure}[H]
\begin{center}
\includegraphics[width=75mm]{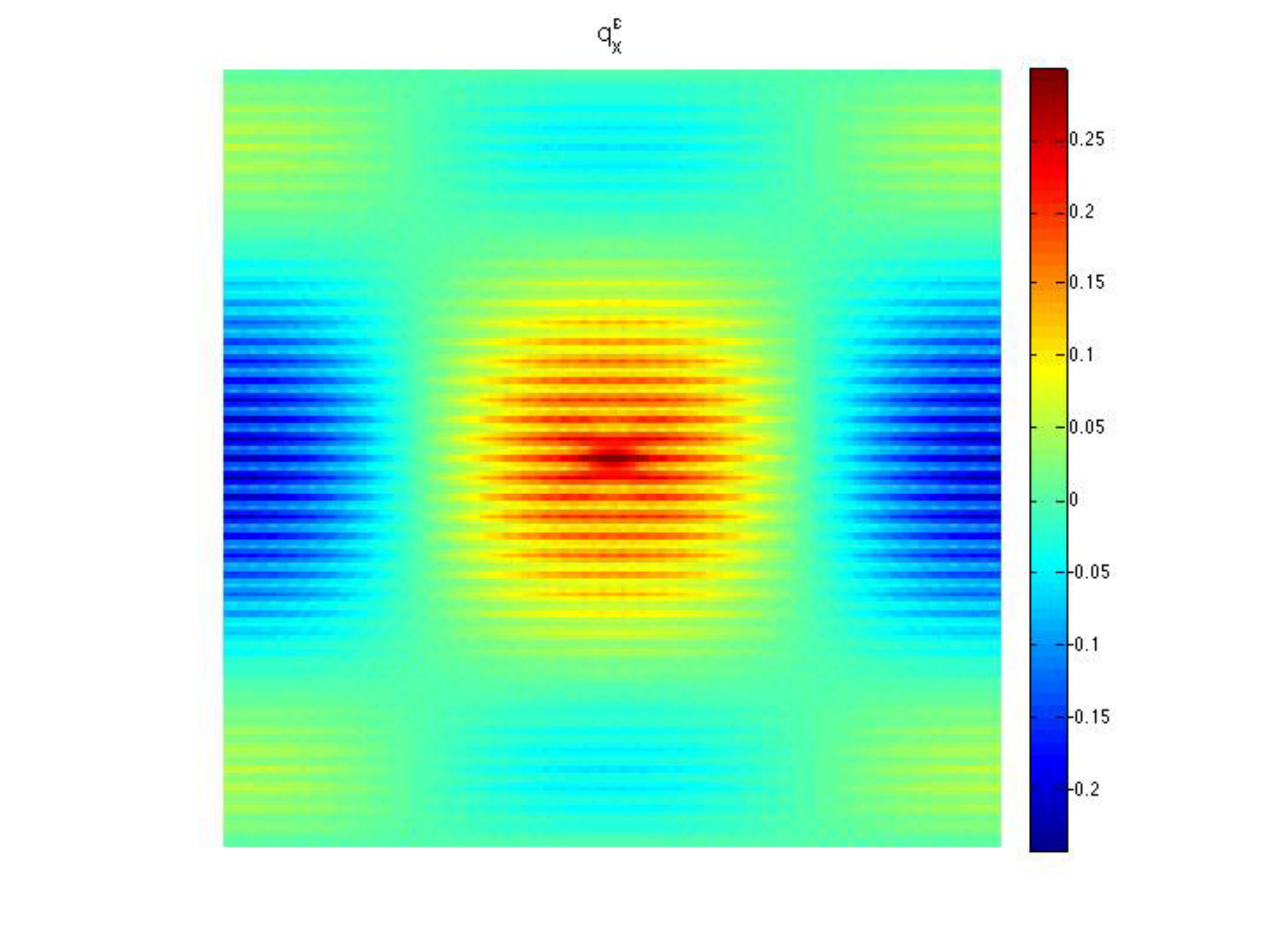}
\includegraphics[width=75mm]{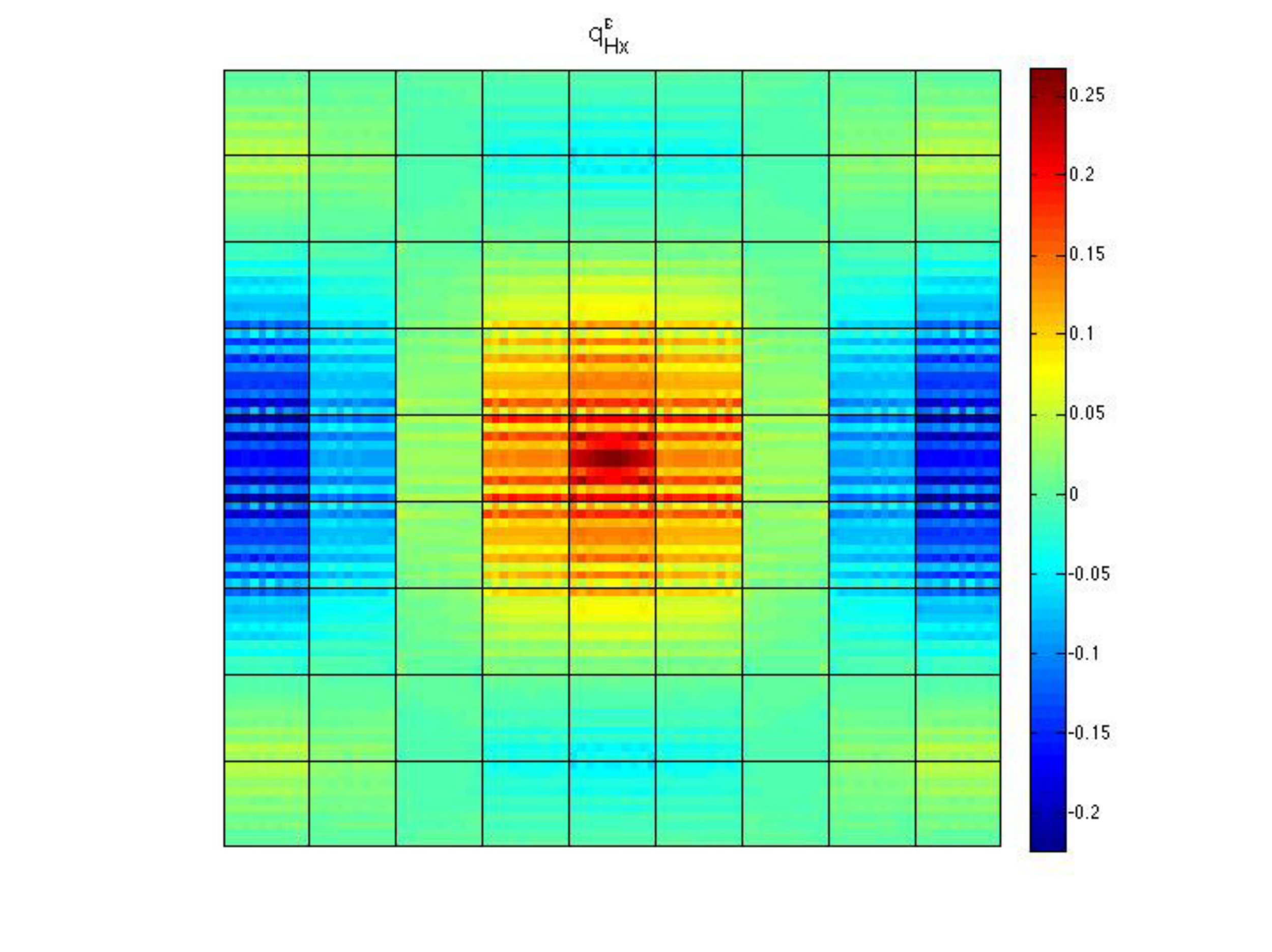}\\
\includegraphics[width=75mm]{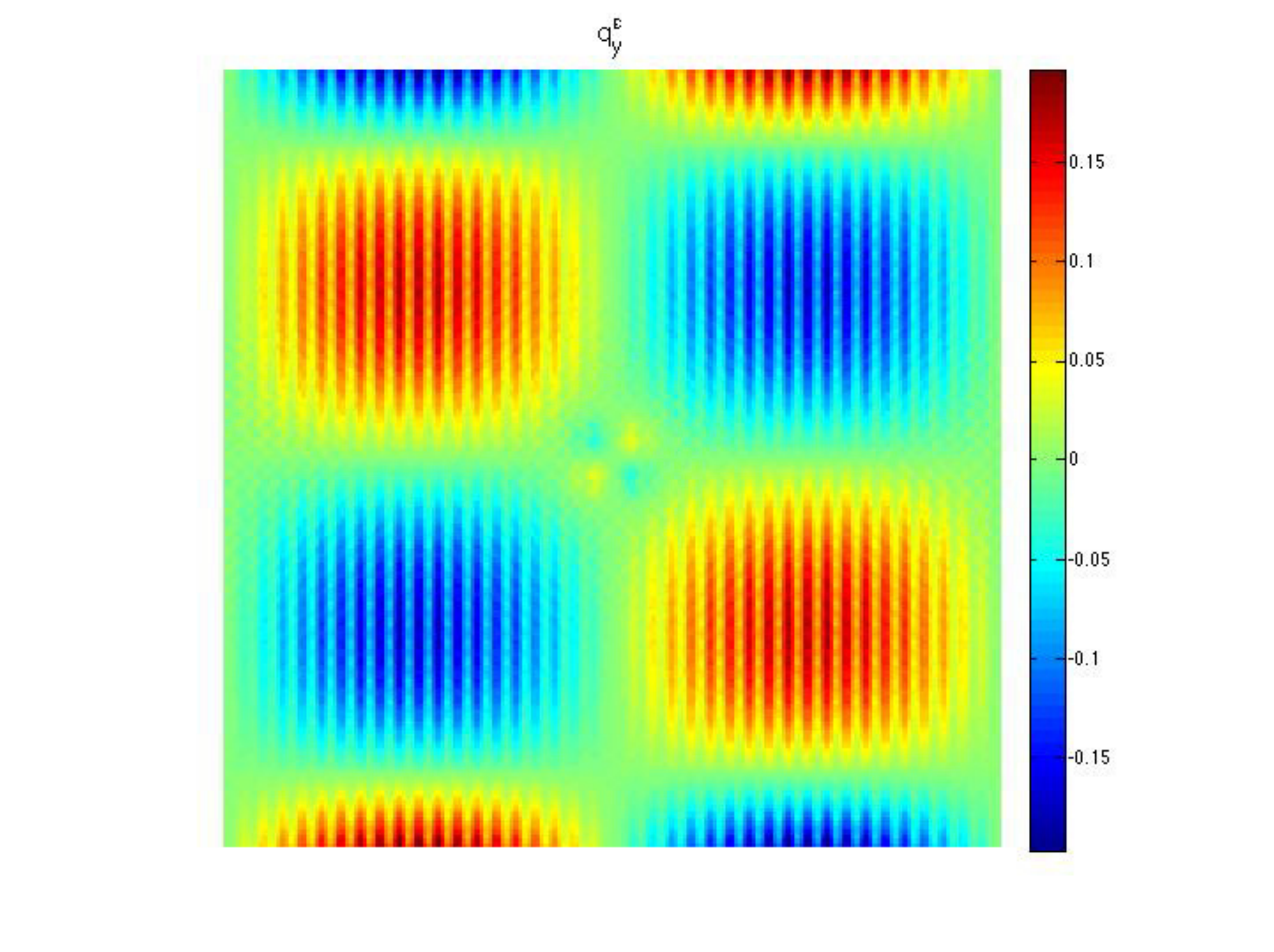}
\includegraphics[width=75mm]{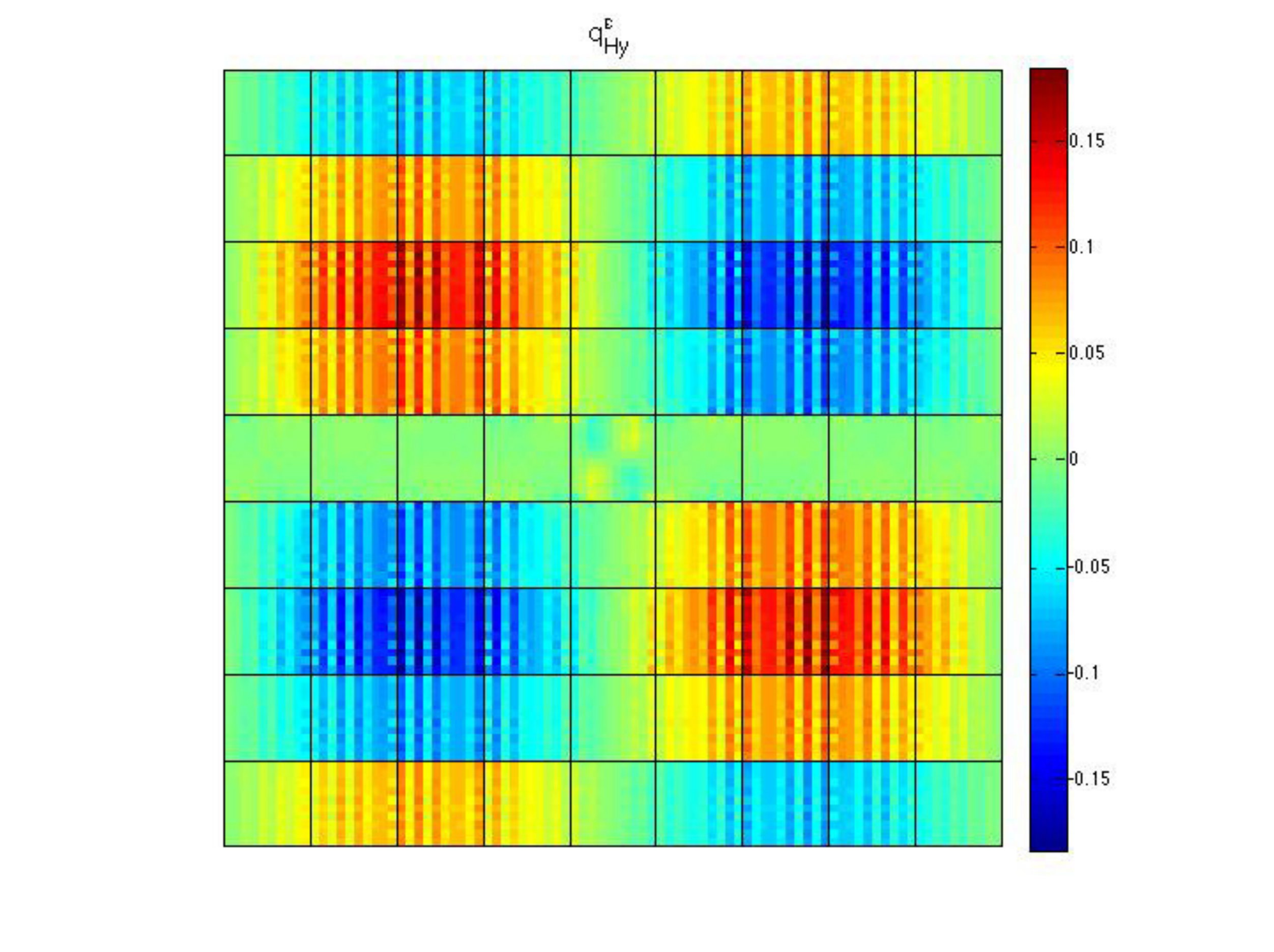}
\end{center}
\caption{Defect problem, sinusoidal loading: exact flux $\bq^\eps$ (left) and MsFEM flux $\bq^\eps_H$ (right) (top row: components $\bq^\eps \cdot \be_1$ and $\bq_H^\eps \cdot \be_1$; bottom row: components $\bq^\eps \cdot \be_2$ and $\bq_H^\eps \cdot \be_2$).}
\label{fig:solMsFEM2Ddefectflux}
\end{figure}

\begin{figure}[H]
\begin{center}
\includegraphics[width=75mm]{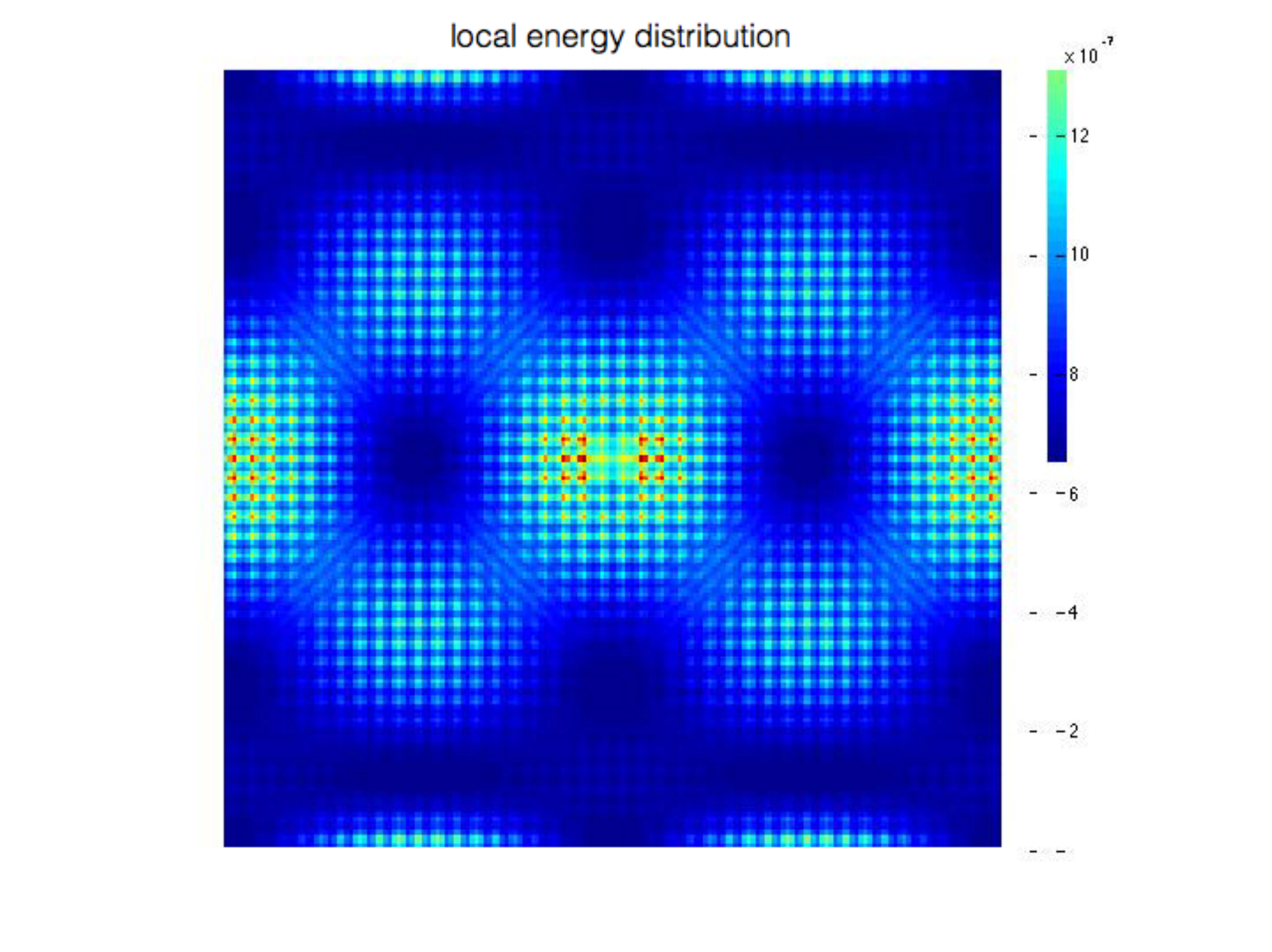}
\includegraphics[width=75mm]{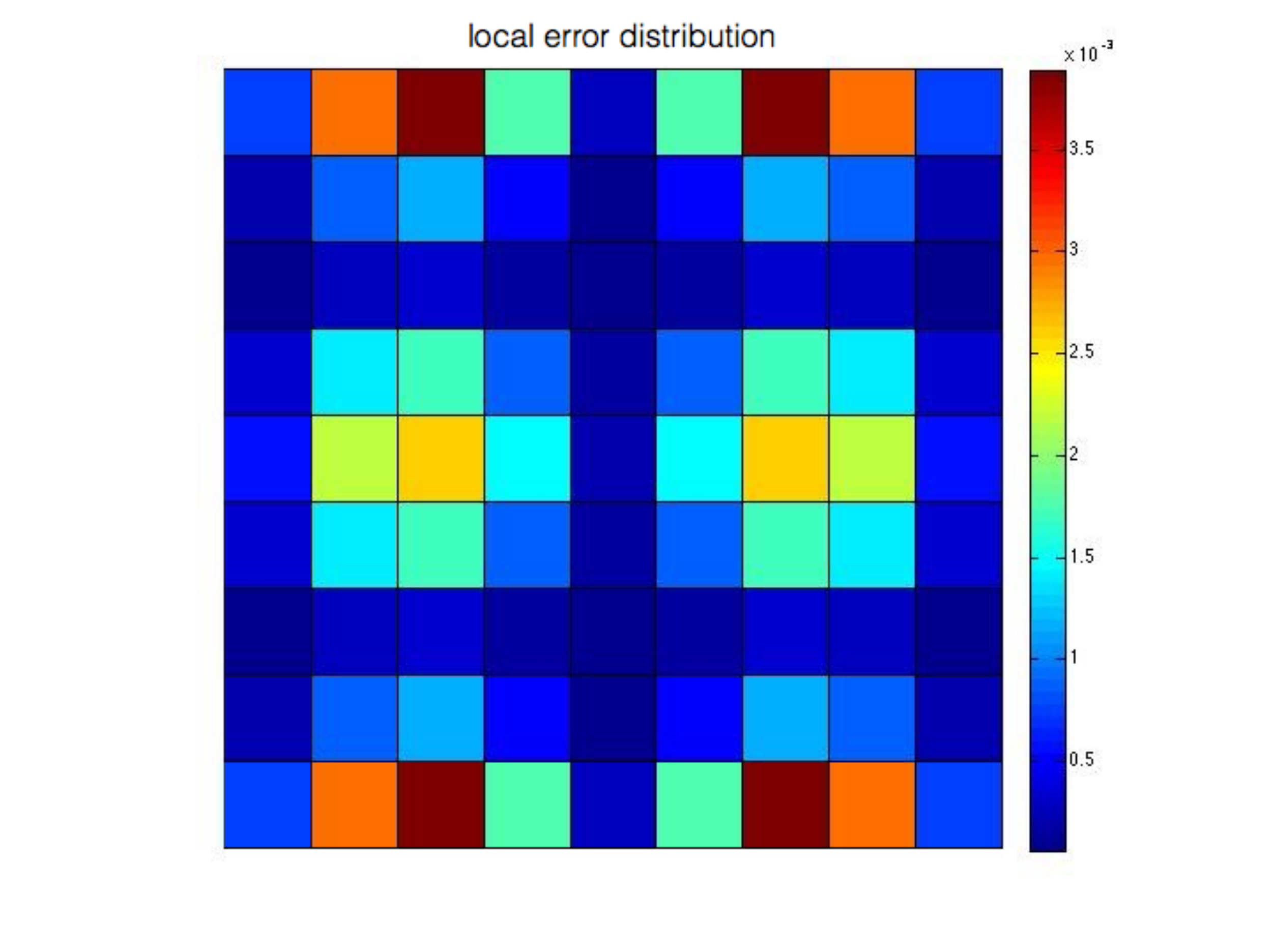}
\end{center}
\caption{Defect problem, sinusoidal loading: distribution of the energy of the exact solution (left) and of the global error (right).}
\label{fig:ener2Ddefect}
\end{figure}

We first consider the loading $f(x_1,x_2)=\sin(\pi x_1)\cos(\pi x_2)$. We take an initial mesh $\mT_H$ made of $9\times 9$ macro elements, with $h_K=\eps/3$ for any $K$ and no oversampling. The approximate MsFEM solution $u^\eps_H$ is represented below and compared to the exact solution $u^\eps$, which is obtained in practice through an expensive computation over a $500 \times 500$ fine mesh. The solutions $u^\eps$ and $u^\eps_H$ are shown on Fig.~\ref{fig:solMsFEM2Ddefectu}, the gradients $\nab u^\eps$ and $\nab u^\eps_H$ are shown on Fig.~\ref{fig:solMsFEM2Ddefectgrad}, and the fluxes $\bq^\eps$ and $\bq^\eps_H$ are shown on Fig.~\ref{fig:solMsFEM2Ddefectflux}. We also show the distribution of the energy of the exact solution and of the global error on Fig.~\ref{fig:ener2Ddefect}.

\begin{figure}[H]
\begin{center}
\includegraphics[width=50mm]{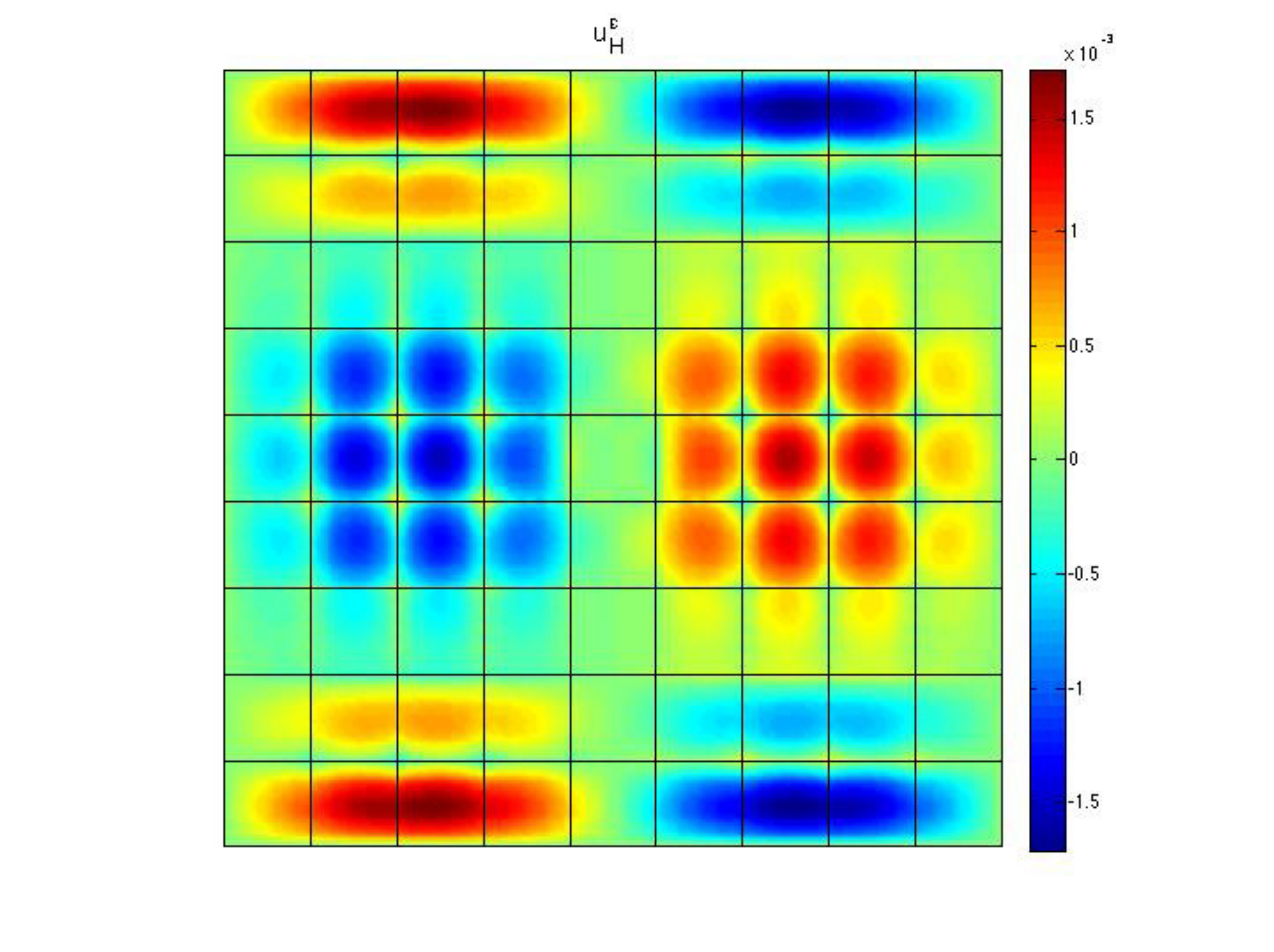}
\includegraphics[width=50mm]{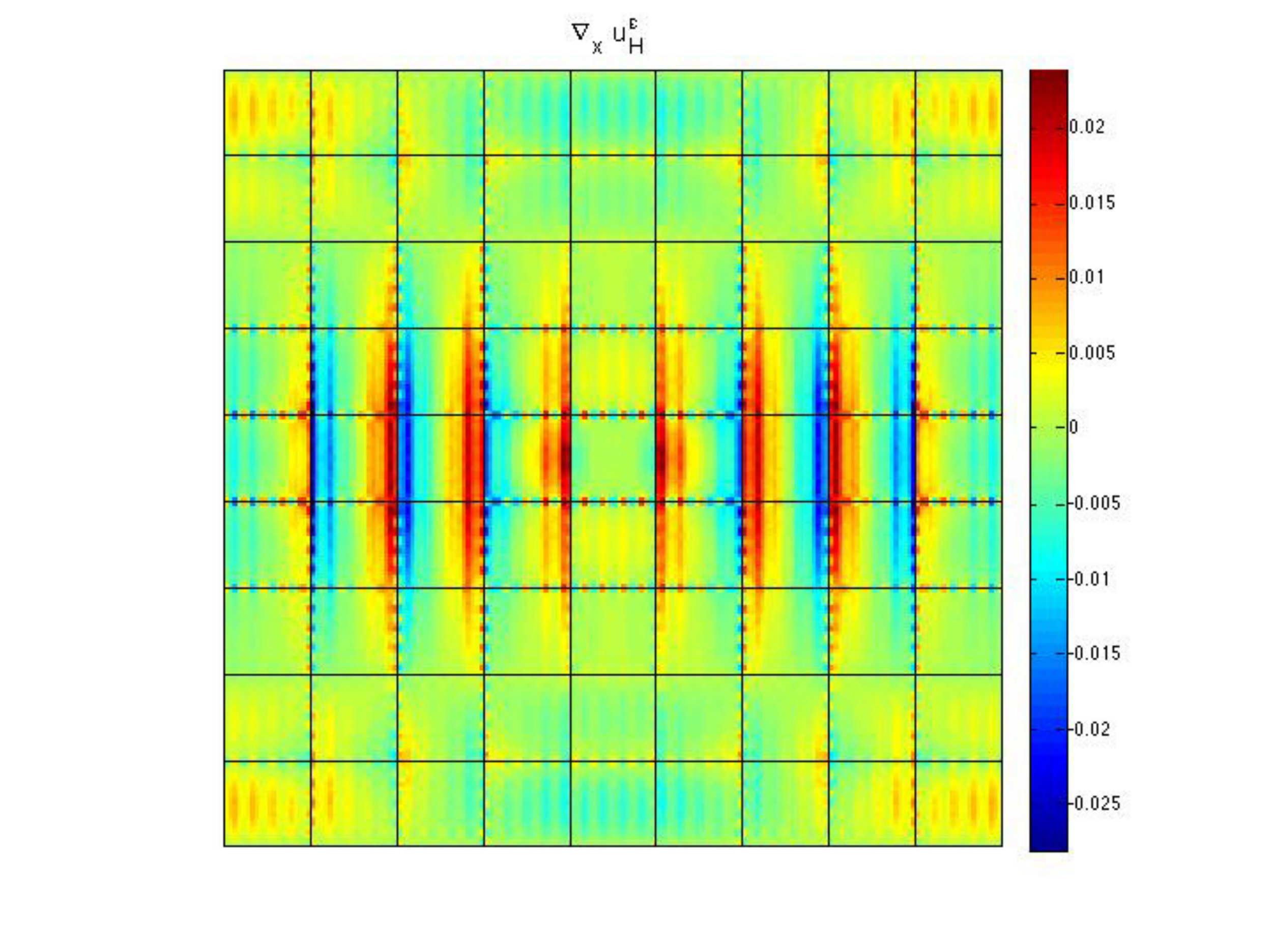}
\includegraphics[width=50mm]{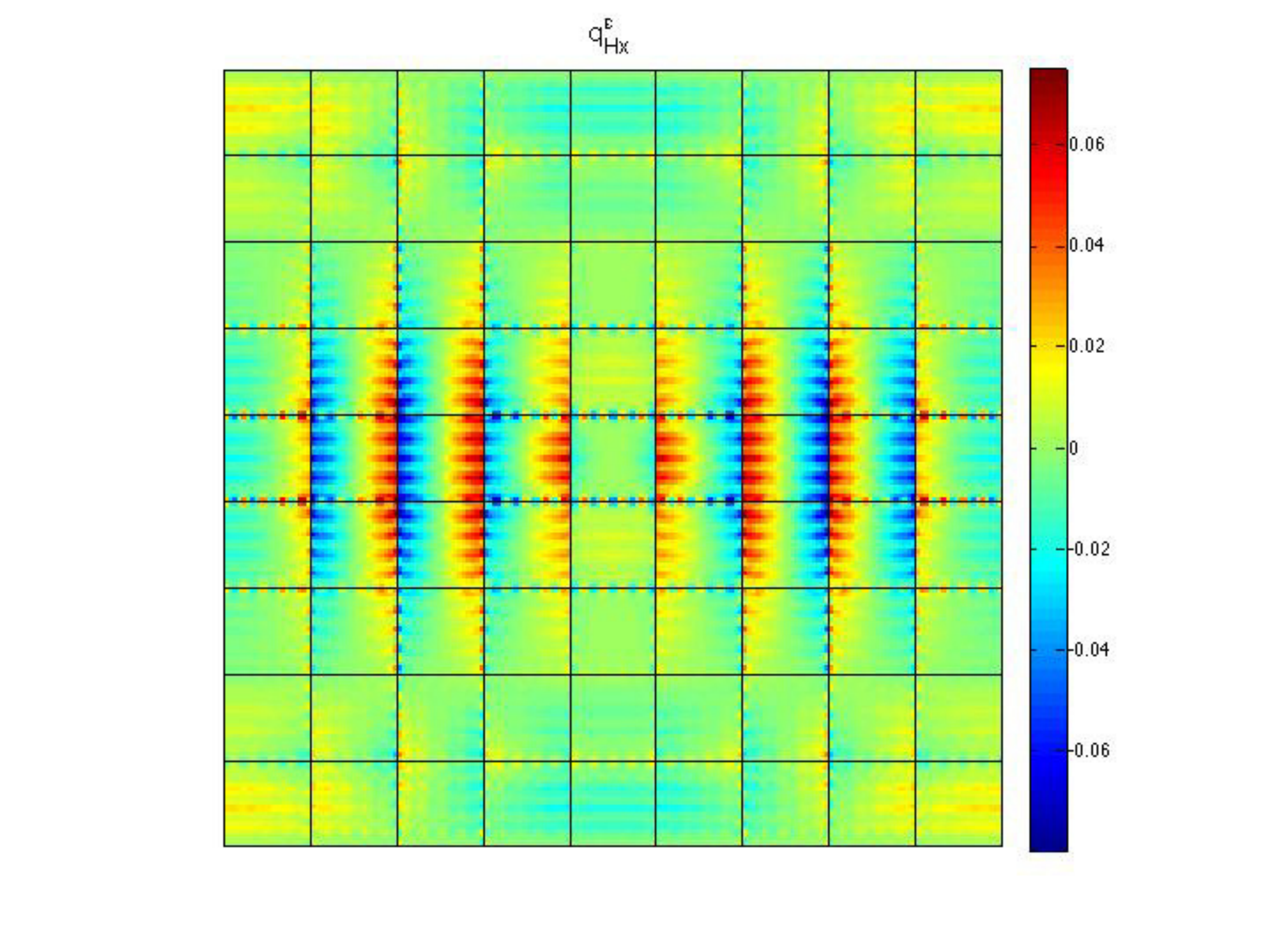}
\end{center}
\caption{Defect problem, sinusoidal loading: error between exact and MsFEM fields for $u^\eps$ (left), $\nab u^\eps \cdot \be_1$ (center) and $\bq^\eps \cdot \be_1$ (right).}
\label{fig:errorloc2Ddefectflux}
\end{figure}

We also show on Fig.~\ref{fig:errorloc2Ddefectflux} the distribution of the error on $u^\eps$, $\nab u^\eps \cdot \be_1$ and $\bq^\eps \cdot \be_1$.

\begin{figure}[H]
\begin{center}
\includegraphics[width=75mm]{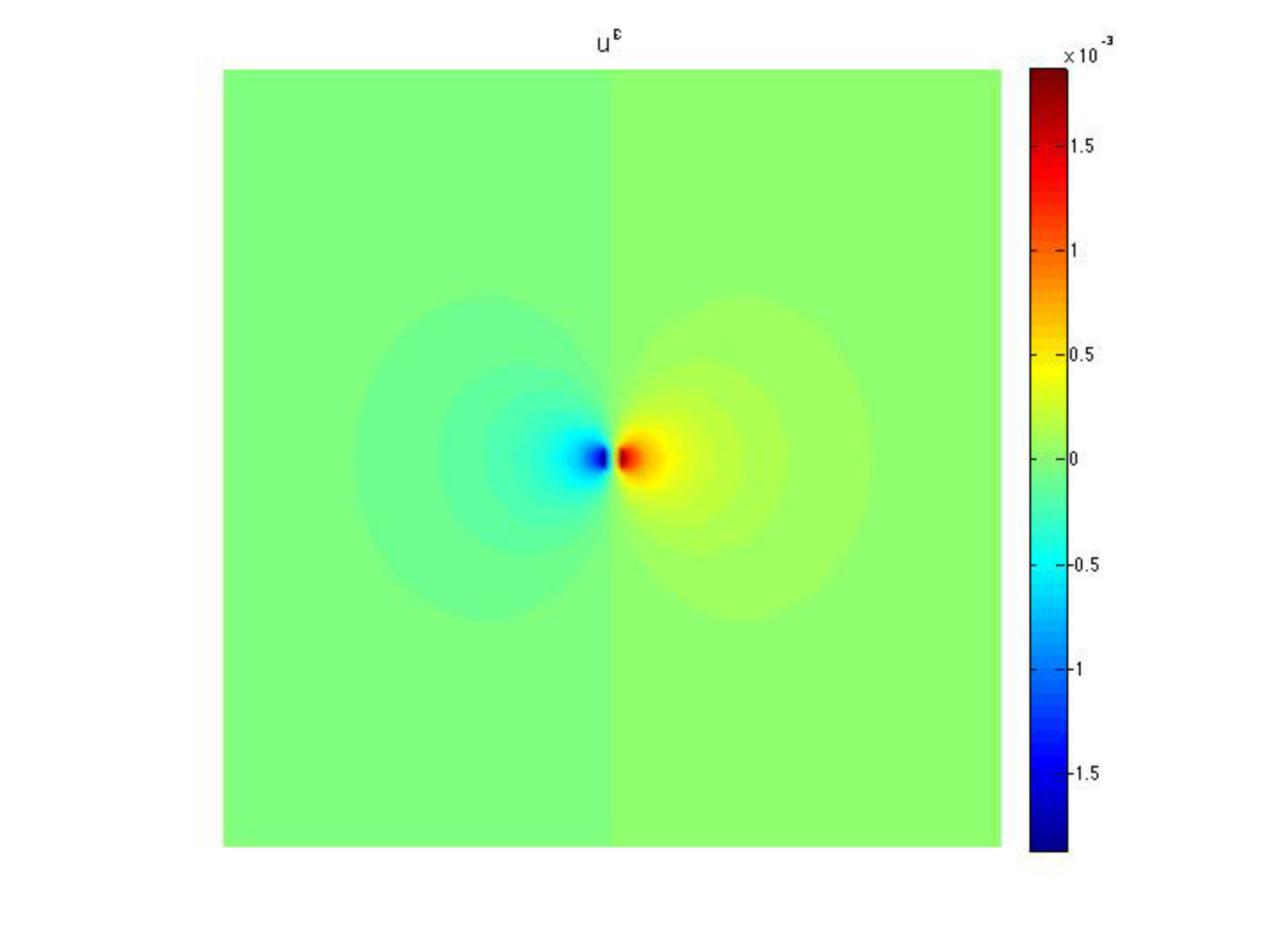}
\includegraphics[width=75mm]{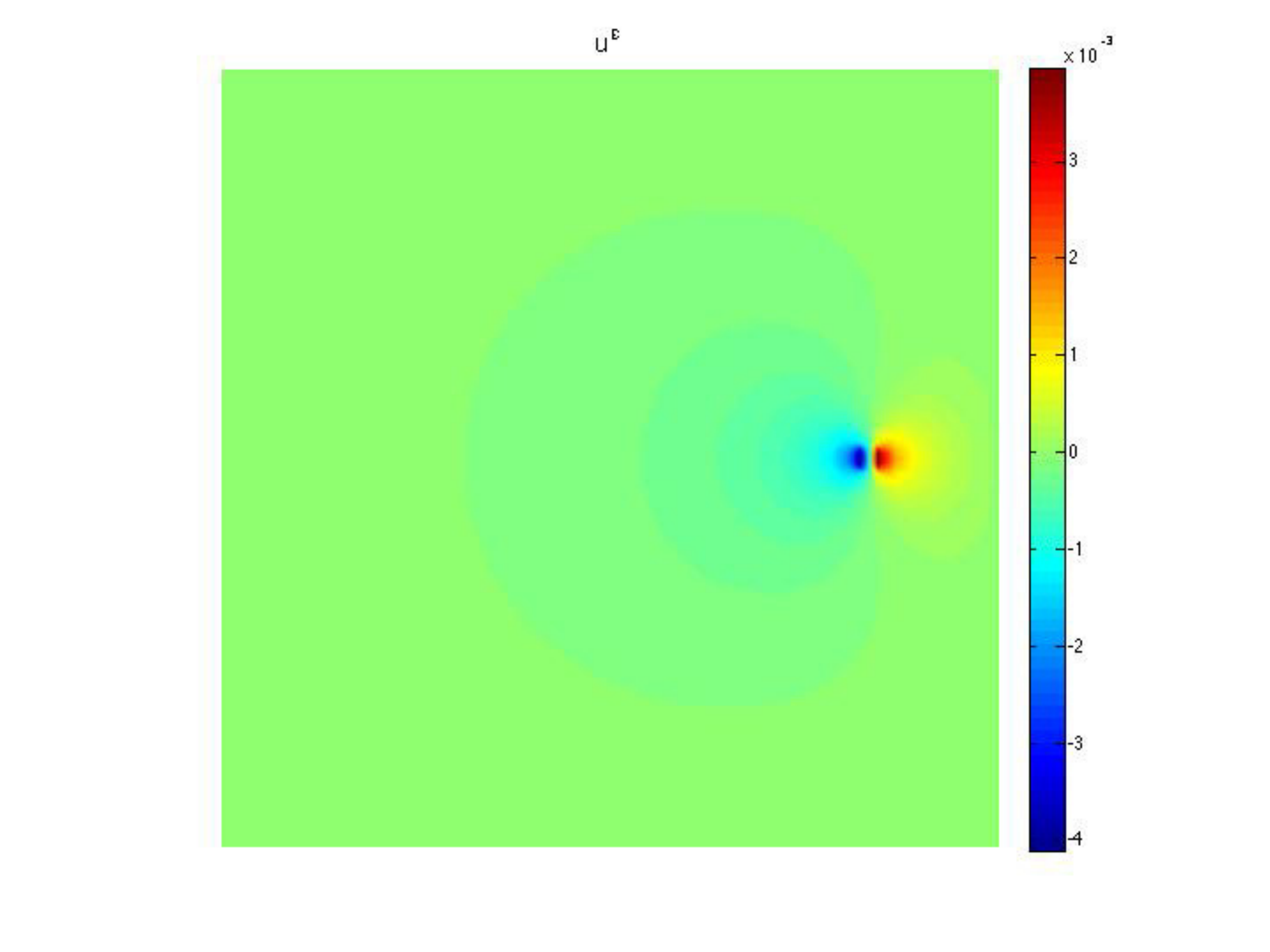}
\end{center}
\caption{Defect problem, sinusoidal loading: exact adjoint solution $\widetilde{u}^\eps$ for $Q_1$ (left) and $Q_2$ (right).}
\label{fig:soladj2Ddefectu}
\end{figure}

\begin{figure}[H]
\begin{center}
\includegraphics[width=75mm]{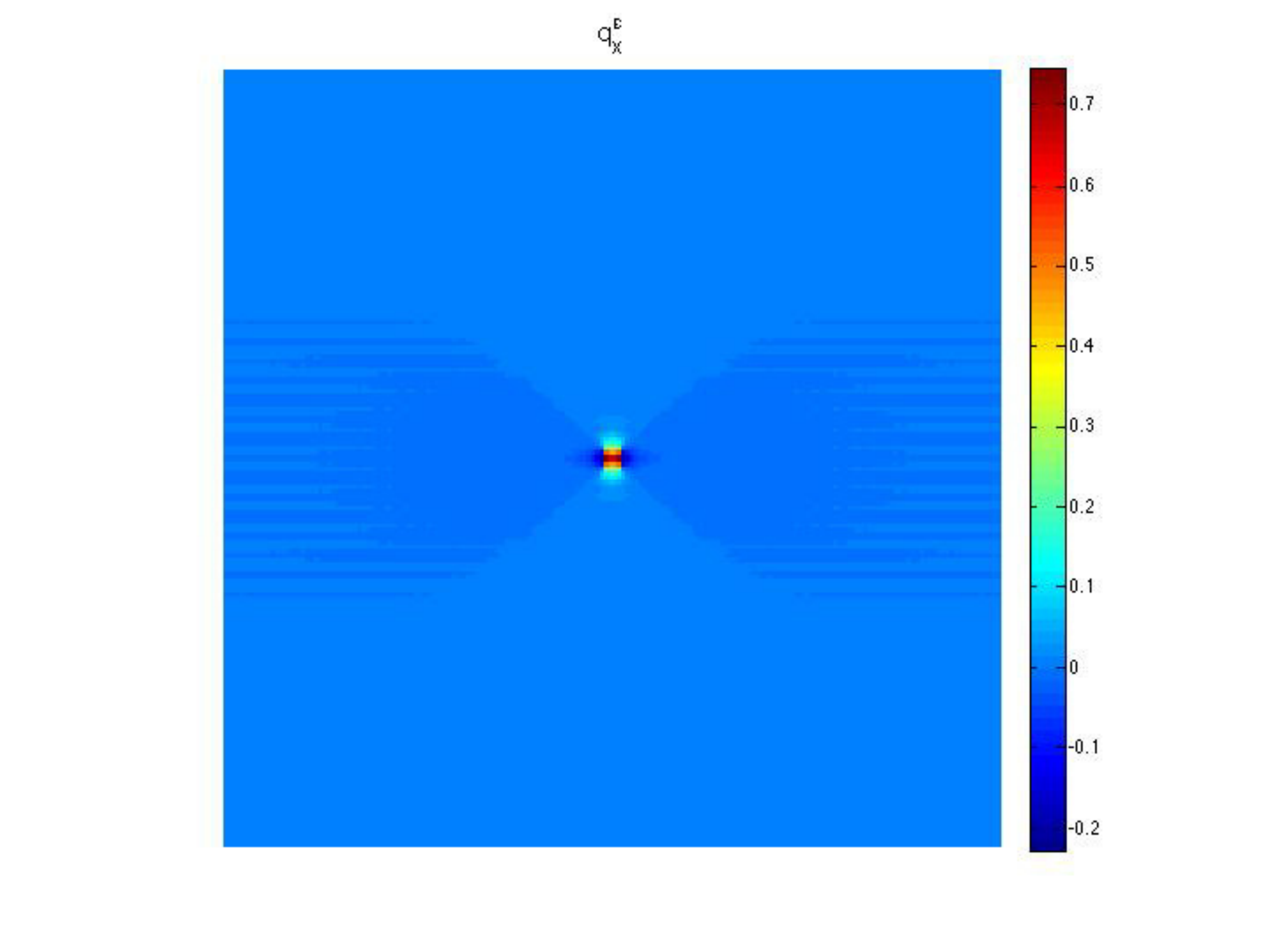}
\includegraphics[width=75mm]{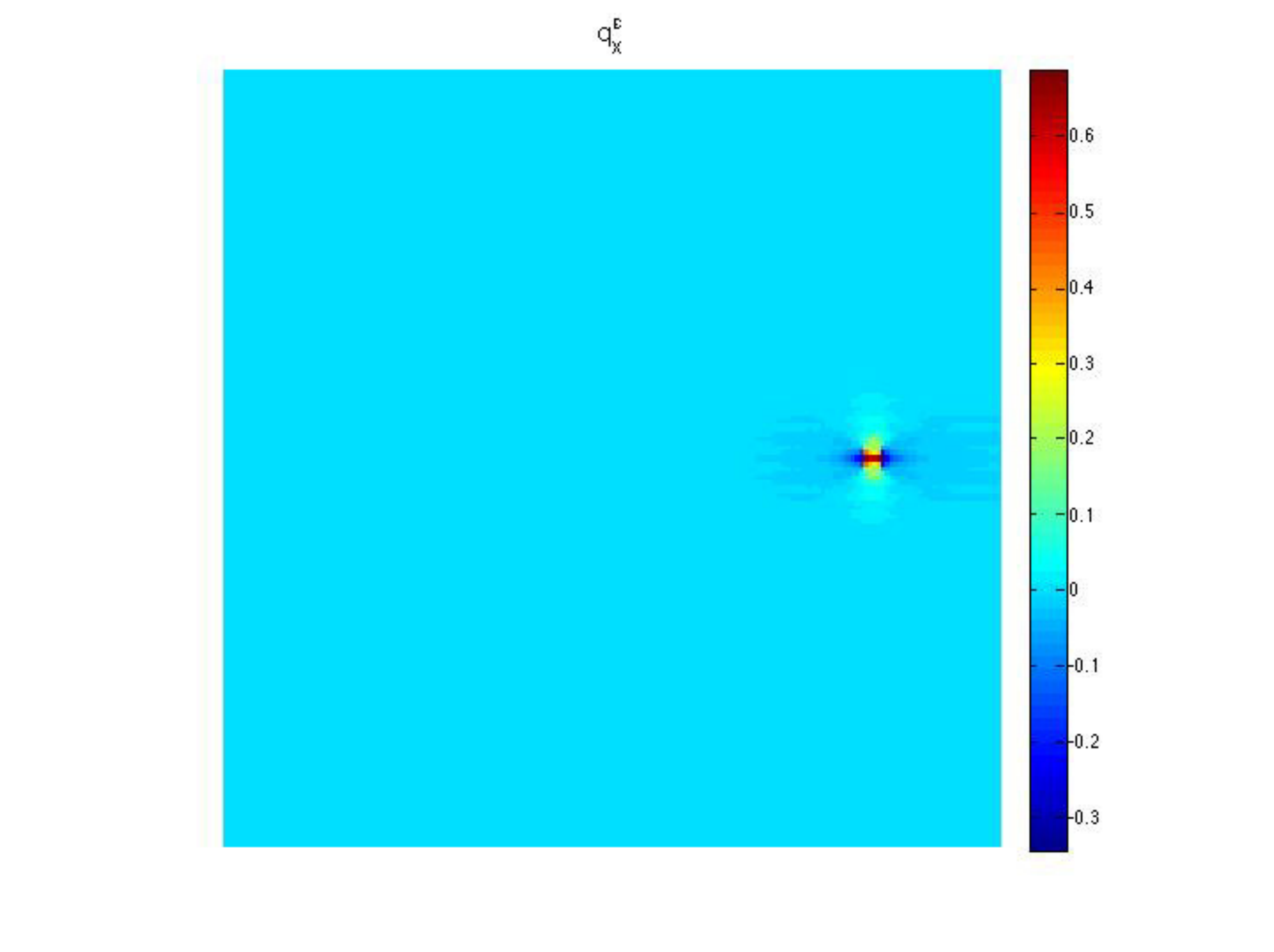}\\
\includegraphics[width=75mm]{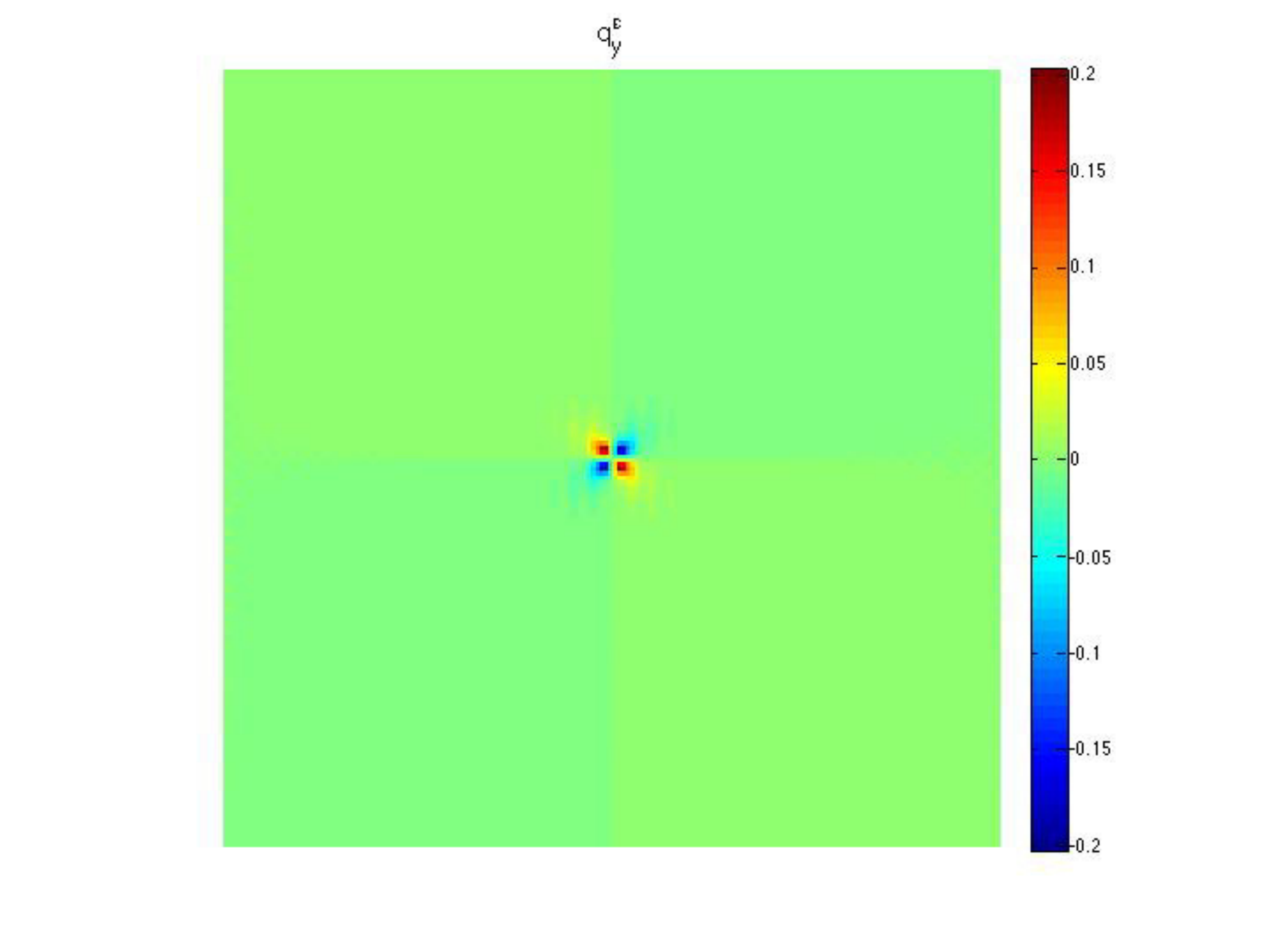}
\includegraphics[width=75mm]{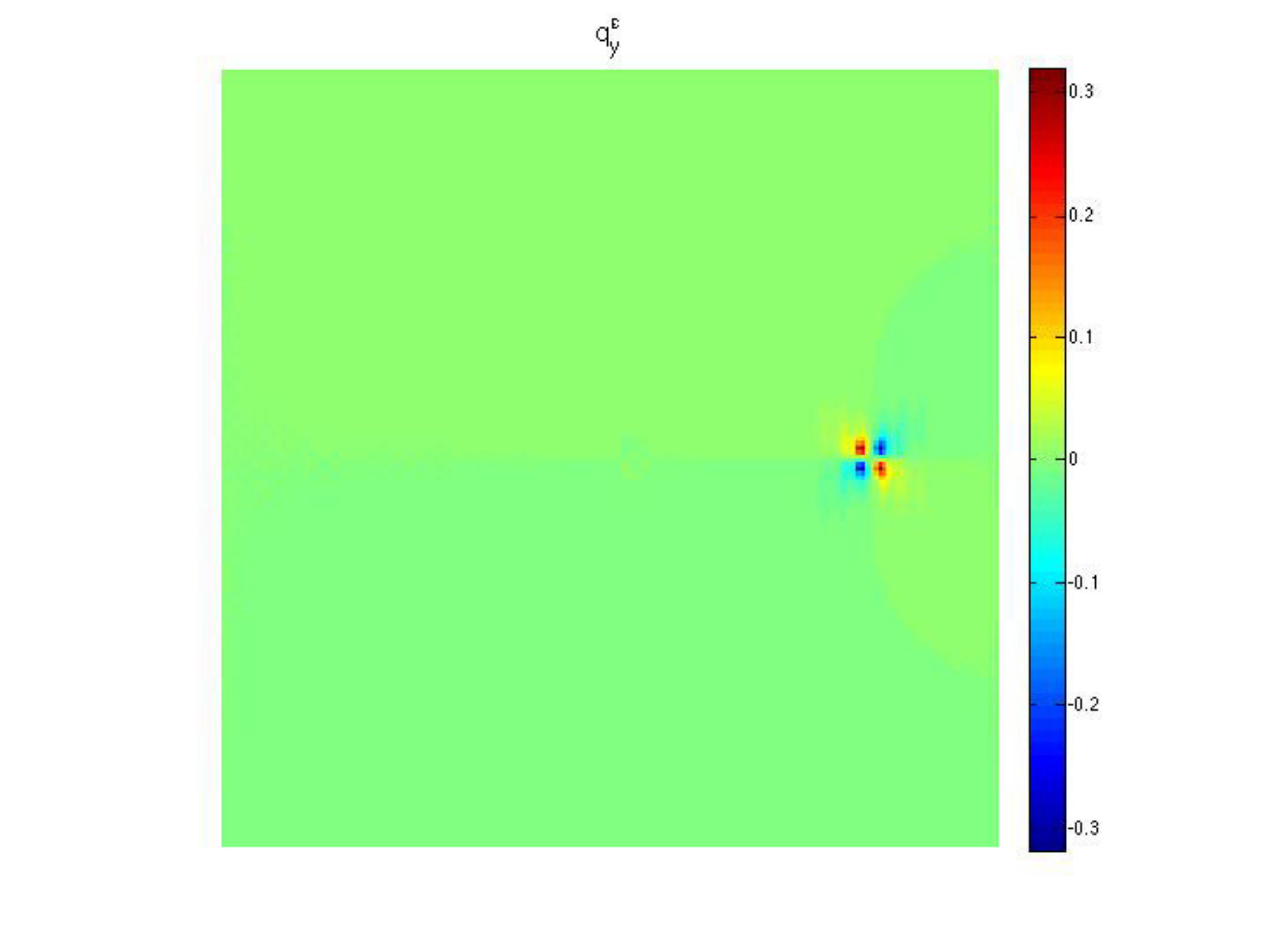}
\end{center}
\caption{Defect problem, sinusoidal loading: exact adjoint flux $\widetilde{\bq}^\eps$ for $Q_1$ (left) and $Q_2$ (right) (top row: component along $\be_1$; bottom row: component along $\be_2$).}
\label{fig:soladj2Ddefectq}
\end{figure}

We choose the quantity of interest $\dis Q(u^\eps) = \frac{1}{|\omega|} \int_{\omega} \bq^\eps \cdot \be_1 = \frac{1}{|\omega|} \int_{\omega} \be_1 \cdot \Aa^\eps \nab u^\eps$ with $|\omega| = 4 \eps \times 4 \eps$. We consider two cases:
\begin{itemize}
\item Case 1: the subregion $\omega$ is centered on the defect ($Q=Q_1$);
\item Case 2: the subregion $\omega$ is far from the defect ($Q=Q_2$).
\end{itemize}
Note that $u^\eps$ is close to zero in the region around the defect (see Fig.~\ref{fig:solMsFEM2Ddefectu}), so considering as quantity of interest the average of $u^\eps$ would not be interesting. In contrast, $\bq^\eps \cdot \be_1$ is not close to zero in that region (see Fig.~\ref{fig:solMsFEM2Ddefectflux}). The adjoint solutions $\widetilde{u}^\eps$ and fluxes $\widetilde{\bq}^\eps$ corresponding to the quantities of interest $Q_1$ and $Q_2$ are shown on Figs.~\ref{fig:soladj2Ddefectu} and~\ref{fig:soladj2Ddefectq}. We observe, as expected, that the adjoint solutions are essentially supported in the domain of interest $\omega$. 

\begin{figure}[H]
\begin{center}
\includegraphics[width=150mm]{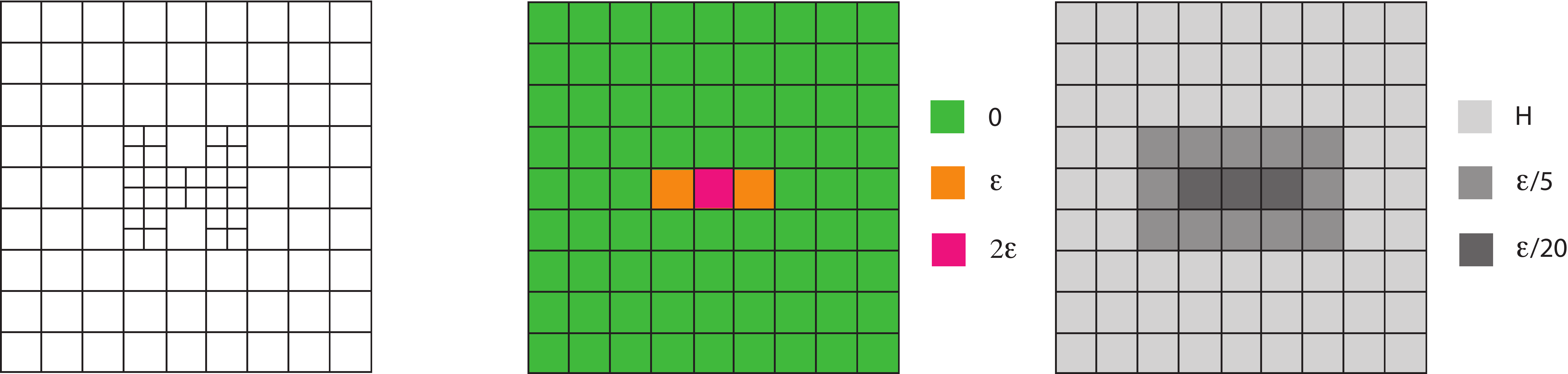} 
\end{center}
\caption{Defect problem, sinusoidal loading: final MsFEM discretization for $Q_1$ (left: coarse mesh; center: oversampling size; right: fine mesh size $h_K$).}
\label{fig:meshconvQ12Ddefect}
\end{figure}

\begin{figure}[H]
\begin{center}
\includegraphics[width=150mm]{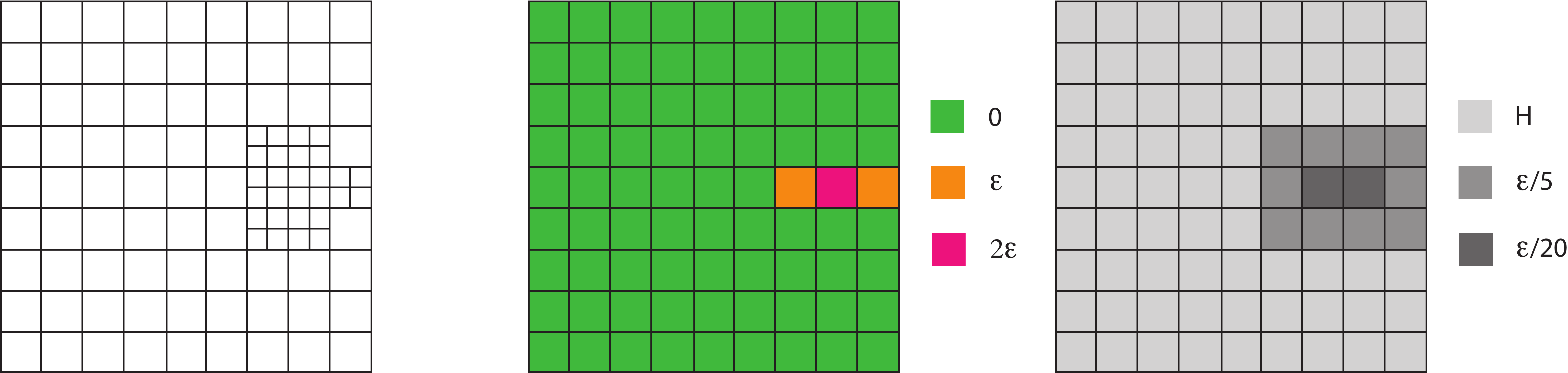}
\end{center}
\caption{Defect problem, sinusoidal loading: final MsFEM discretization for $Q_2$.}
\label{fig:meshconvQ22Ddefect}
\end{figure}

\begin{figure}[H]
\begin{center}
\includegraphics[width=150mm]{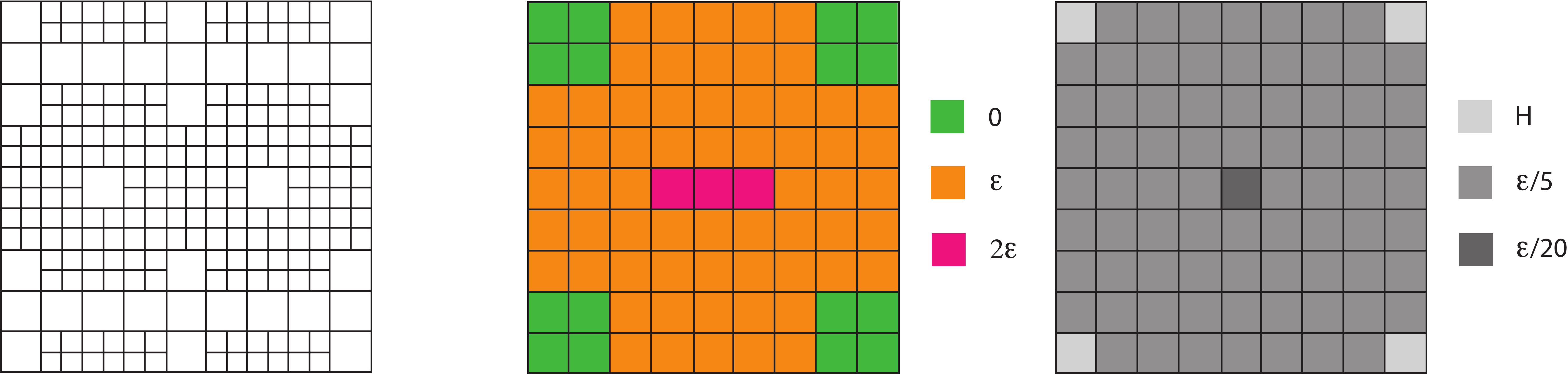} 
\end{center}
\caption{Defect problem, sinusoidal loading: Final MsFEM discretization when controlling the global error.}
\label{fig:meshconvglob2Ddefect}
\end{figure}

The error threshold is fixed to $1\%$. We show on Fig.~\ref{fig:meshconvQ12Ddefect} (resp. Fig.~\ref{fig:meshconvQ22Ddefect}) the final MsFEM discretization obtained after using the adaptive algorithm presented in Section~\ref{section:algorithm}, in order to respect the error threshold for the quantity of interest $Q_1$ (resp. $Q_2$). This discretization is to be compared with that obtained when controlling the global error (in energy norm) using the tools presented in~\cite{CHA16b}, which is shown on Fig.~\ref{fig:meshconvglob2Ddefect}. Obviously, the mesh on Figs.~\ref{fig:meshconvQ12Ddefect} and~\ref{fig:meshconvQ22Ddefect} is much coarser than the mesh on Fig.~\ref{fig:meshconvglob2Ddefect}. This shows, as expected, that it is worth working with a goal-oriented adaptive strategy, rather than an adaptive strategy based on the global error: for the same error, the computational cost is much lower. Comparing Figs.~\ref{fig:meshconvQ12Ddefect} and~\ref{fig:meshconvQ22Ddefect}, we also observe that the mesh is only refined close to the region of interest $\omega$. In the second example, when this region is far from the defect, it is not needed to adapt the mesh close to the defect.

As regards the accuracy of the error estimate (with respect to the actual error), we observe that $\eta^Q/|Q(u^\eps)-Q(u^\eps_H)| = 1.11$ (resp. 1.13) at the first iteration of the adaptive algorithm, and $\eta^Q/|Q(u^\eps)-Q(u^\eps_H)| = 1.07$ (resp. 1.06) at the final iteration of the adaptive algorithm, for the quantities of interest $Q_1$ and $Q_2$, respectively. As pointed out above, the error estimate thus only very slightly overestimates the actual error. 

\begin{figure}[H]
\begin{center}
\includegraphics[width=75mm]{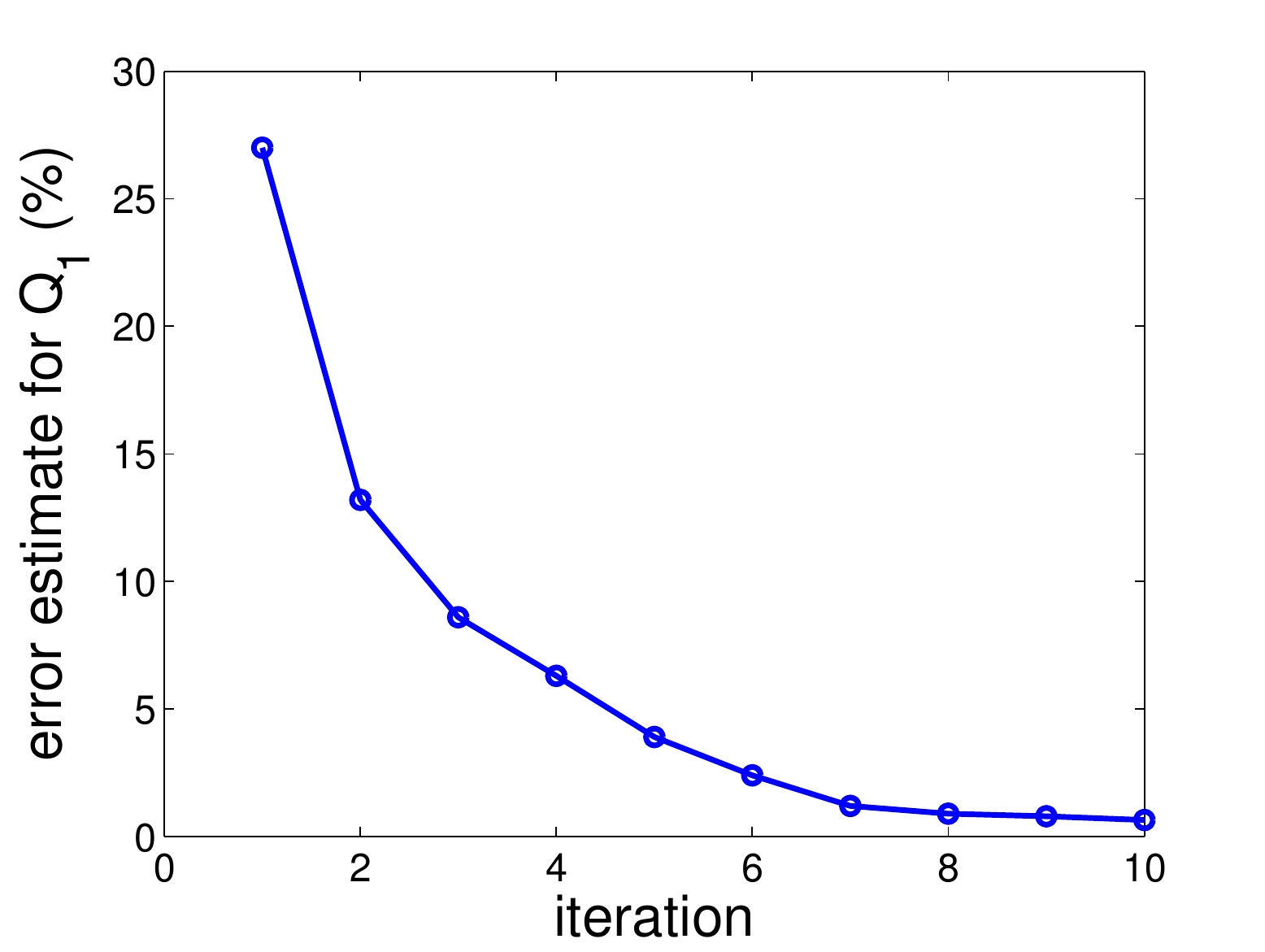} 
\includegraphics[width=75mm]{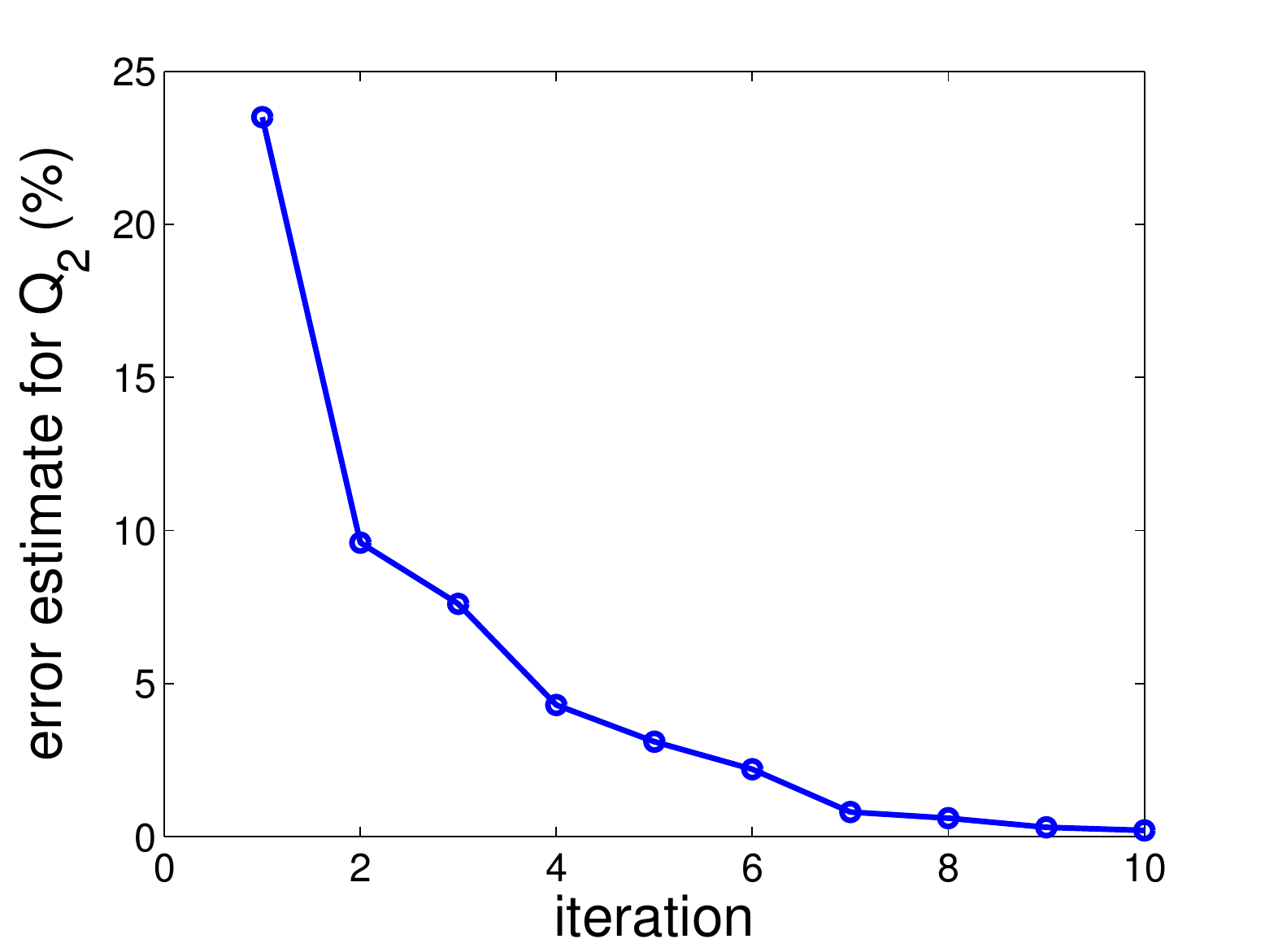} 
\end{center}
\caption{Defect problem, sinusoidal loading: convergence of the error estimate for $Q_1$ (left) and $Q_2$ (right).}
\label{fig:converrorQ2Ddefect}
\end{figure}

\begin{figure}[H]
\begin{center}
\includegraphics[width=75mm]{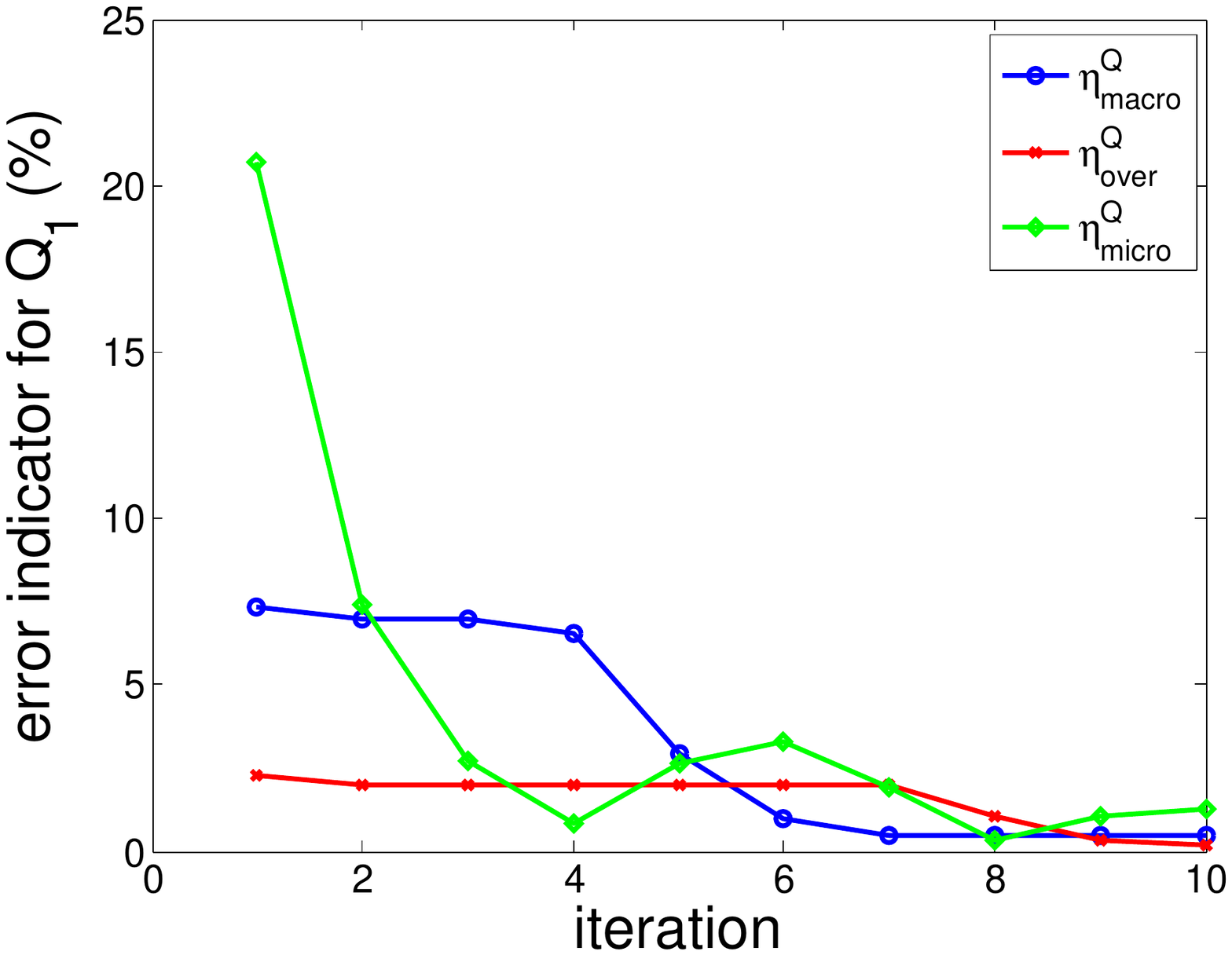} 
\includegraphics[width=75mm]{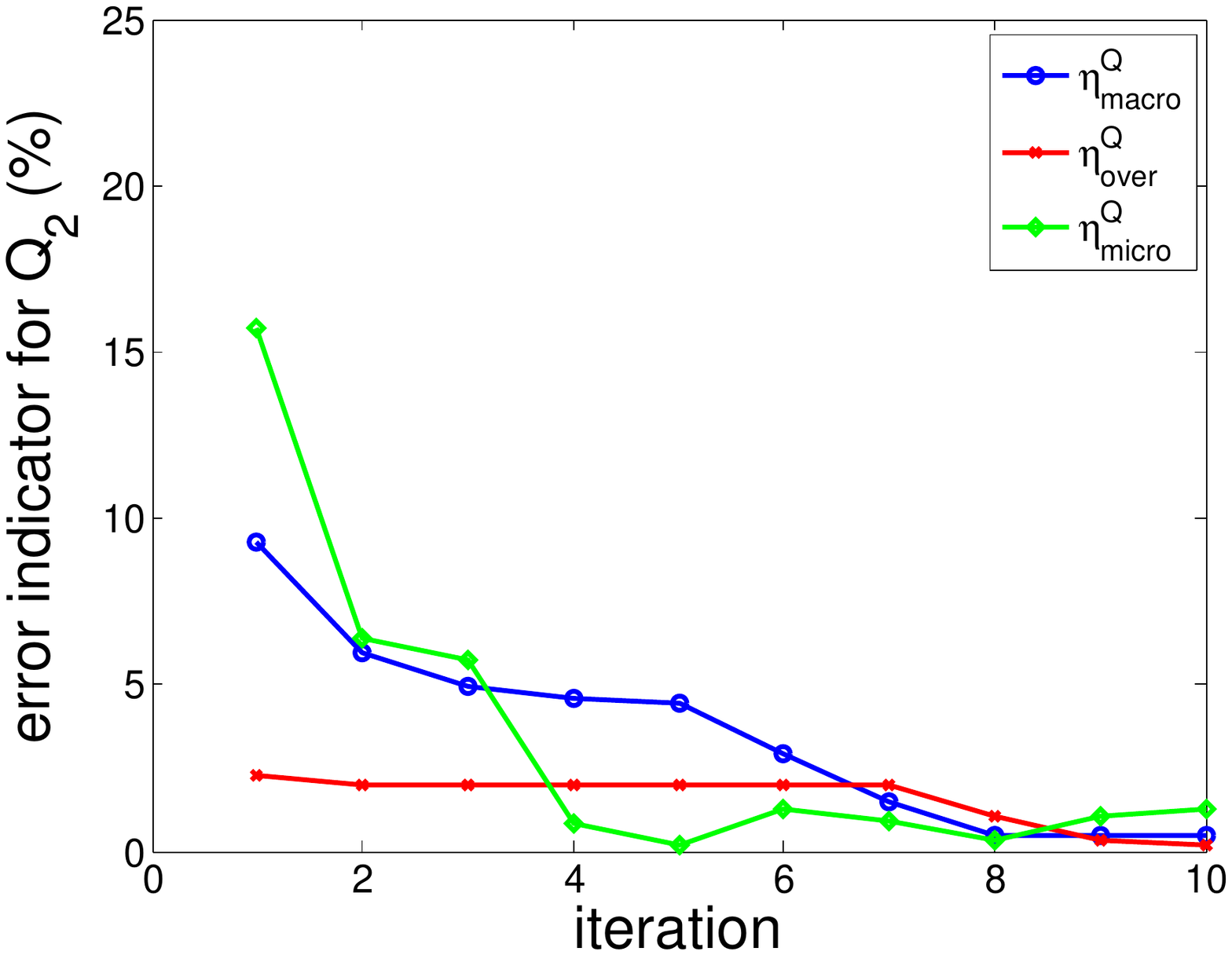}
\end{center}
\caption{Defect problem, sinusoidal loading: convergence of the error indicators for $Q_1$ (left) and $Q_2$ (right).}
\label{fig:convindicQ2Ddefect}
\end{figure}

We eventually show on Fig.~\ref{fig:converrorQ2Ddefect} the convergence of the error estimate on $Q$ with respect to the iterations of the adaptive process. As the iterations proceed, the error (which is initially close to 25\%, due purposely to the coarse mesh that is used) monotically decreases until reaching the prescribed tolerance of 1\% in about 10 iterations. We also show on Fig.~\ref{fig:convindicQ2Ddefect} the evolution of the error indicators along the adaptive process. We observe that, for this case, the fine mesh is first refined before adapting the coarse mesh and finally extending the oversampling zone.

\subsubsection{Exponential loading}

\begin{figure}[H]
\begin{center}
\includegraphics[width=75mm]{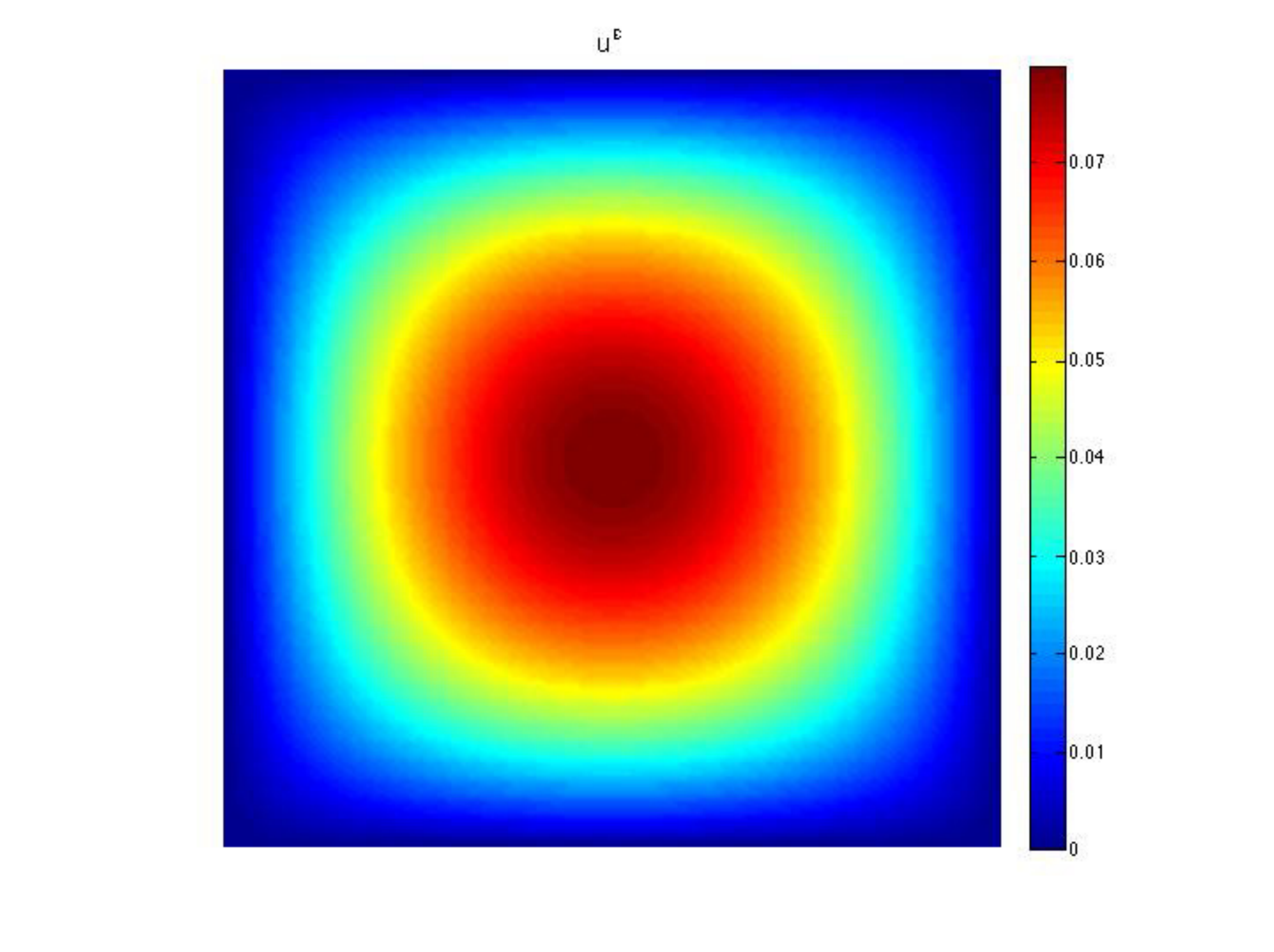}
\includegraphics[width=75mm]{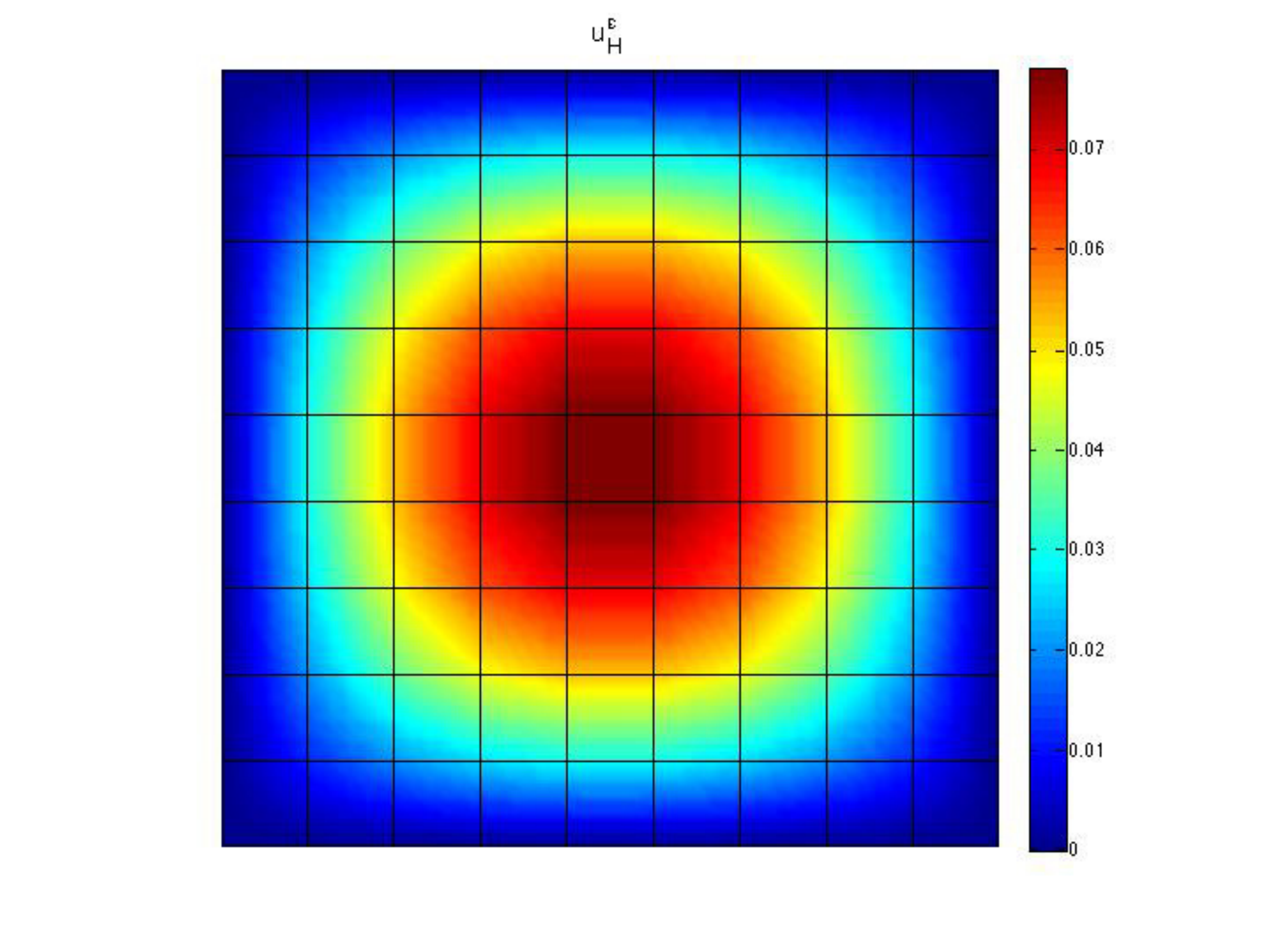}
\end{center}
\caption{Defect problem, exponential loading: exact solution $u^\eps$ (left) and MsFEM solution $u^\eps_H$ (right).}
\label{fig:solMsFEM2Ddefect2u}
\end{figure}

\begin{figure}[H]
\begin{center}
\includegraphics[width=75mm]{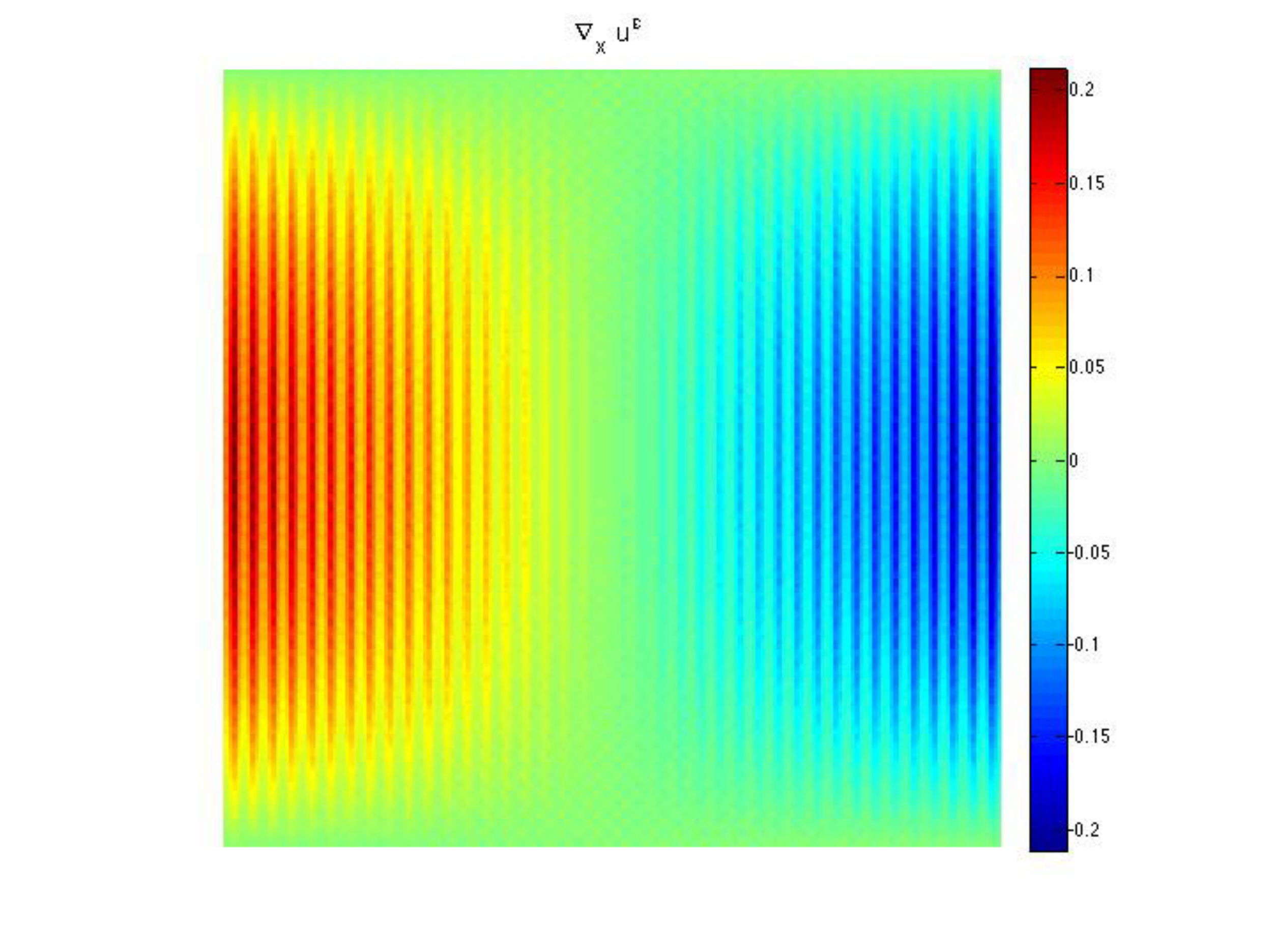}
\includegraphics[width=75mm]{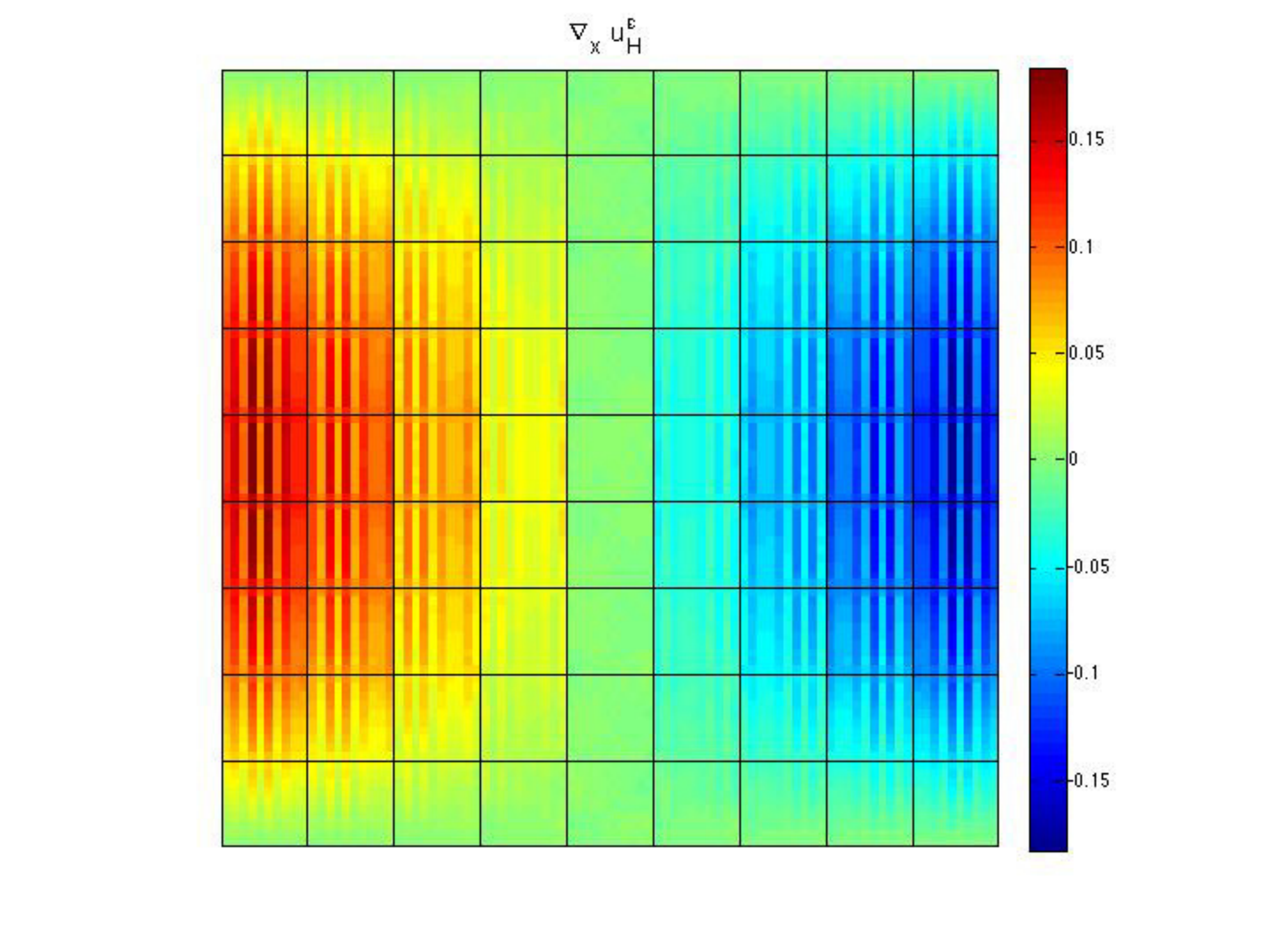}\\
\includegraphics[width=75mm]{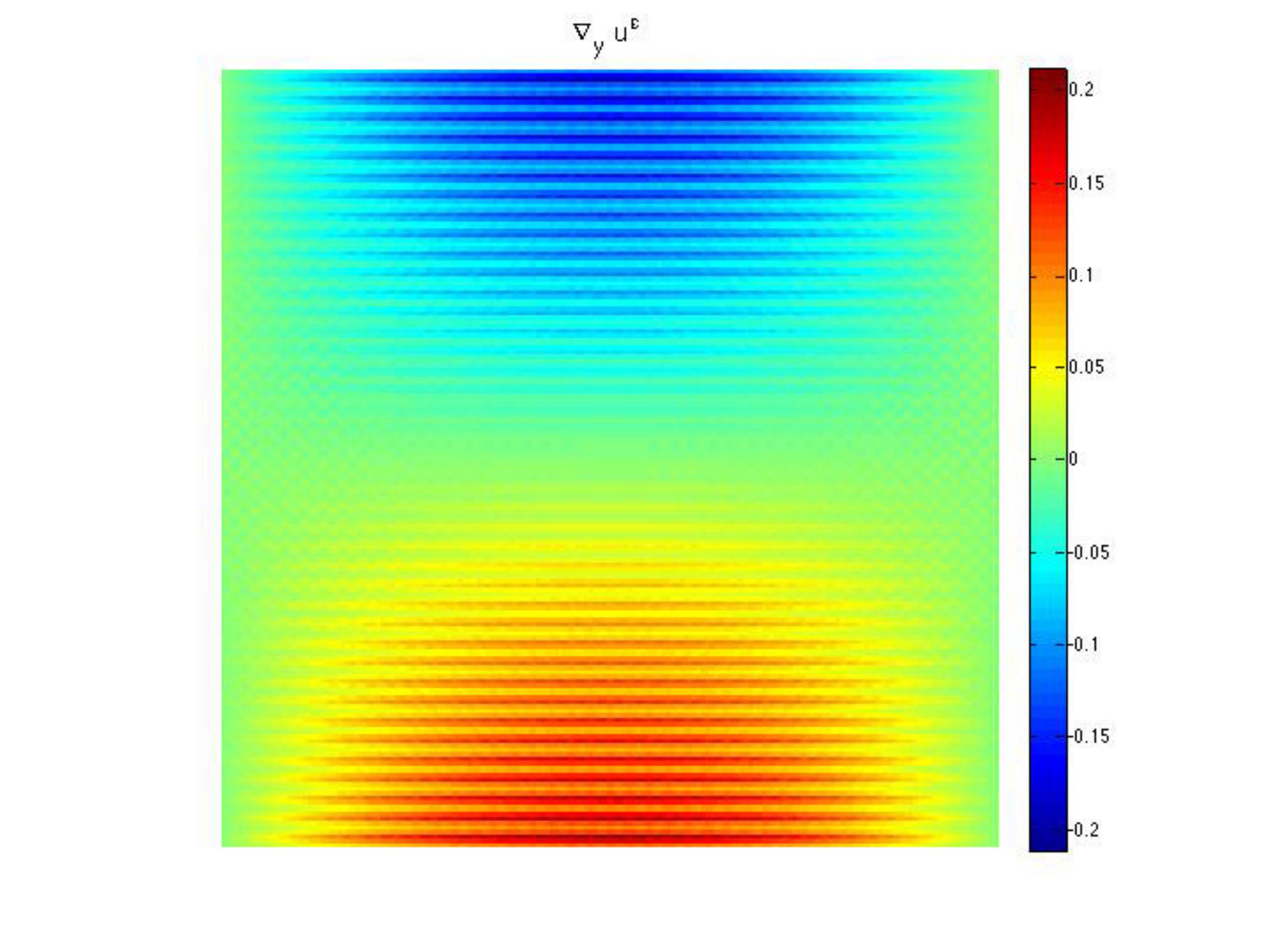}
\includegraphics[width=75mm]{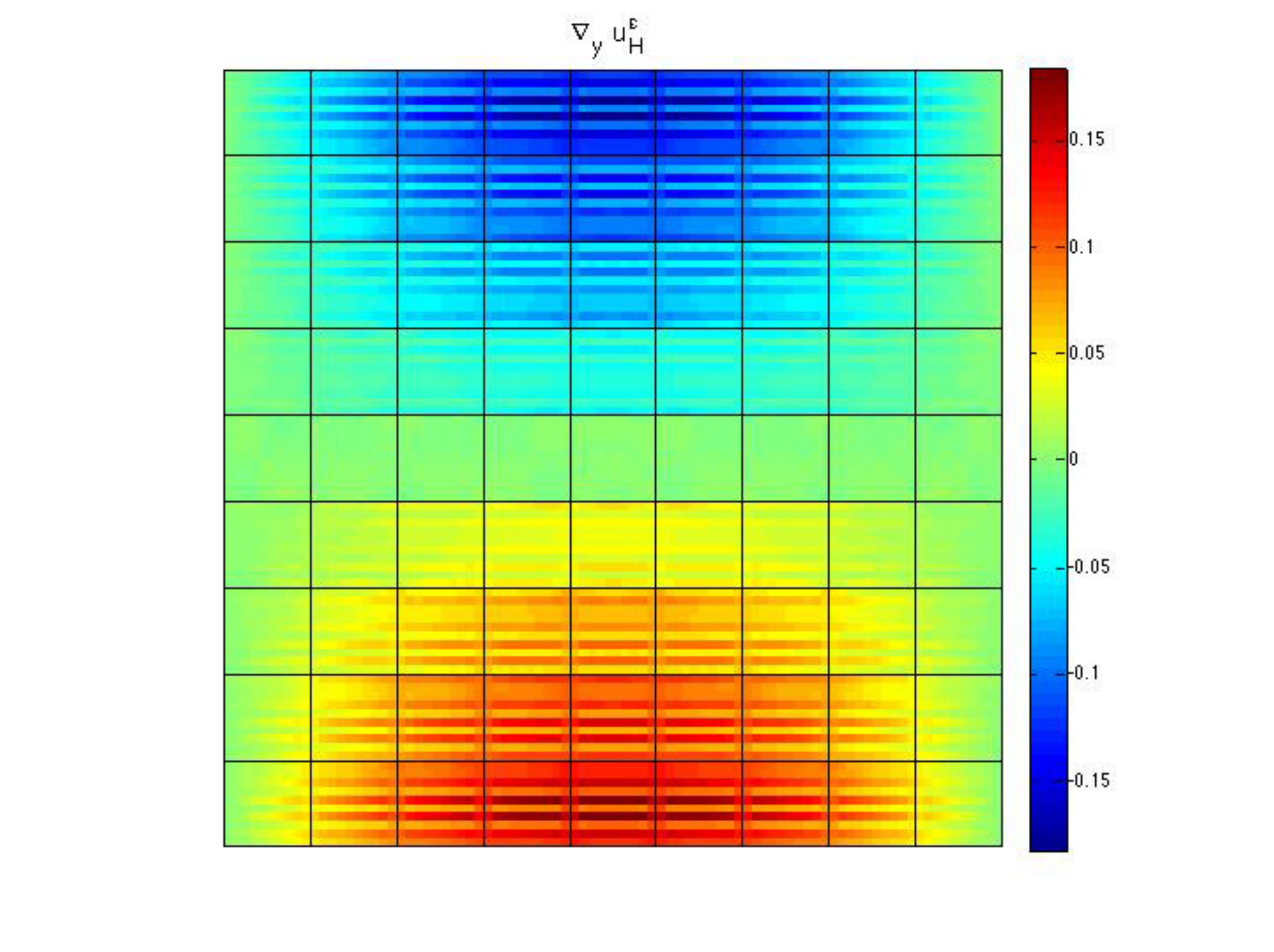}
\end{center}
\caption{Defect problem, exponential loading: exact gradient $\nab u^\eps$ (left) and MsFEM gradient $\nab u^\eps_H$ (right) (top row: components $\nab u^\eps \cdot \be_1$ and $\nab u_H^\eps \cdot \be_1$; bottom row: components $\nab u^\eps \cdot \be_2$ and $\nab u_H^\eps \cdot \be_2$).}
\label{fig:solMsFEM2Ddefect2grad}
\end{figure}

\begin{figure}[H]
\begin{center}
\includegraphics[width=75mm]{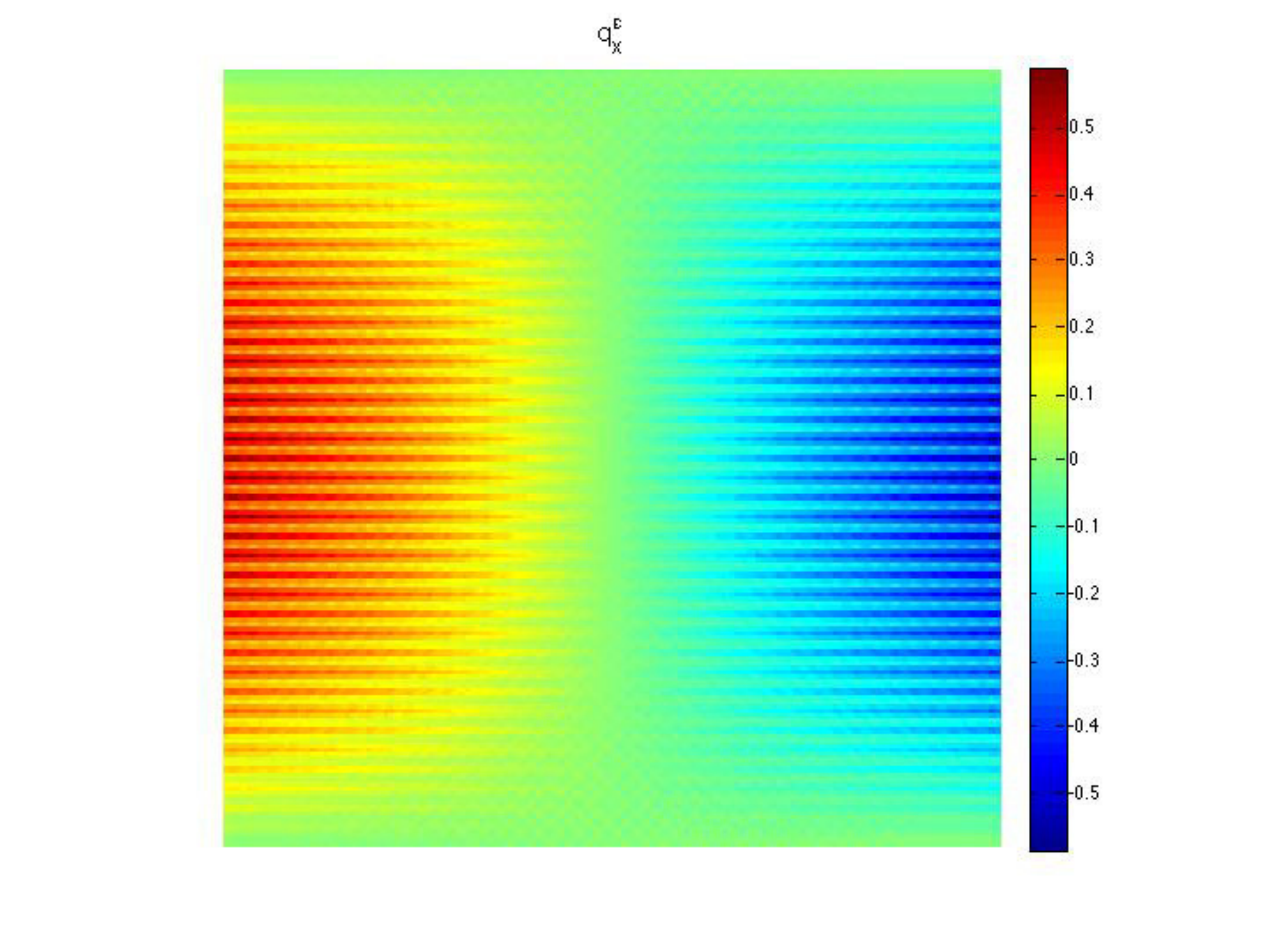}
\includegraphics[width=75mm]{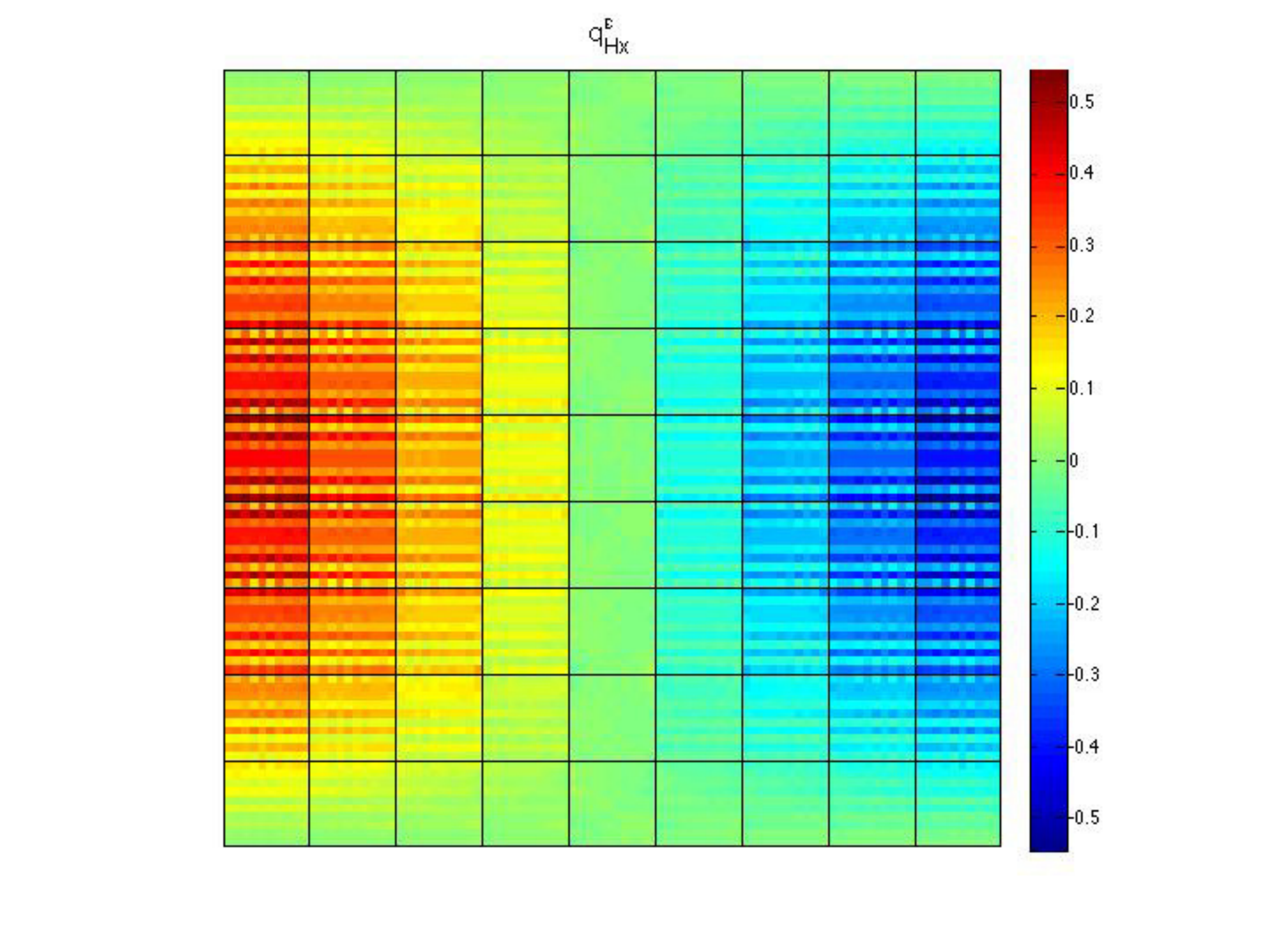}\\
\includegraphics[width=75mm]{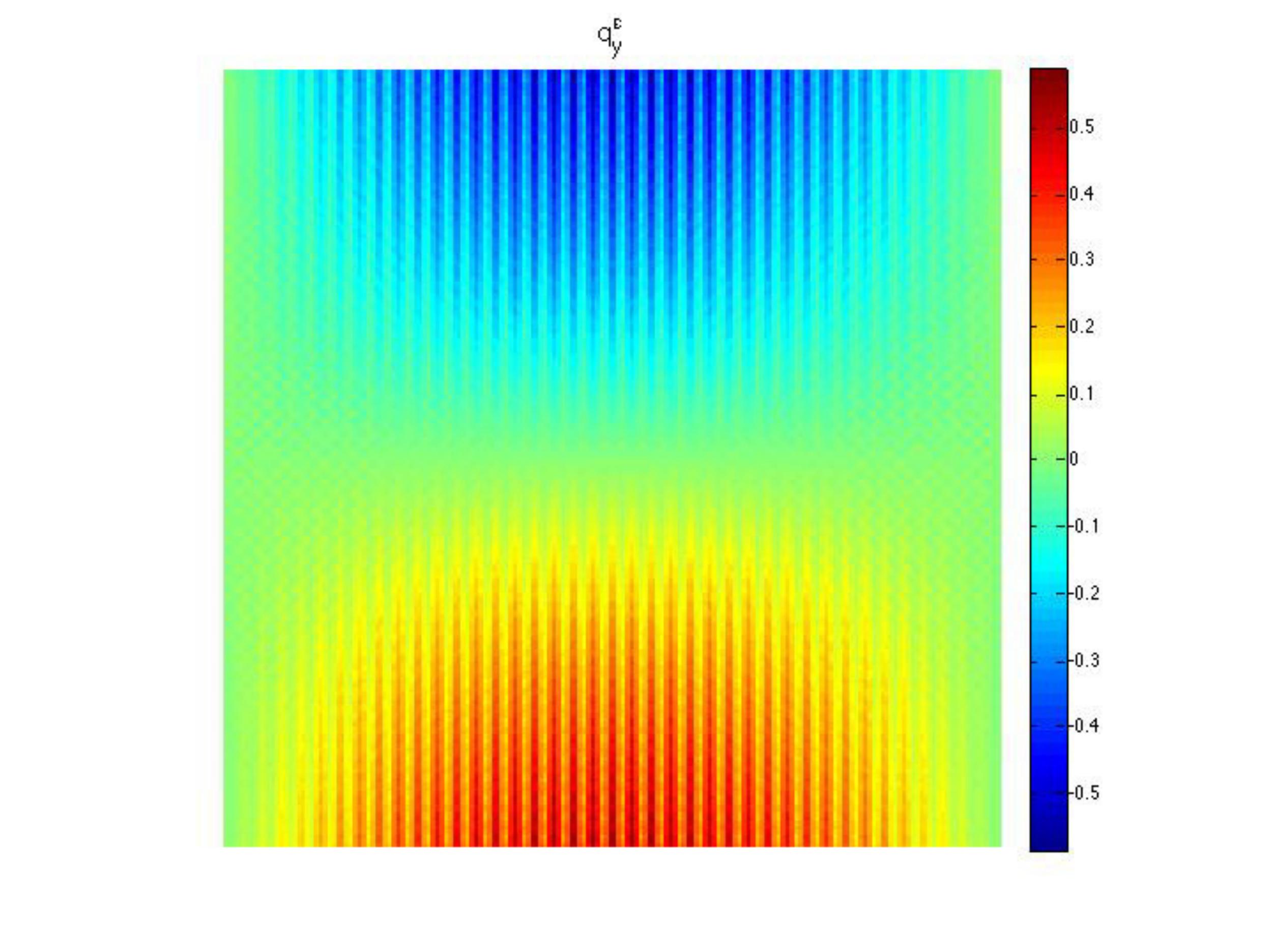}
\includegraphics[width=75mm]{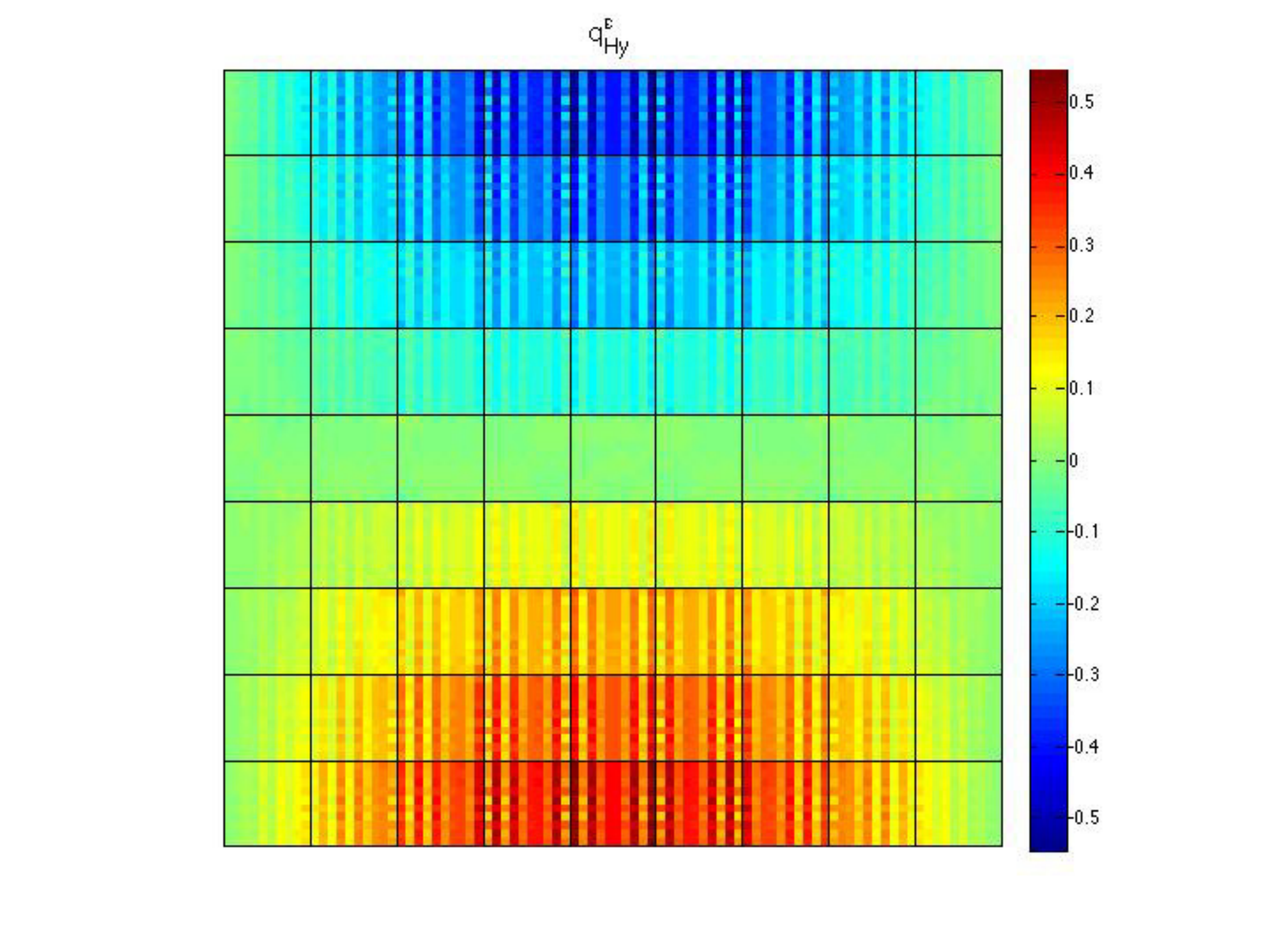}
\end{center}
\caption{Defect problem, exponential loading: exact flux $\bq^\eps$ (left) and MsFEM flux $\bq^\eps_H$ (right) (top row: components $\bq^\eps \cdot \be_1$ and $\bq_H^\eps \cdot \be_1$; bottom row: components $\bq^\eps \cdot \be_2$ and $\bq_H^\eps \cdot \be_2$).}
\label{fig:solMsFEM2Ddefect2flux}
\end{figure}

\begin{figure}[H]
\begin{center}
\includegraphics[width=75mm]{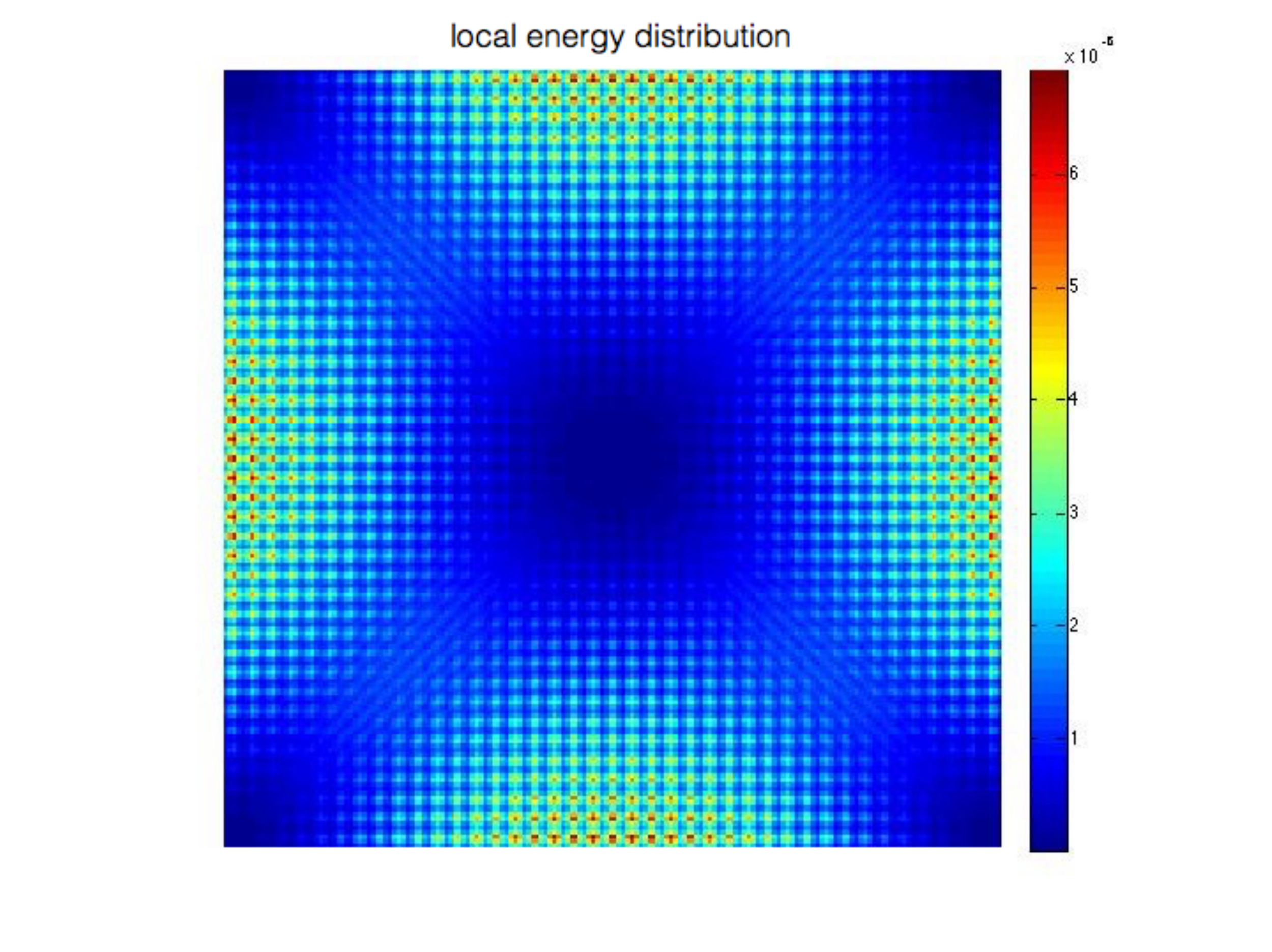}
\includegraphics[width=75mm]{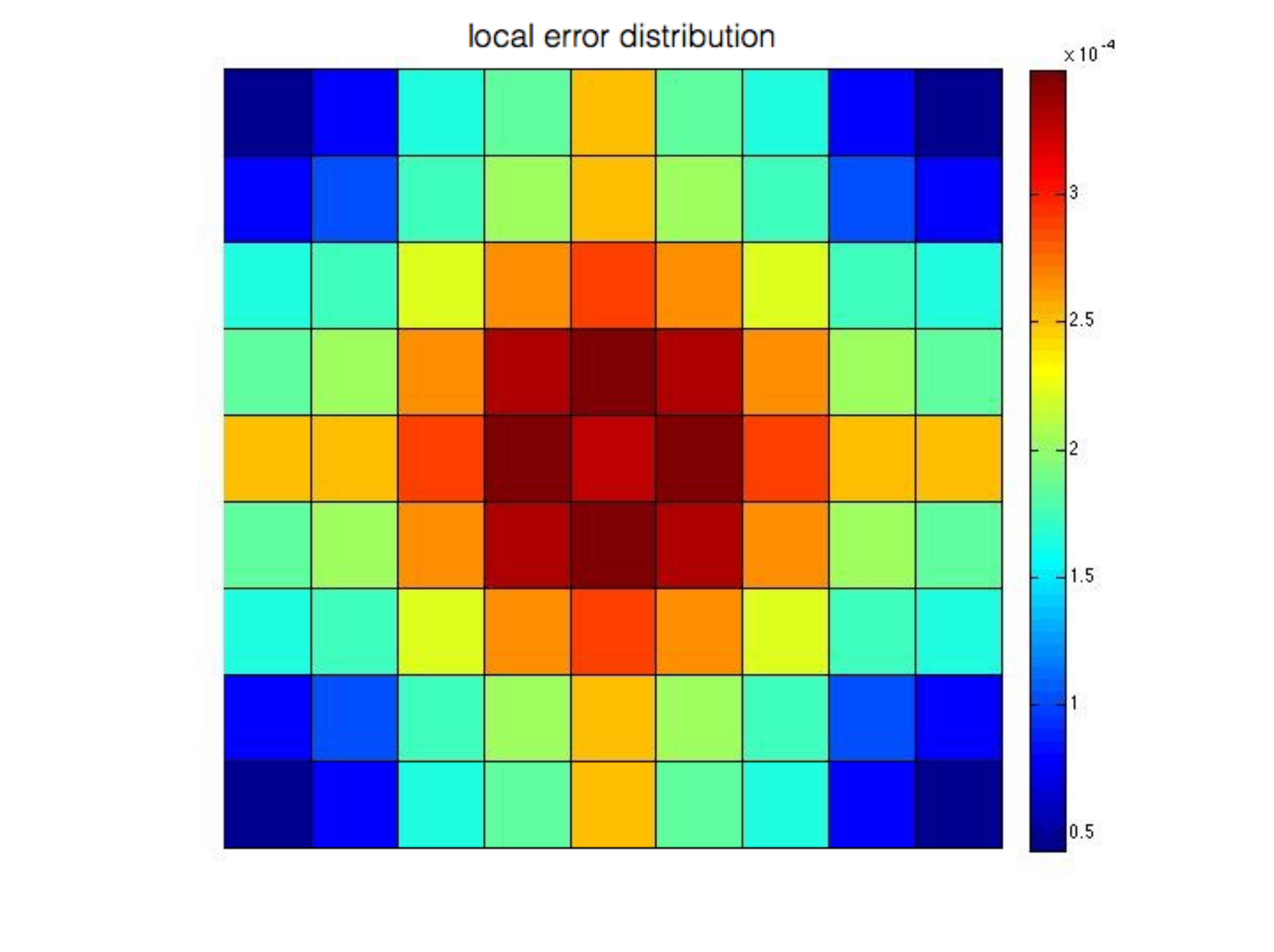}
\end{center}
\caption{Defect problem, exponential loading: distribution of the energy of the exact solution (left) and of the global error (right).}
\label{fig:ener2Ddefect2}
\end{figure}

With the aim to consider a quantity of interest that depends on $u^\eps$ (and not on the flux $\bq^\eps$ as in the previous example), we consider here a different loading: $f(x_1,x_2)= \exp(-(x_1^2+x_2^2))$. We again take an initial mesh $\mT_H$ made of $9\times 9$ macro elements, with $h_K=\eps/3$ for any $K$ and without any oversampling. The reference solution is computed on a $500 \times 500$ fine mesh. The MsFEM solution $u^\eps_H$ is compared to the reference solution $u^\eps$ on Fig.~\ref{fig:solMsFEM2Ddefect2u}, their gradients $\nab u^\eps_H$ and $\nab u^\eps$ on Fig.~\ref{fig:solMsFEM2Ddefect2grad}, and fluxes $\bq^\eps_H$ and $\bq^\eps$ on Fig.~\ref{fig:solMsFEM2Ddefect2flux}. We also represent the distribution of the energy of the exact solution and of the global error on Fig.~\ref{fig:ener2Ddefect2}. In contrast to the previous test-case, $u^\eps$ is not close to 0 in the zone where the defect is located (see Fig.~\ref{fig:solMsFEM2Ddefect2u}). This is thus a relevant test-case for considering as quantity of interest the average of $u^\eps$ in that zone.

\begin{figure}[H]
\begin{center}
\includegraphics[width=50mm]{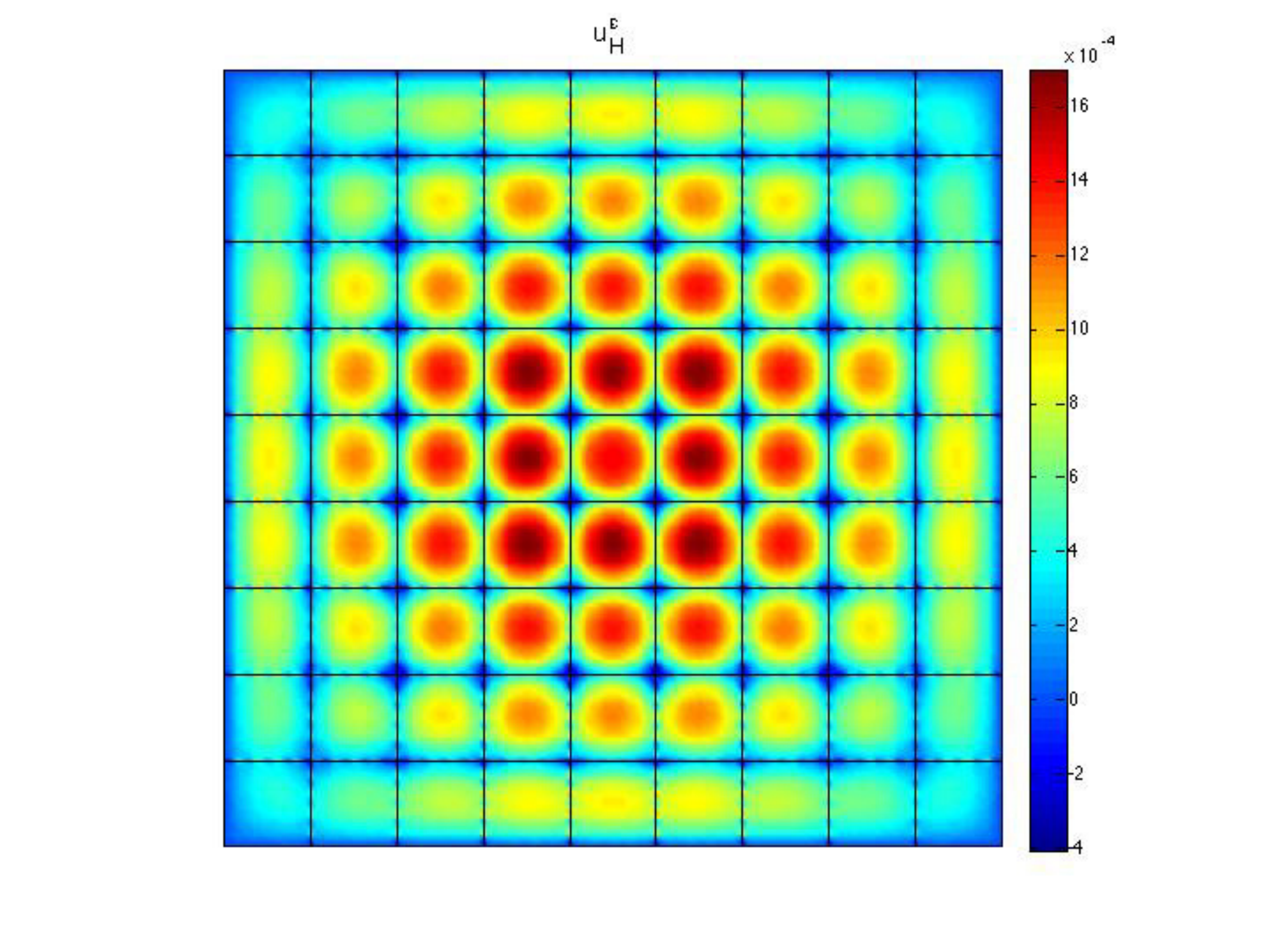}
\includegraphics[width=50mm]{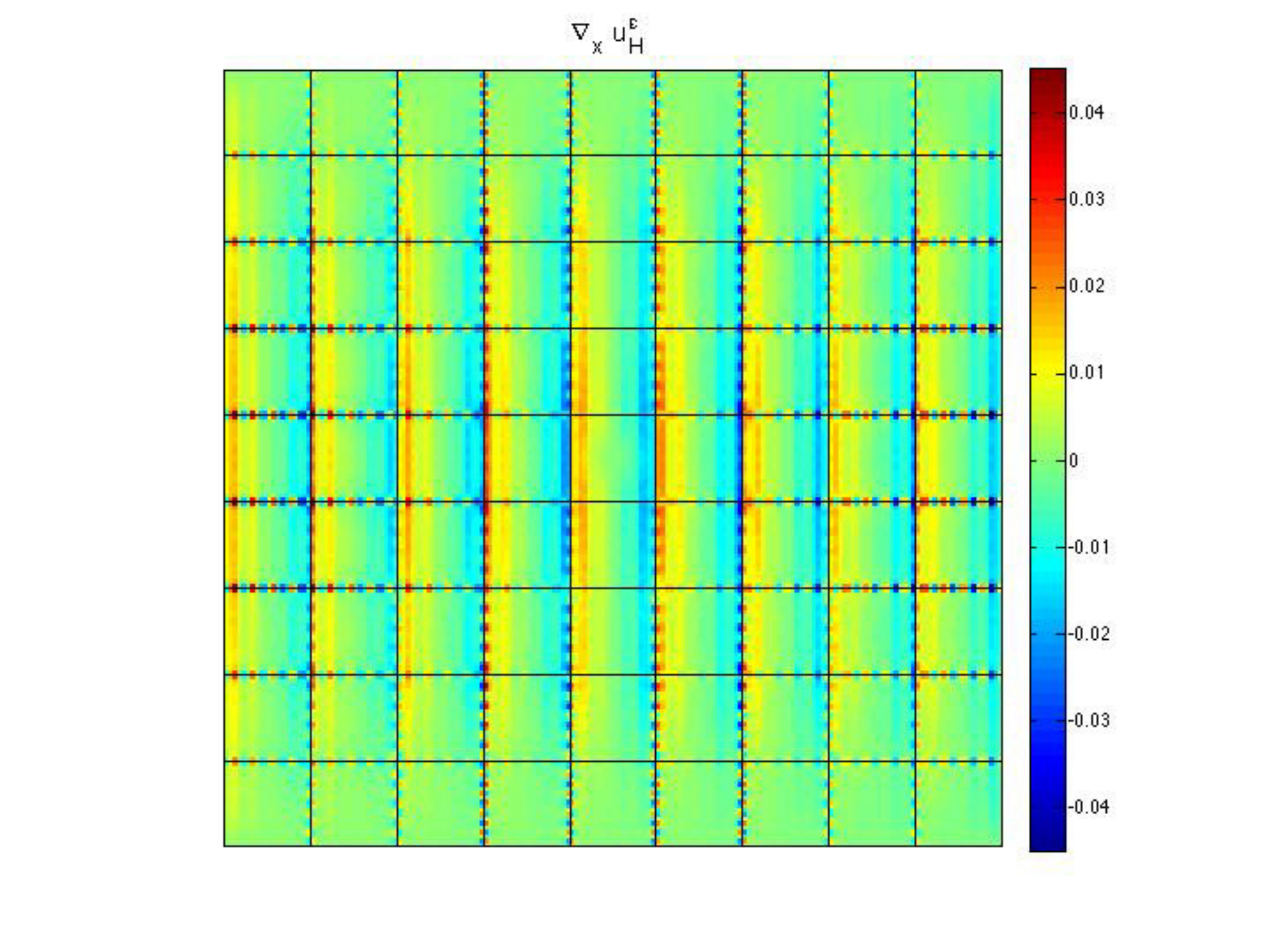}
\includegraphics[width=50mm]{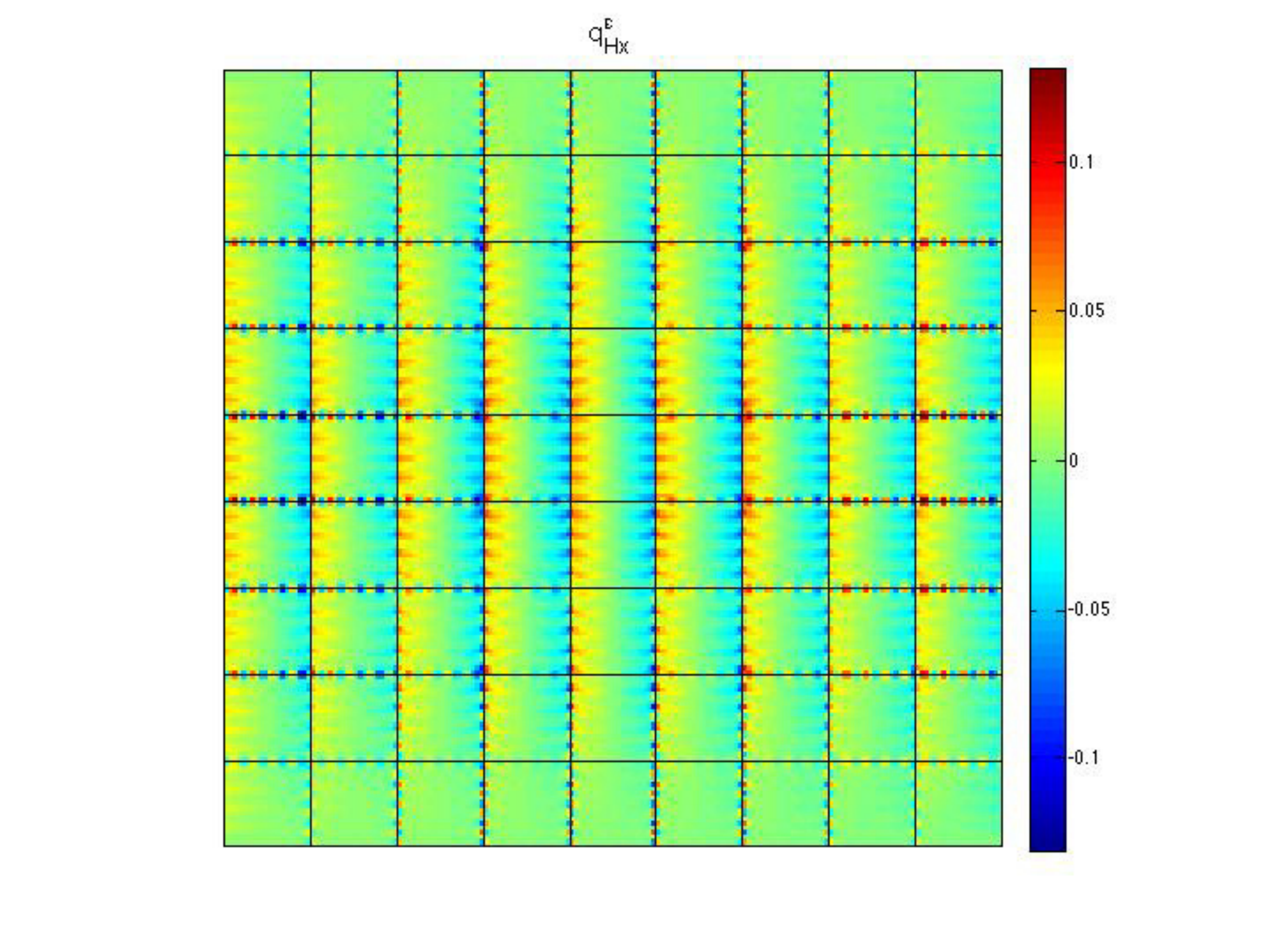}
\end{center}
\caption{Defect problem, exponential loading: error between exact and MsFEM fields for $u^\eps$ (left), $\nab u^\eps \cdot \be_1$ (center) and $\bq^\eps \cdot \be_1$ (right).}
\label{fig:errorloc2Ddefect2flux}
\end{figure}

We also plot on Fig.~\ref{fig:errorloc2Ddefect2flux} the distribution of the error on $u^\eps$, $\nab u^\eps \cdot \be_1$ and $\bq^\eps \cdot \be_1$.

\begin{figure}[H]
\begin{center}
\includegraphics[width=75mm]{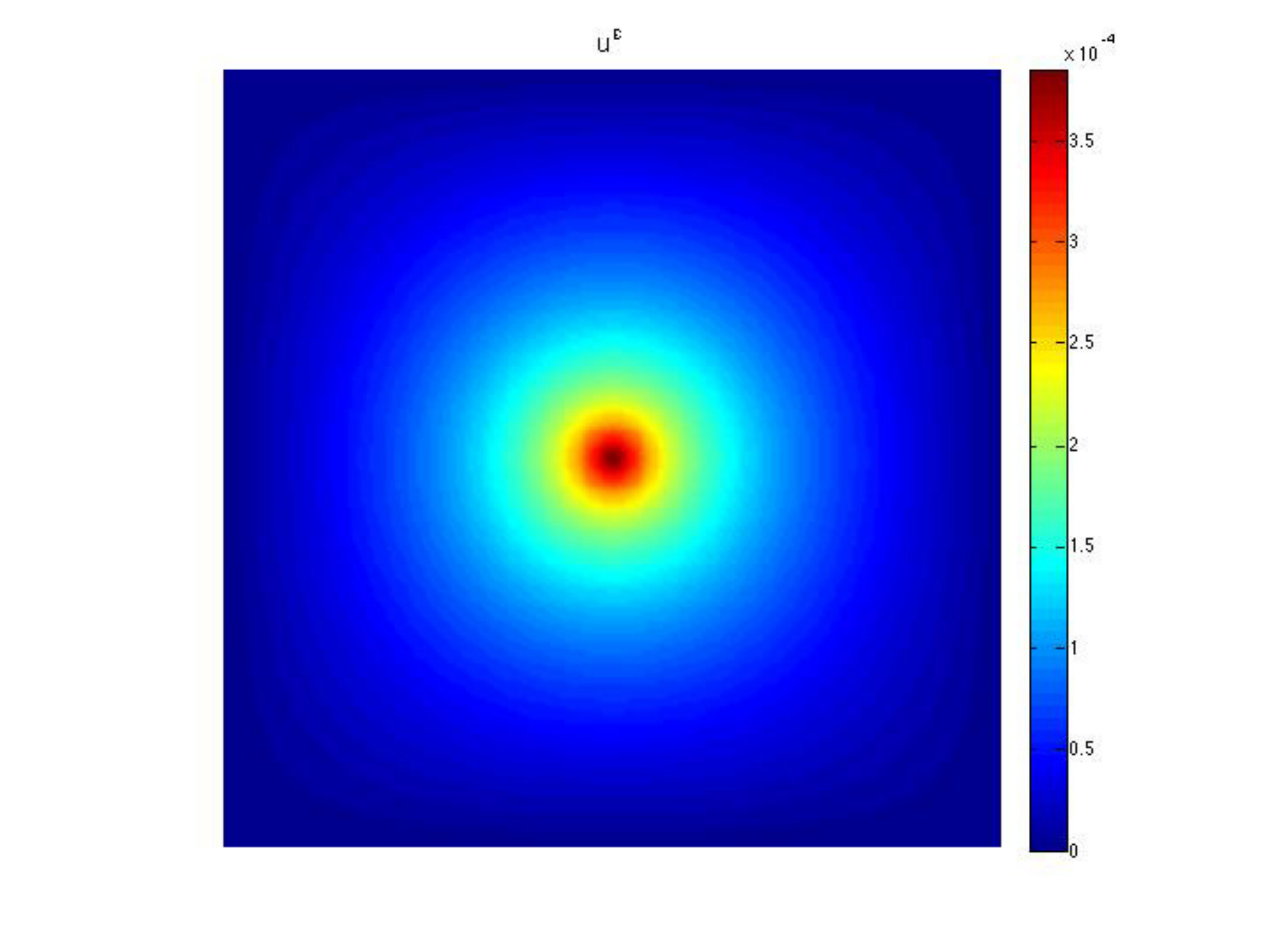}
\includegraphics[width=75mm]{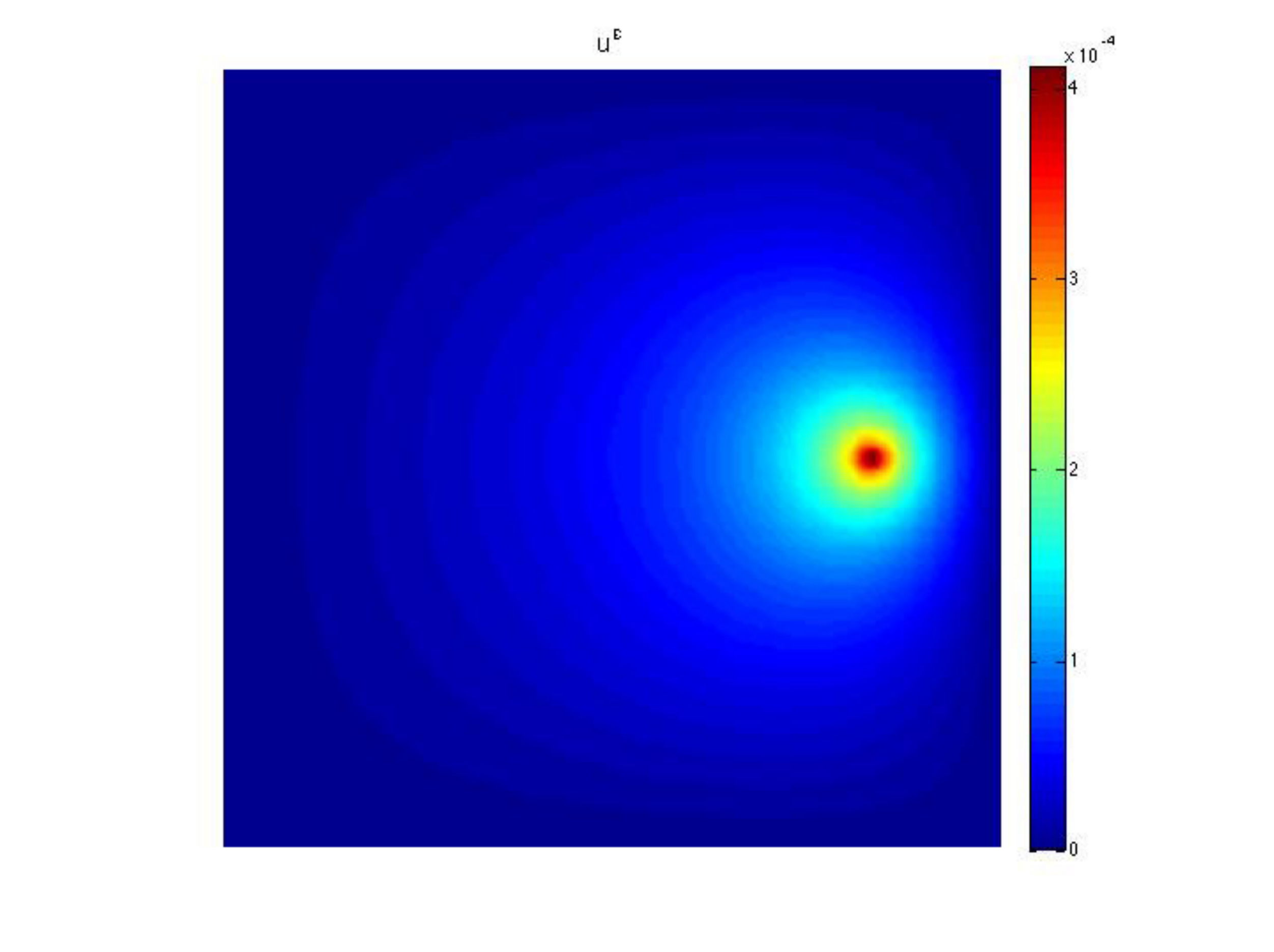}
\end{center}
\caption{Defect problem, exponential loading: exact adjoint solution $\widetilde{u}^\eps$ for $Q_1$ (left) and $Q_2$ (right).}
\label{fig:soladj2Ddefect2u}
\end{figure}

\begin{figure}[H]
\begin{center}
\includegraphics[width=75mm]{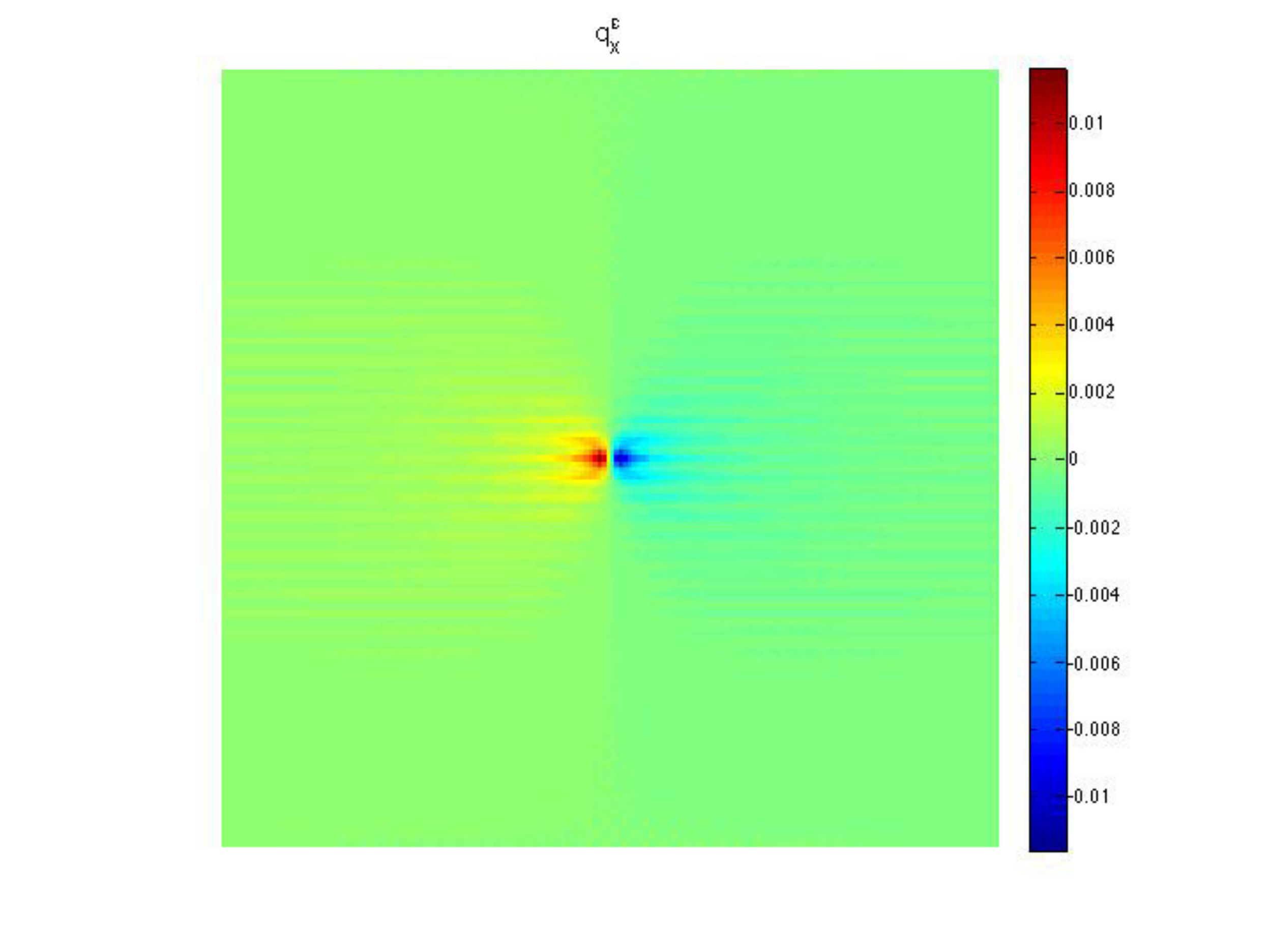}
\includegraphics[width=75mm]{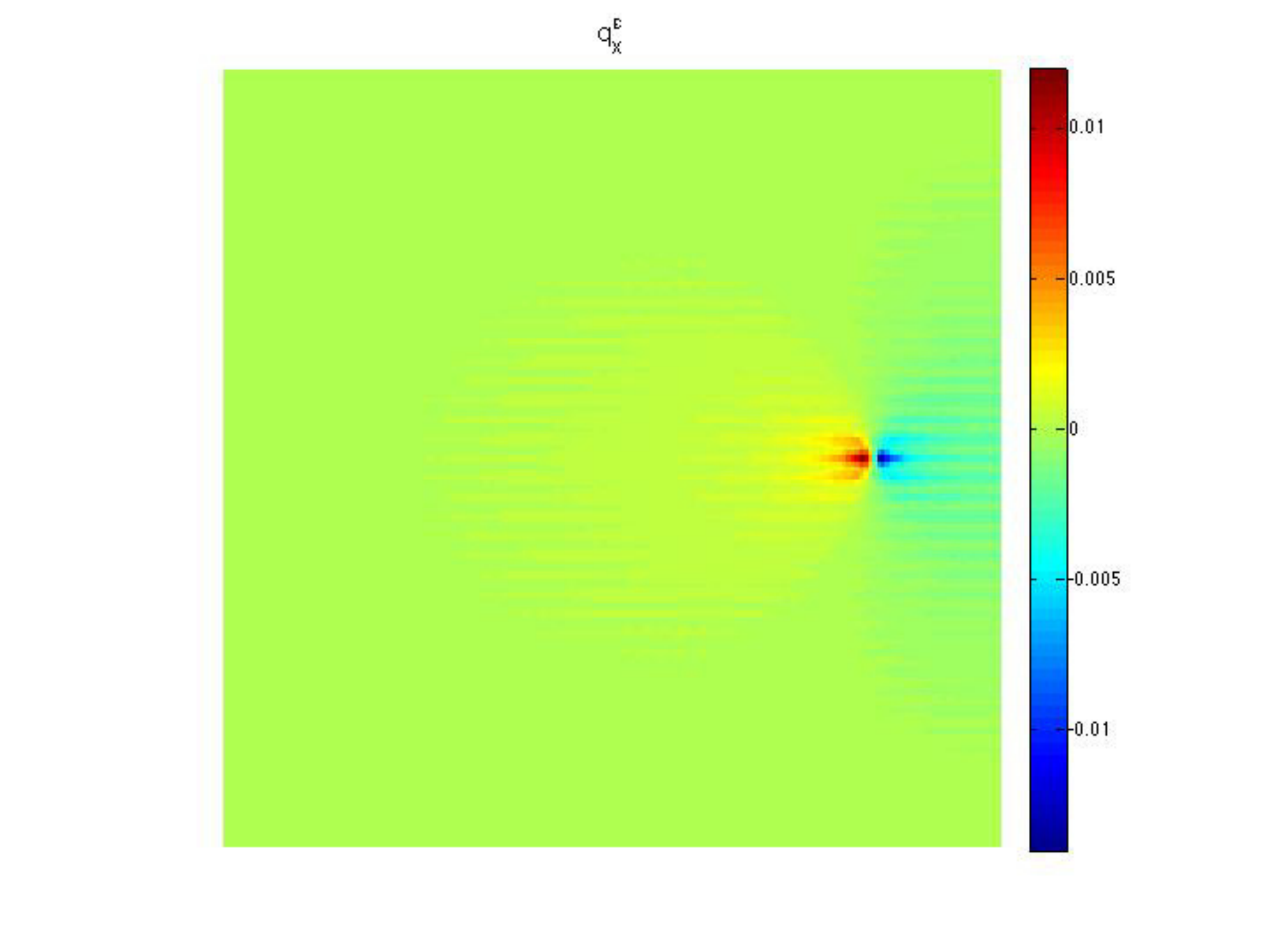}\\
\includegraphics[width=75mm]{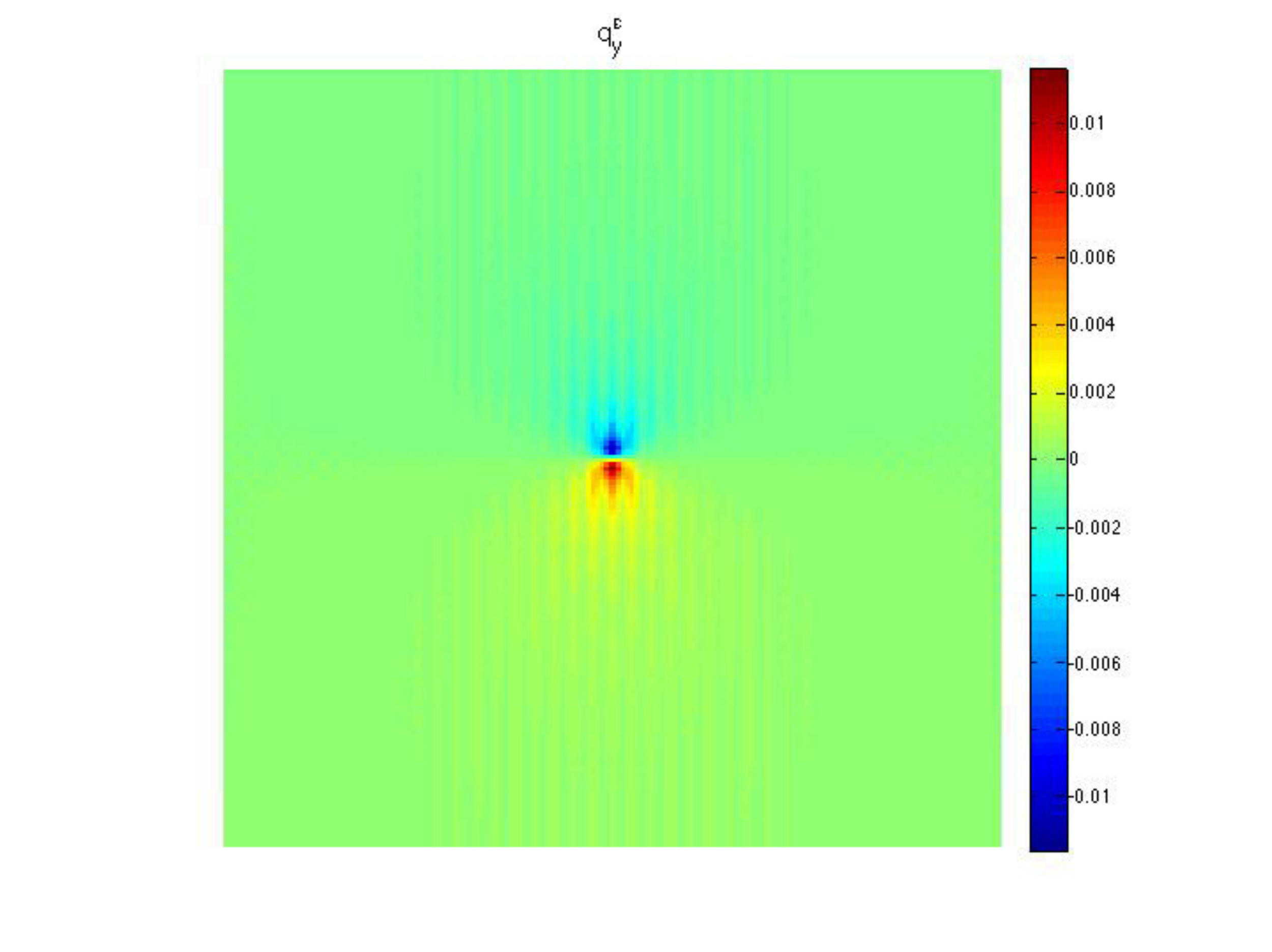}
\includegraphics[width=75mm]{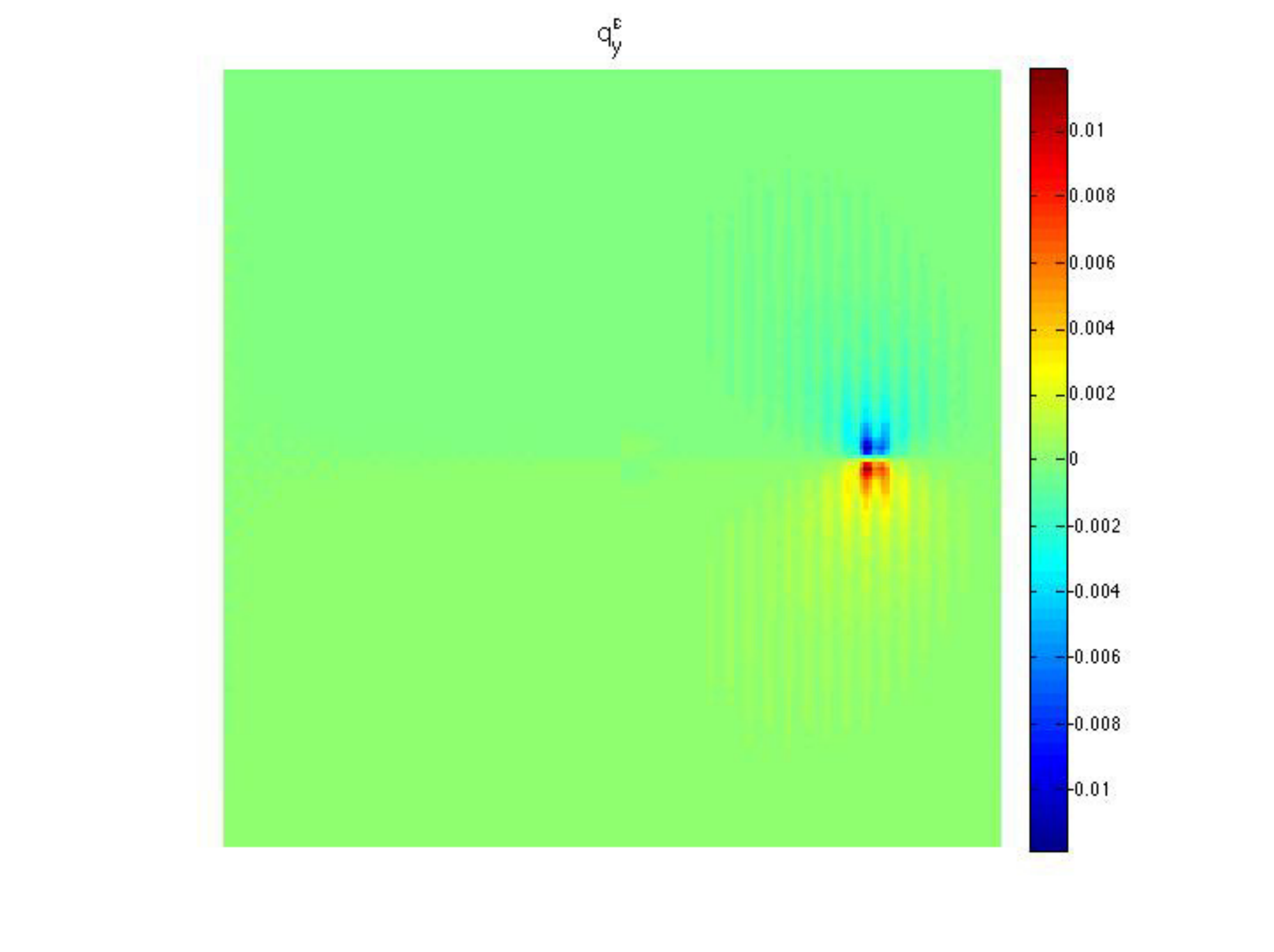}
\end{center}
\caption{Defect problem, exponential loading: exact adjoint flux $\widetilde{\bq}^\eps$ for $Q_1$ (left) and $Q_2$ (right) (top row: component along $\be_1$; bottom row: component along $\be_2$).}
\label{fig:soladj2Ddefect2q}
\end{figure}

As pointed out above, we choose the quantity of interest $\dis Q(u^\eps)=\frac{1}{|\omega|}\int_{\omega} u^\eps$ with $|\omega| = 4 \eps \times 4 \eps$. We consider two cases:
\begin{itemize}
\item Case 1: the subregion $\omega$ is centered on the defect ($Q=Q_1$);
\item Case 2: the subregion $\omega$ is far from the defect ($Q=Q_2$).
\end{itemize}
The corresponding adjoint solutions $\widetilde{u}^\eps$ and fluxes $\widetilde{\bq}^\eps$ are shown on Figs.~\ref{fig:soladj2Ddefect2u} and~\ref{fig:soladj2Ddefect2q}. Again, as expected, the adjoint solutions are essentially supported in the domain of interest $\omega$.

\begin{figure}[H]
\begin{center}
\includegraphics[width=150mm]{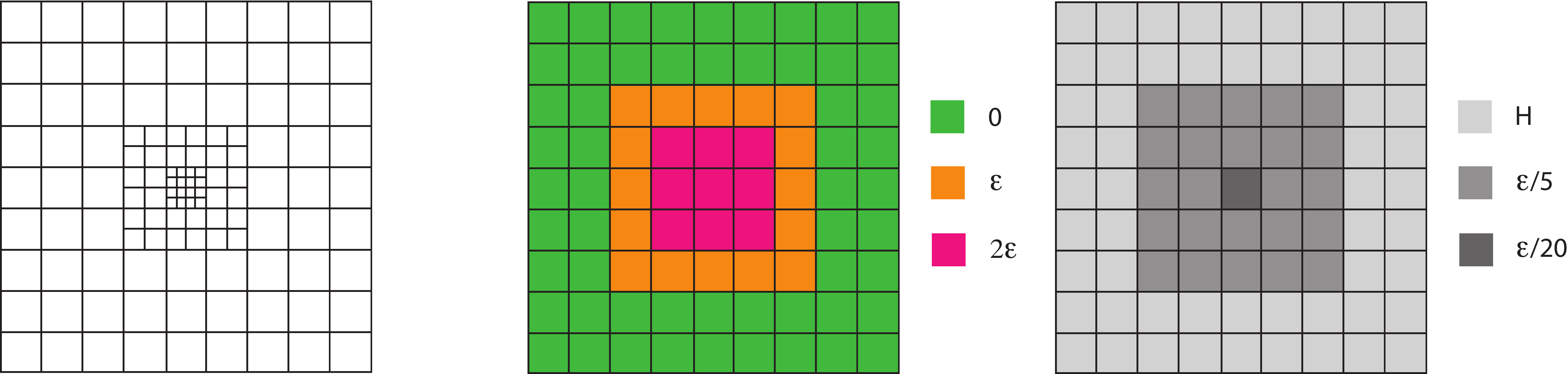} 
\end{center}
\caption{Defect problem, exponential loading: final MsFEM discretization for $Q_1$.}
\label{fig:meshconvQ12Ddefect2}
\end{figure}

\begin{figure}[H]
\begin{center}
\includegraphics[width=150mm]{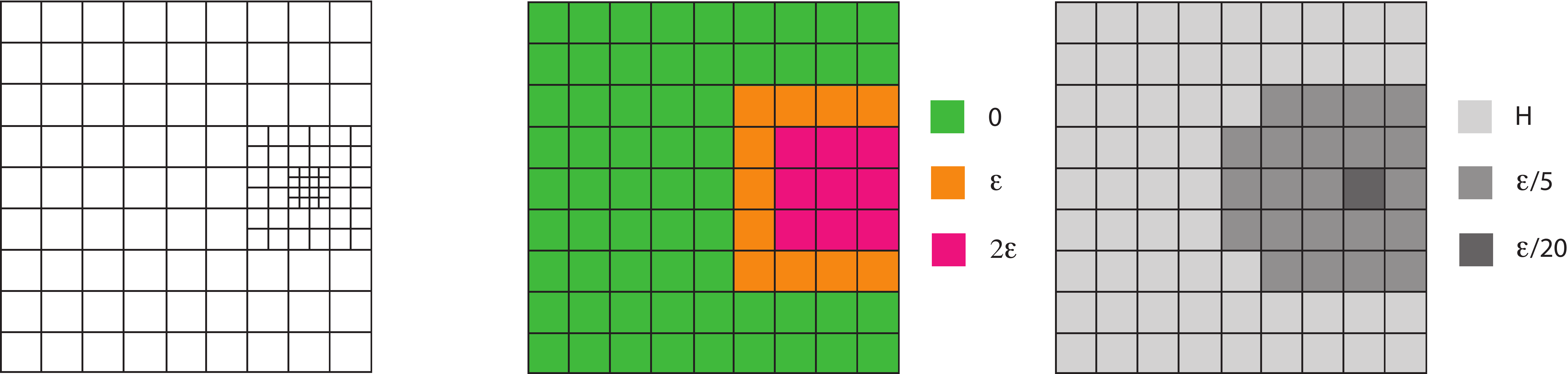}
\end{center}
\caption{Defect problem, exponential loading: final MsFEM discretization for $Q_2$.}
\label{fig:meshconvQ22Ddefect2}
\end{figure}

\begin{figure}[H]
\begin{center}
\includegraphics[width=150mm]{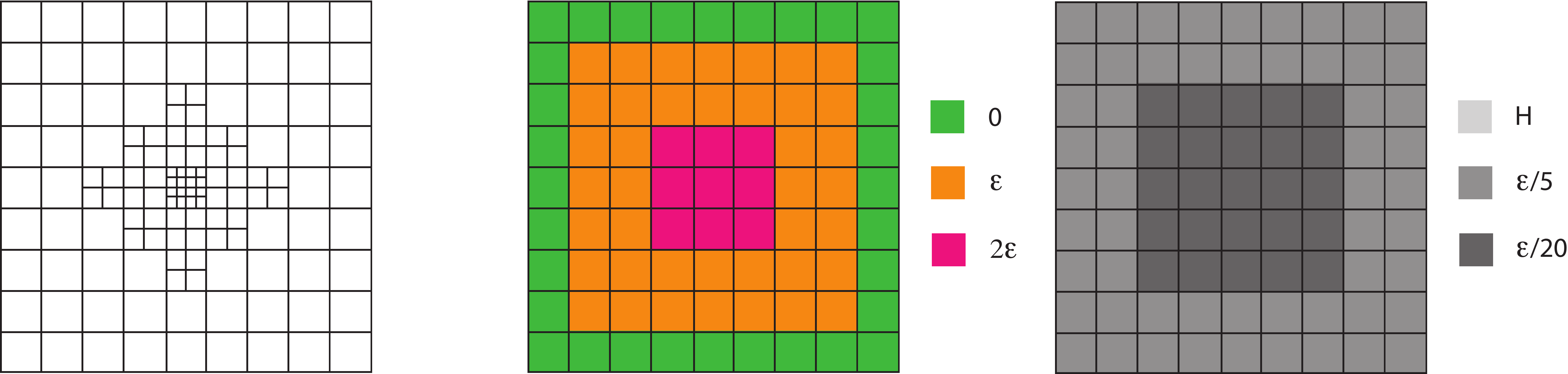} 
\end{center}
\caption{Defect problem, exponential loading: final MsFEM discretization when controlling the global error.}
\label{fig:meshconvglob2Ddefect2}
\end{figure}

The error threshold is fixed to $1\%$. We show on Fig.~\ref{fig:meshconvQ12Ddefect2} (resp. Fig.~\ref{fig:meshconvQ22Ddefect2}) the final MsFEM discretization obtained after using the adaptive algorithm presented in Section~\ref{section:algorithm}, in order to respect the error threshold for the quantity of interest $Q_1$ (resp. $Q_2$). This discretization can be compared with that obtained when controlling the global error (in energy norm) using the tools presented in~\cite{CHA16b}, which is shown on Fig.~\ref{fig:meshconvglob2Ddefect2}. The conclusions are the same as in the case of sinusoidal loading (see Section~\ref{sec:sinus}): (i) it is worth working with a goal-oriented adaptive strategy rather than with an adaptive strategy based on the global error and (ii) the mesh needs to be refined only in the region close to the region of interest $\omega$.

As regards the accuracy of the error estimate, we again observe that it is very good with only a slight overestimation of the error. Computations yield $\eta^Q/|Q(u^\eps)-Q(u^\eps_H)| = 1.14$ (resp. 1.15) at the first iteration of the adaptive algorithm, and $\eta^Q/|Q(u^\eps)-Q(u^\eps_H)| = 1.07$ (resp. 1.04) at the final iteration of the adaptive algorithm, for the quantities of interest $Q_1$ and $Q_2$, respectively.

\begin{figure}[H]
\begin{center}
\includegraphics[width=75mm]{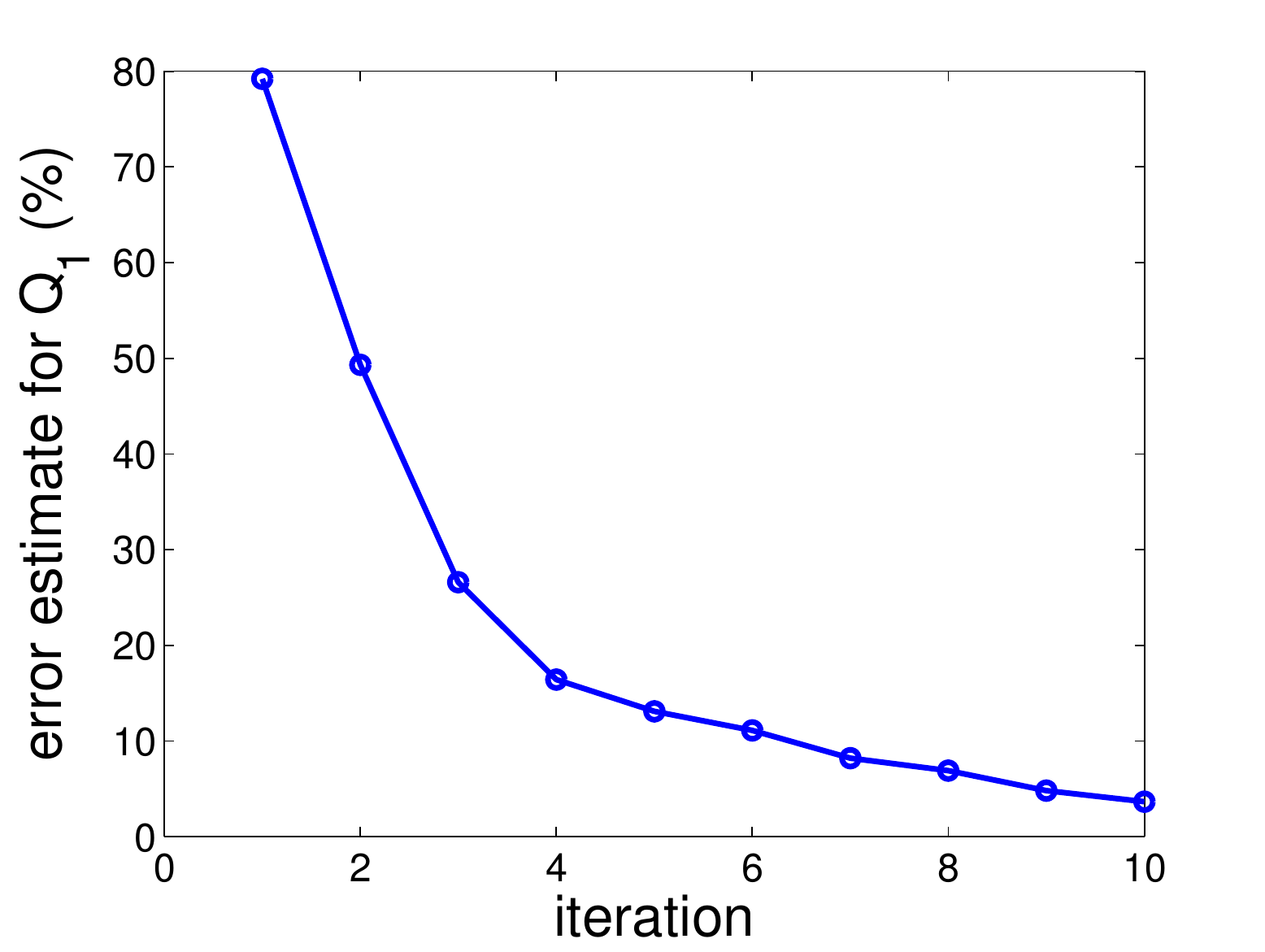} 
\includegraphics[width=75mm]{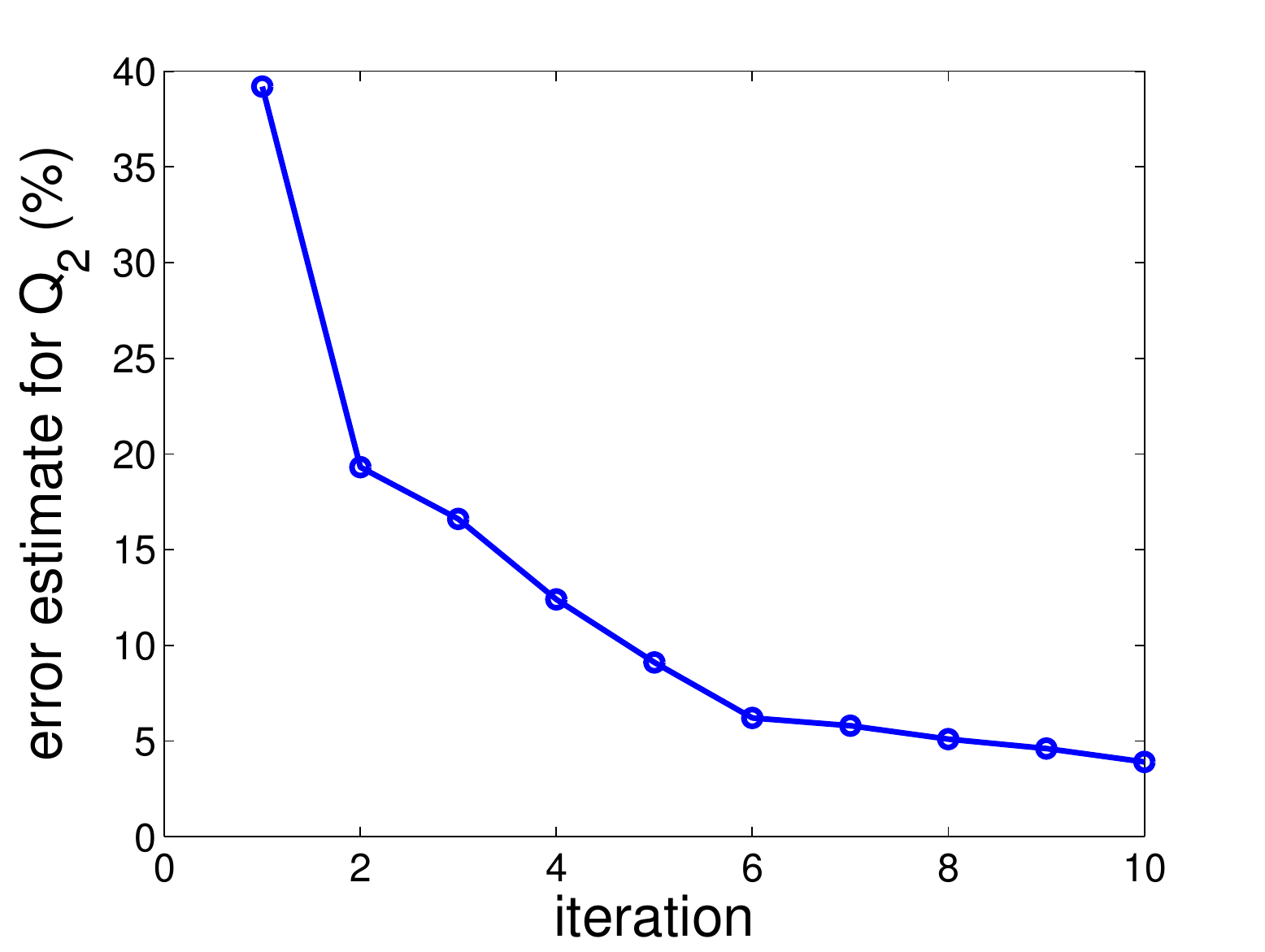} 
\end{center}
\caption{Defect problem, exponential loading: convergence of the error estimate for $Q_1$ (left) and $Q_2$ (right).}
\label{fig:converrorQ2Ddefect2}
\end{figure}

We show on Fig.~\ref{fig:converrorQ2Ddefect2} how the error estimate on $Q$ decreases as the adaptive iterations of the mesh proceed. For instance, after the first adaptation, the error has been divided by a factor 2, going from 80\% to 50\% for $Q_1$ (resp. 40\% to 20\% for $Q_2$).

\subsection{Flow problem}\label{section:flow}

We now consider a flow problem in a fractured porous medium, with a source and a sink. This problem has already been considered in~\cite{CHU16}, with different multiscale techniques. The conductivity within the fractures can be several orders of magnitude higher than the conductivity within the matrix. In flow applications, one needs to obtain a good approximation of the pressure in a neighborhood of the wells.

We consider the domain $\Omega = (0,1)^2$ and the inflow-outflow source term $f = \chi_{K_1} - \chi_{K_2}$ , where $K_1 = [0.1, 0.2]\times[0.8, 0.9]$ and $K_2 = [0.8, 0.9]\times[0.1, 0.2]$. The pressure $u^\eps$ is solution to 
$$
-\nab \cdot (A^\eps \nab u^\eps) = f \ \ \text{in $\Omega$}, \qquad u^\eps = 0 \ \ \text{on $\partial \Omega$},
$$
and we consider the quantity of interest $\dis Q(u^\eps) = \frac{1}{|K_2|}\int_{K_2}u^\eps$, that is the average value of the solution (i.e. the pressure) on the outflow region $K_2$ (where the well is located).

\begin{figure}[H]
\begin{center}
\includegraphics[width=75mm]{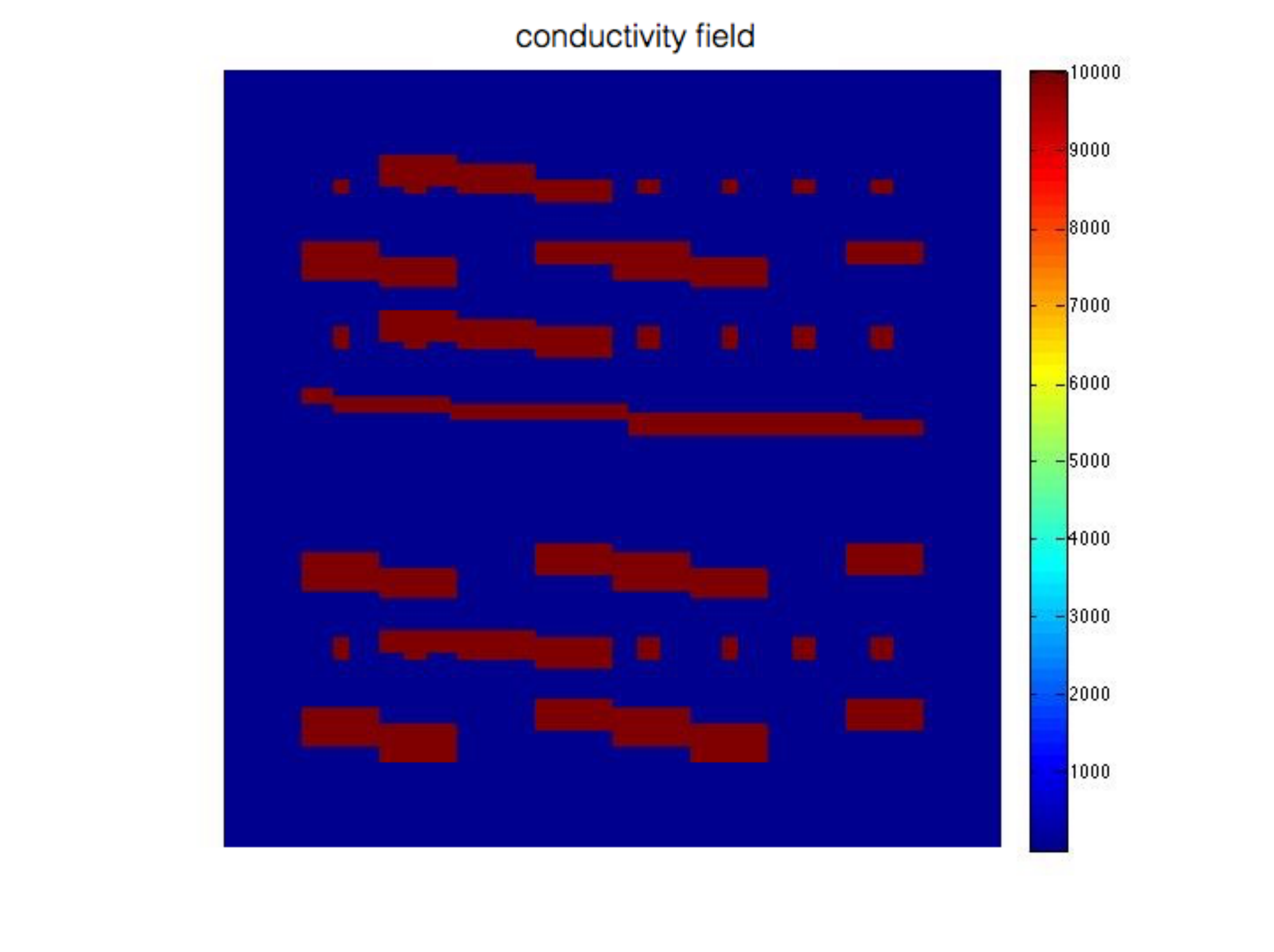}
\end{center}
\caption{Flow problem: evolution of $A^\eps$ in the domain $\Omega$.}
\label{fig:evolAE2Dflow}
\end{figure}

The scalar-valued high-contrast function $A^\eps(\bx)$ is shown on Fig.~\ref{fig:evolAE2Dflow}. We observe a high-conductivity channel crossing the domain $\Omega$, and separating the inflow and outflow regions. In the blue region, we have $A^\eps=1$, and in the red region we set $A^\eps=10^4$. The contrast between the conductivities in the blue and red regions is thus very large.

\begin{figure}[H]
\begin{center}
\includegraphics[width=75mm]{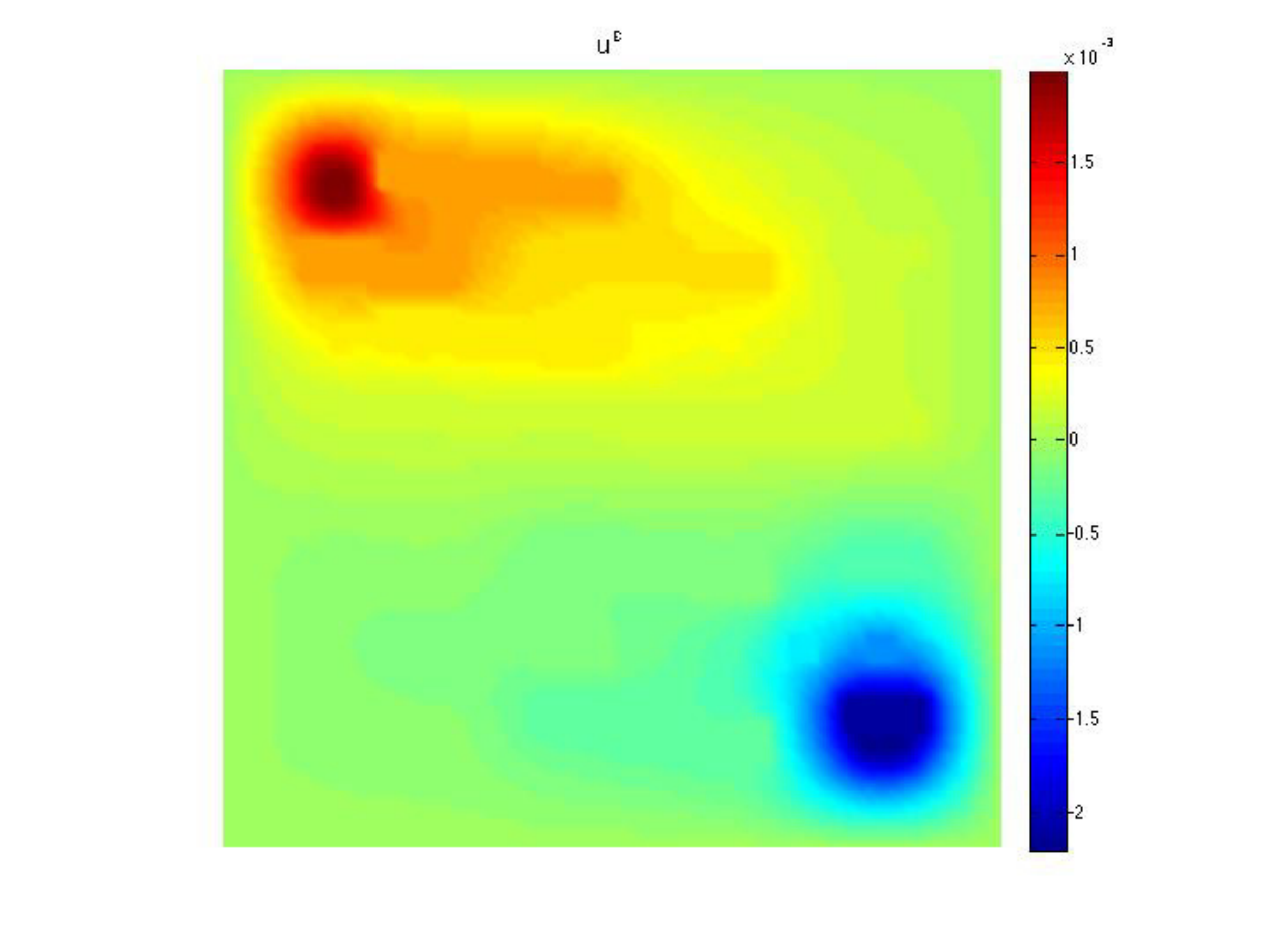}
\includegraphics[width=75mm]{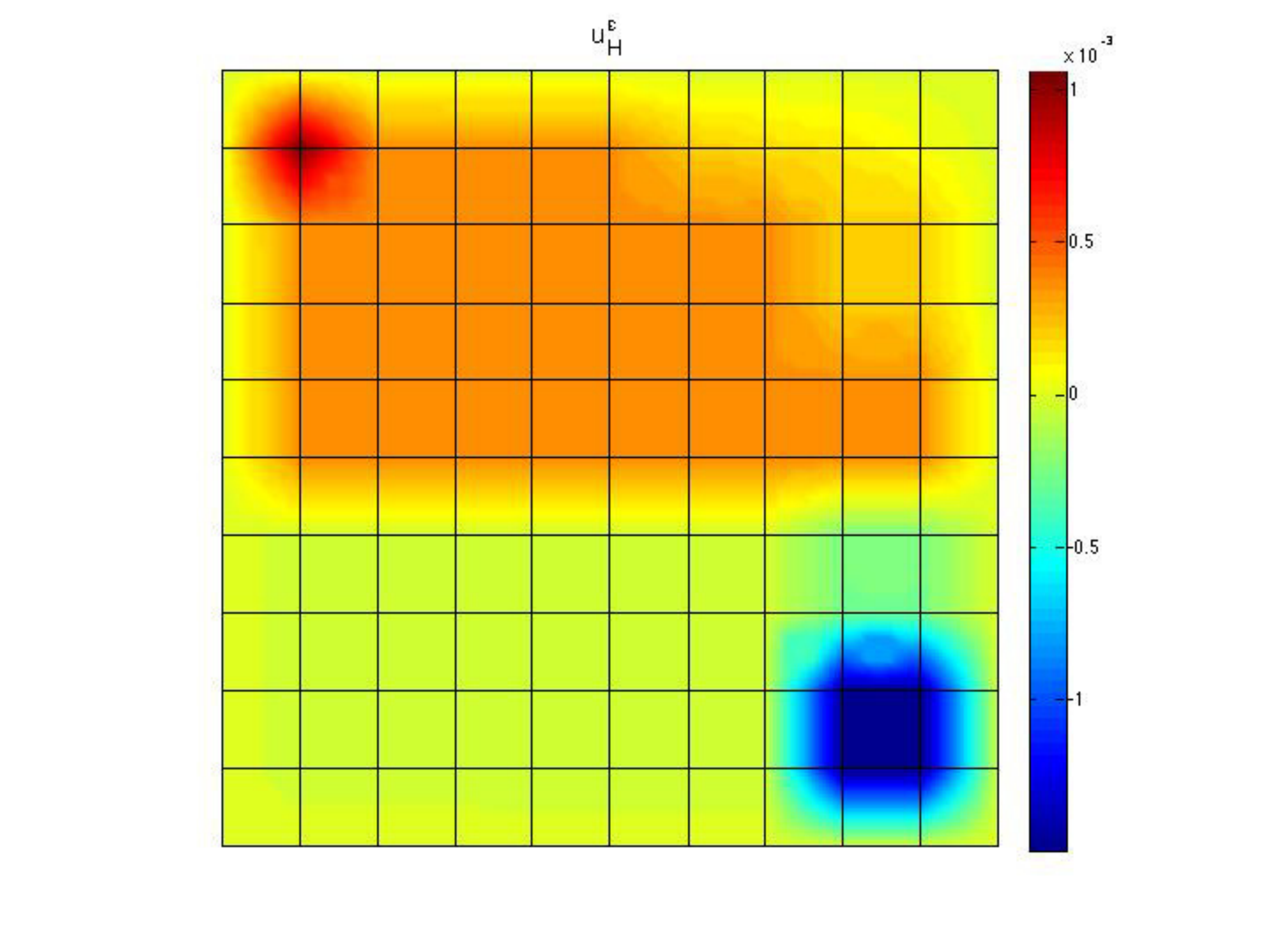}
\end{center}
\caption{Flow problem: exact solution $u^\eps$ (left) and MsFEM solution $u^\eps_H$ (right).}
\label{fig:solMsFEM2Dflowu}
\end{figure}

\begin{figure}[H]
\begin{center}
\includegraphics[width=75mm]{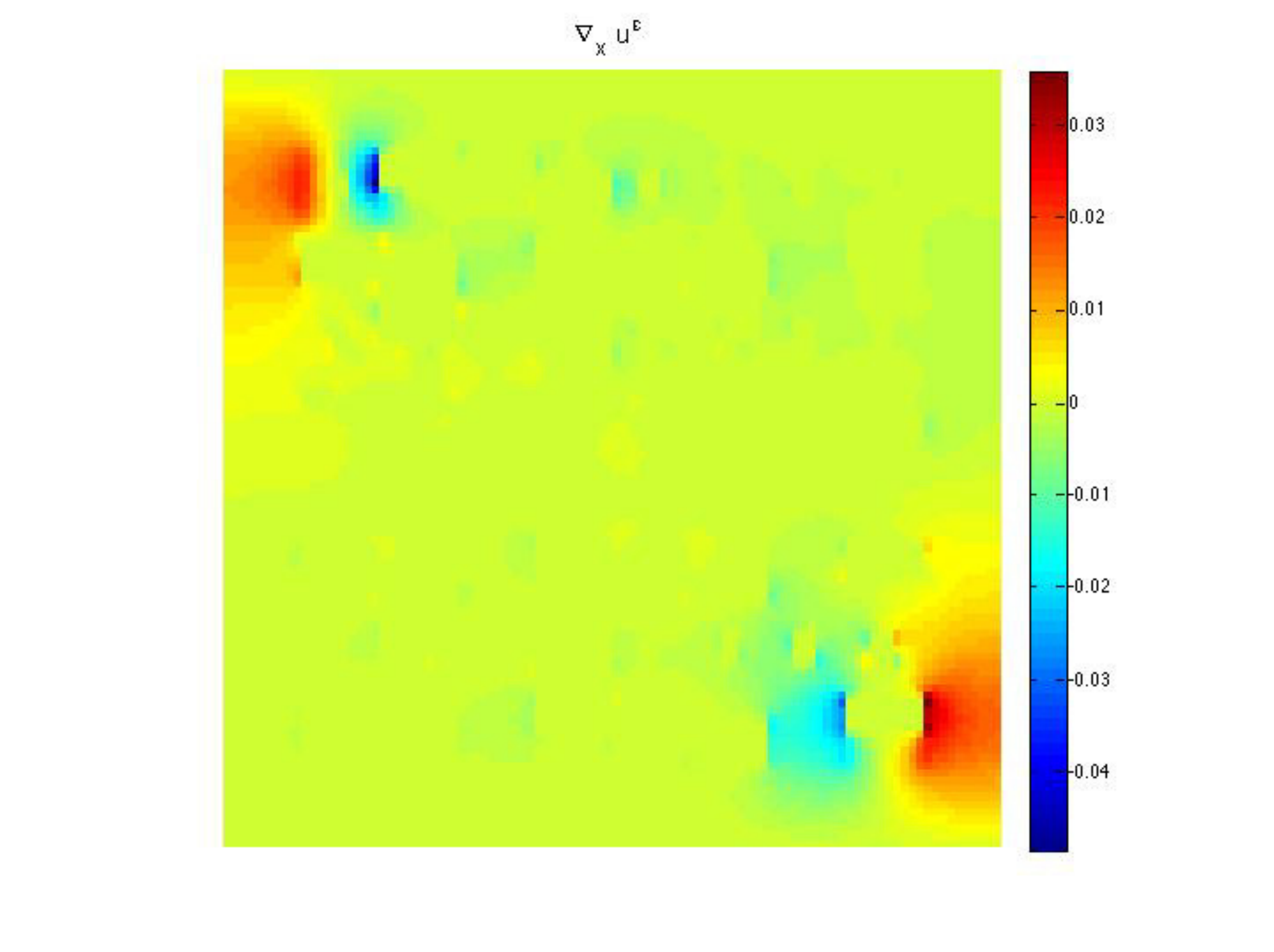}
\includegraphics[width=75mm]{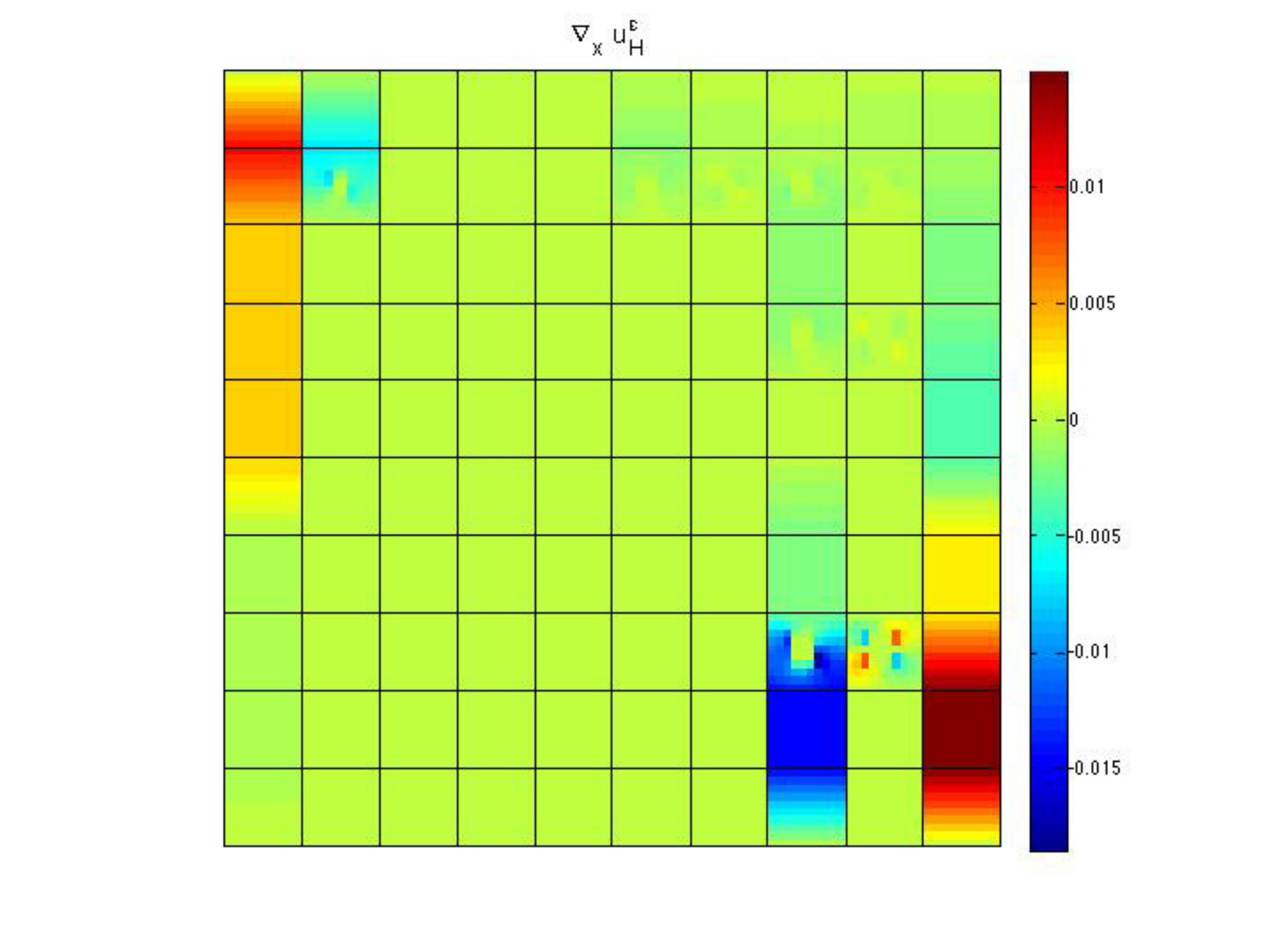}\\
\includegraphics[width=75mm]{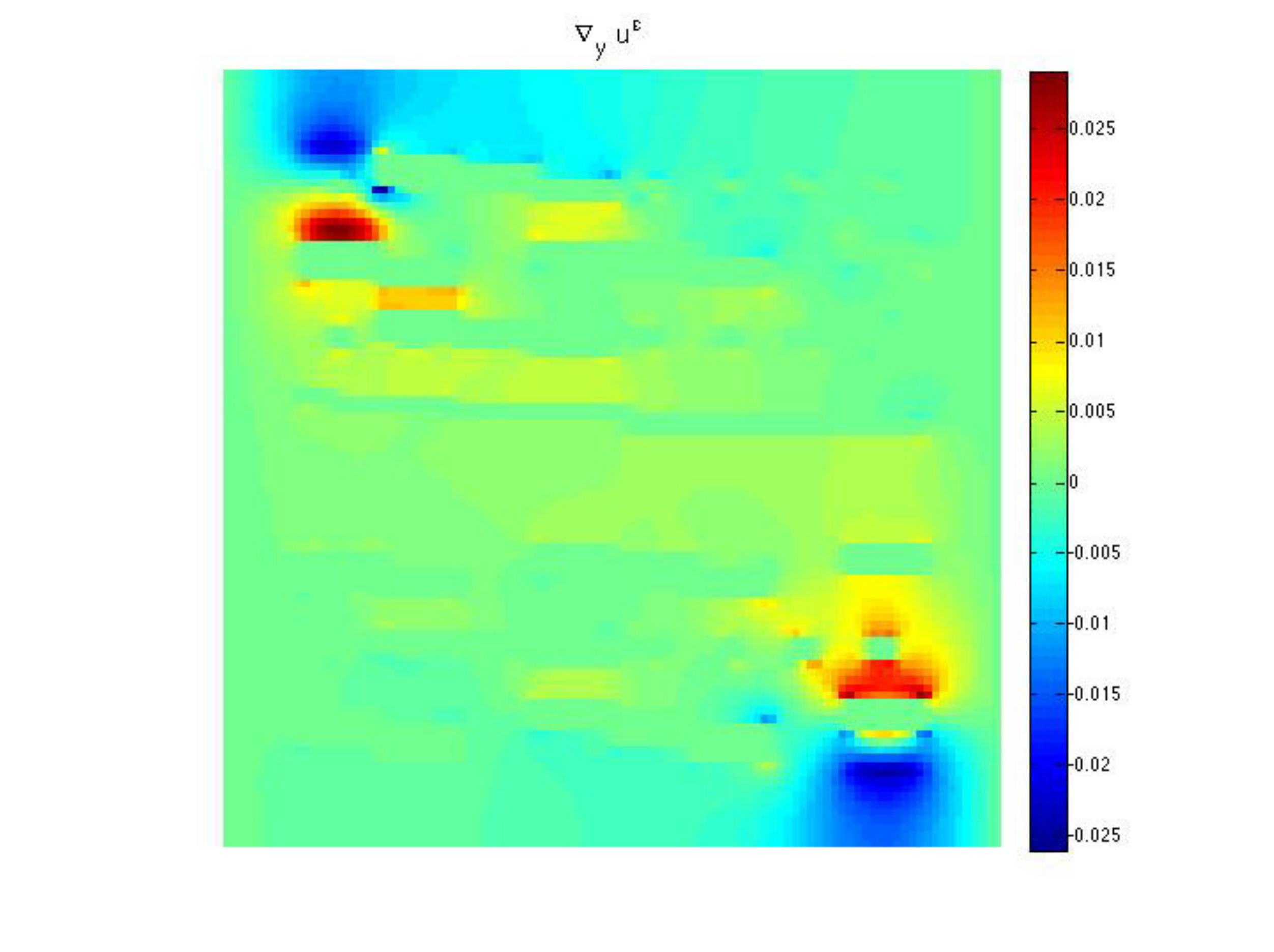}
\includegraphics[width=75mm]{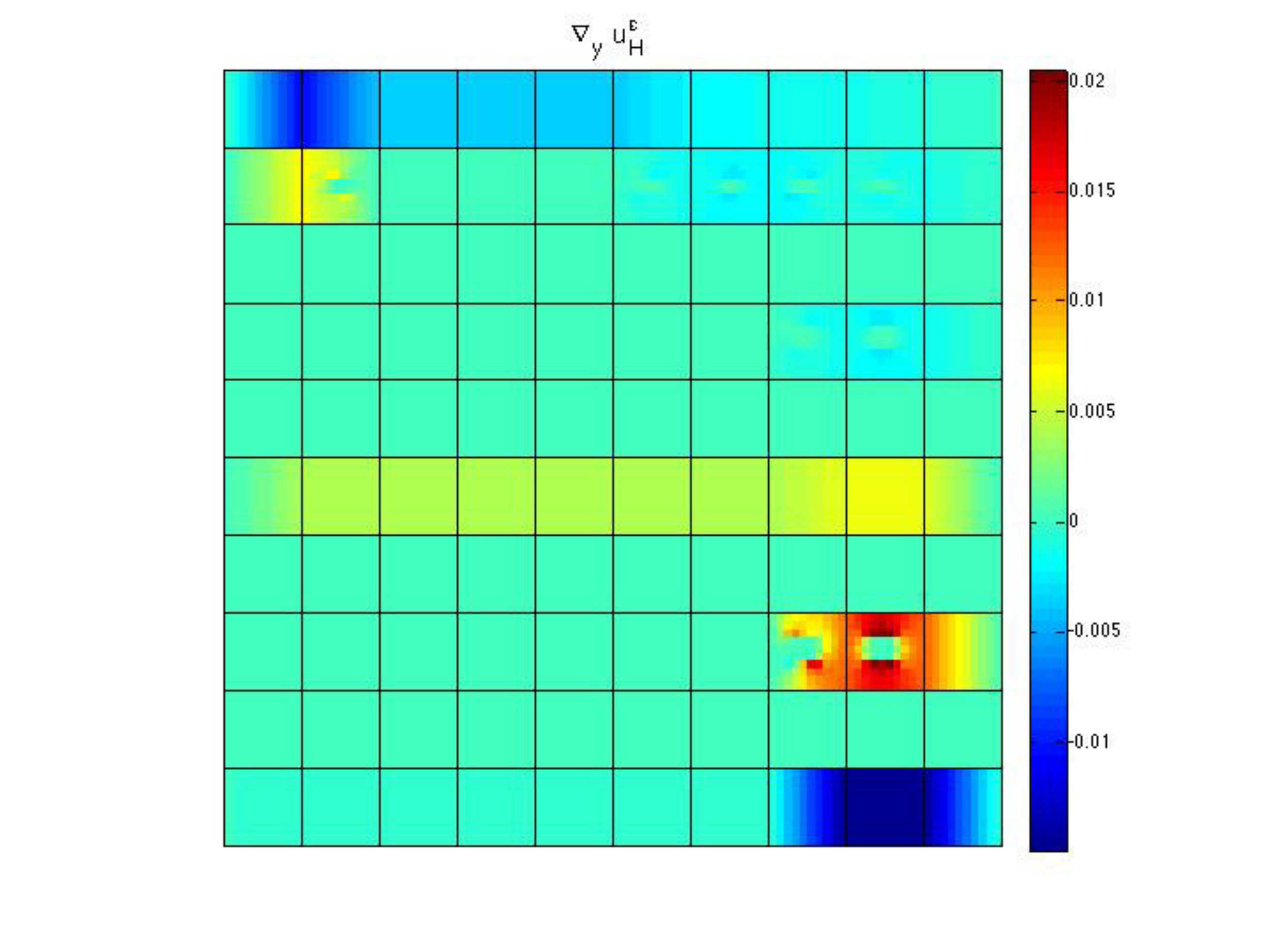}
\end{center}
\caption{Flow problem: exact gradient $\nab u^\eps$ (left) and MsFEM gradient $\nab u^\eps_H$ (right) (top row: components $\nab u^\eps \cdot \be_1$ and $\nab u_H^\eps \cdot \be_1$; bottom row: components $\nab u^\eps \cdot \be_2$ and $\nab u_H^\eps \cdot \be_2$).}
\label{fig:solMsFEM2Dflowgrad}
\end{figure}

\begin{figure}[H]
\begin{center}
\includegraphics[width=75mm]{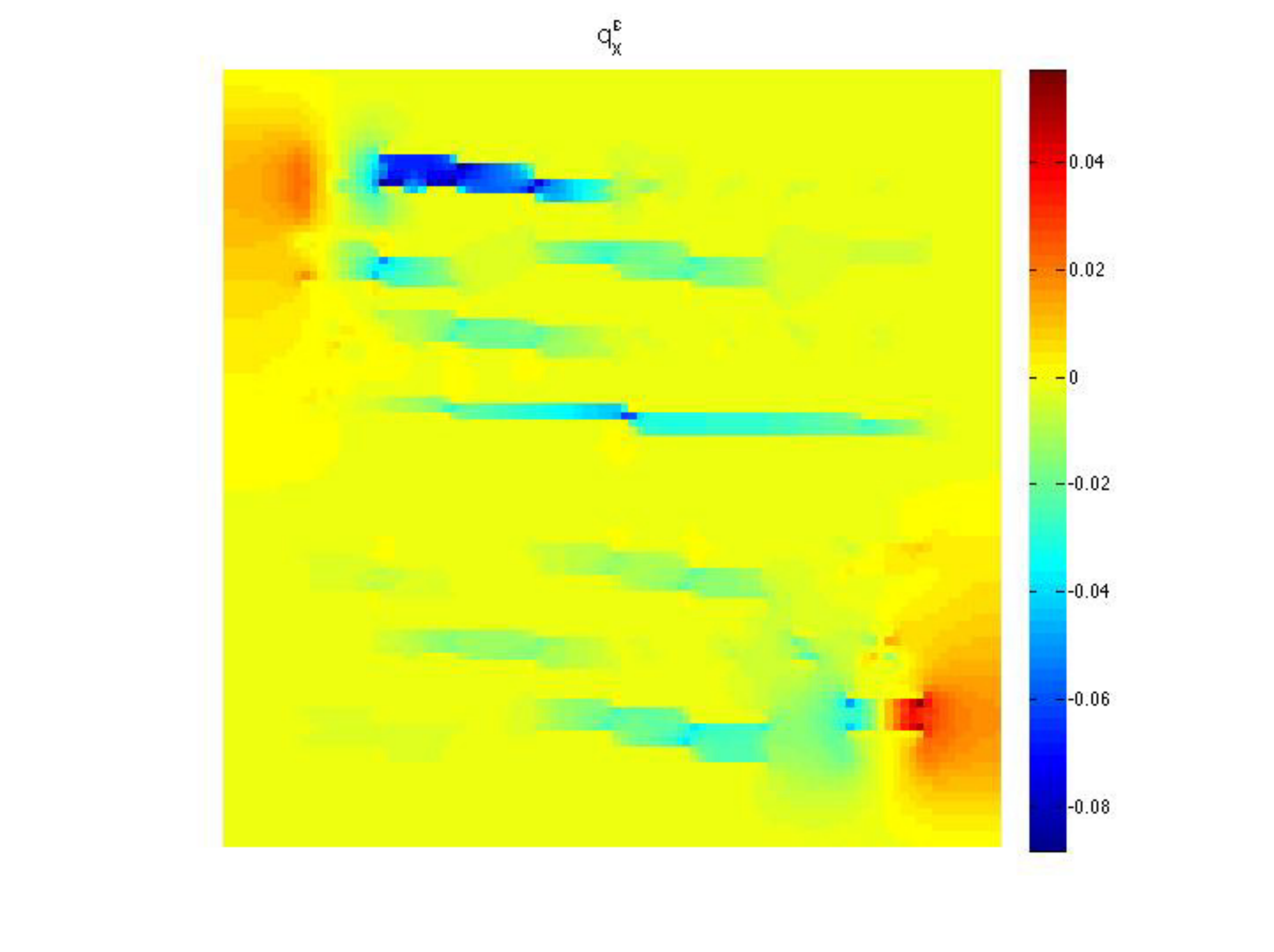}
\includegraphics[width=75mm]{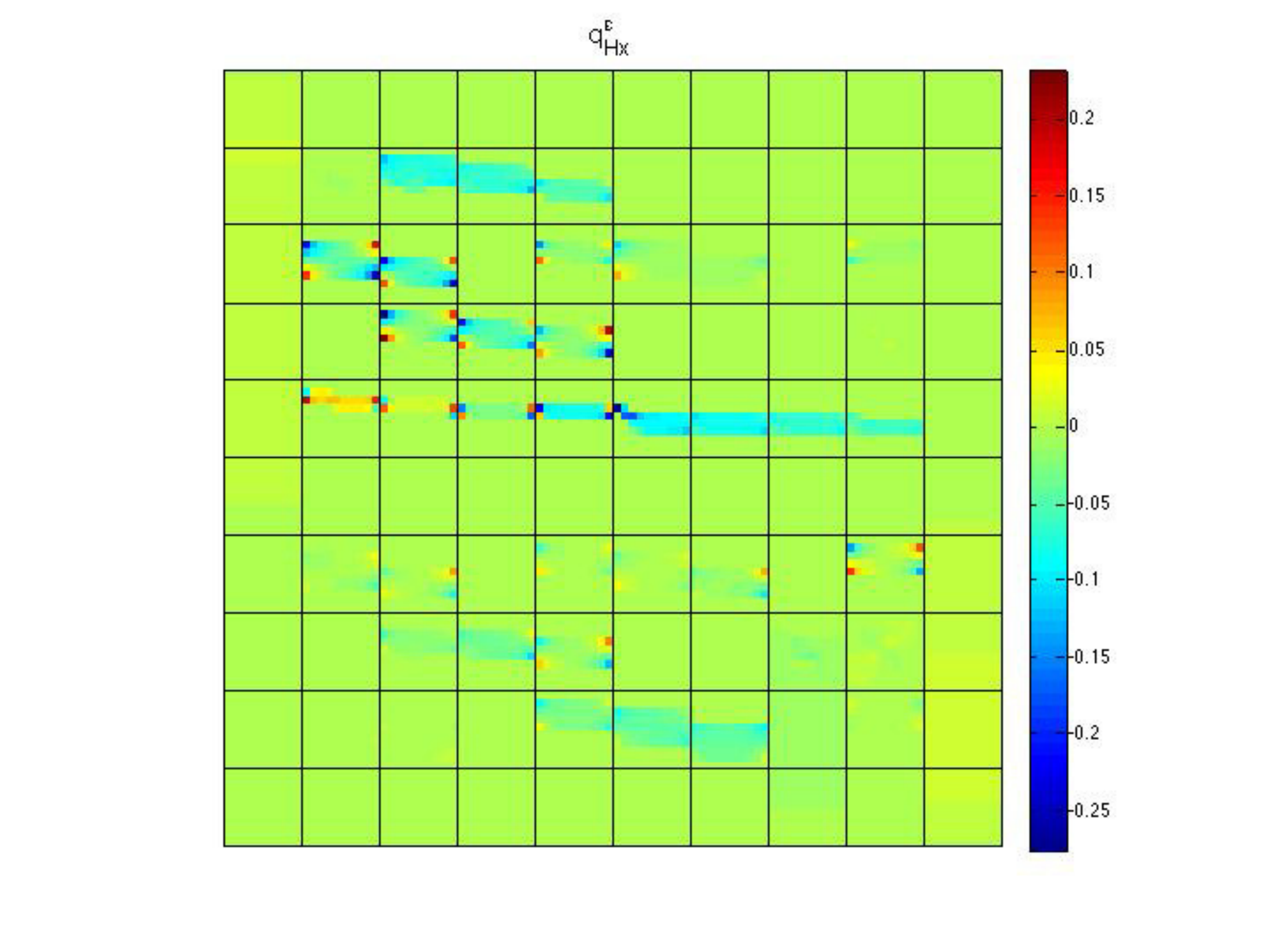}\\
\includegraphics[width=75mm]{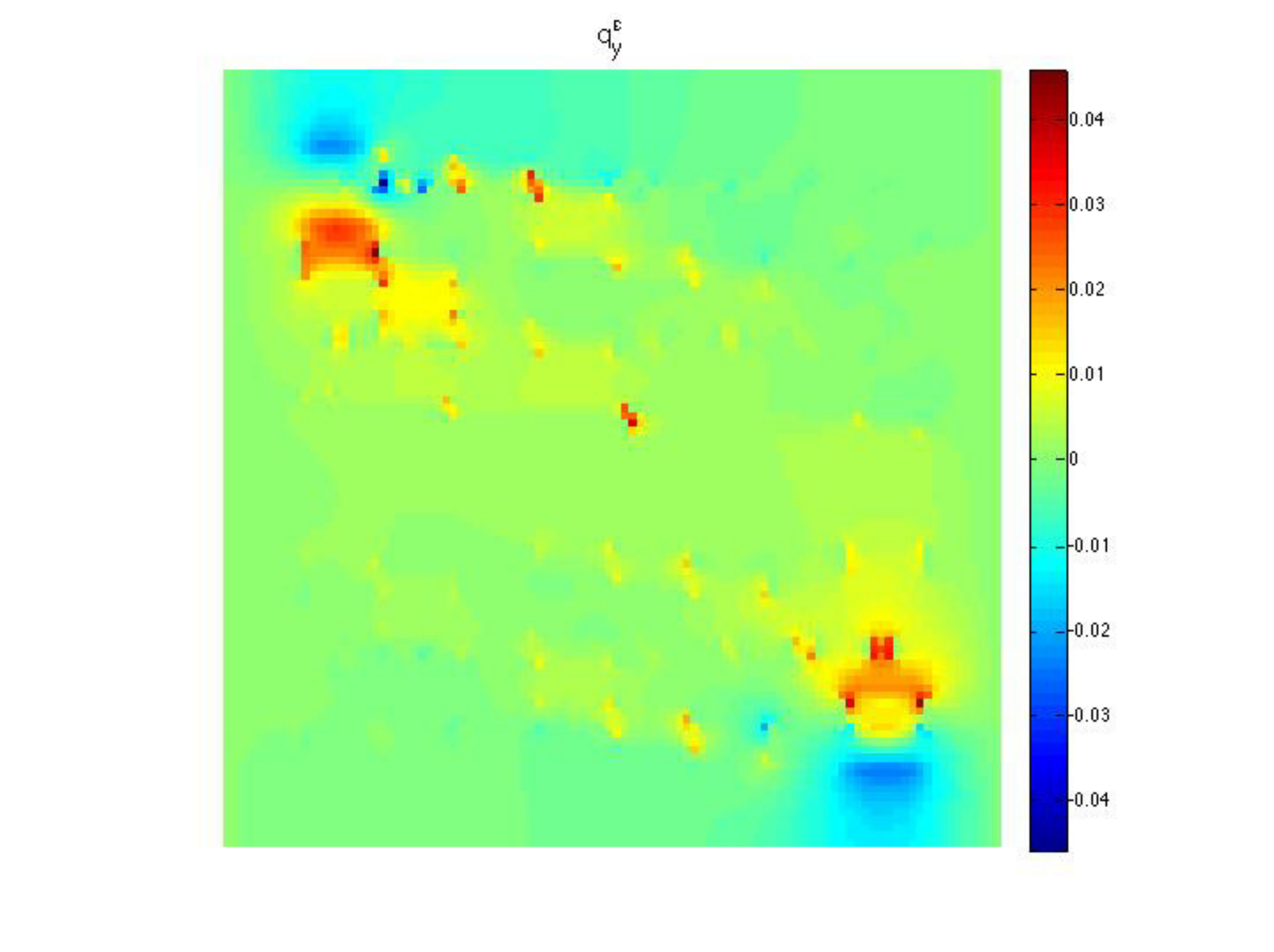}
\includegraphics[width=75mm]{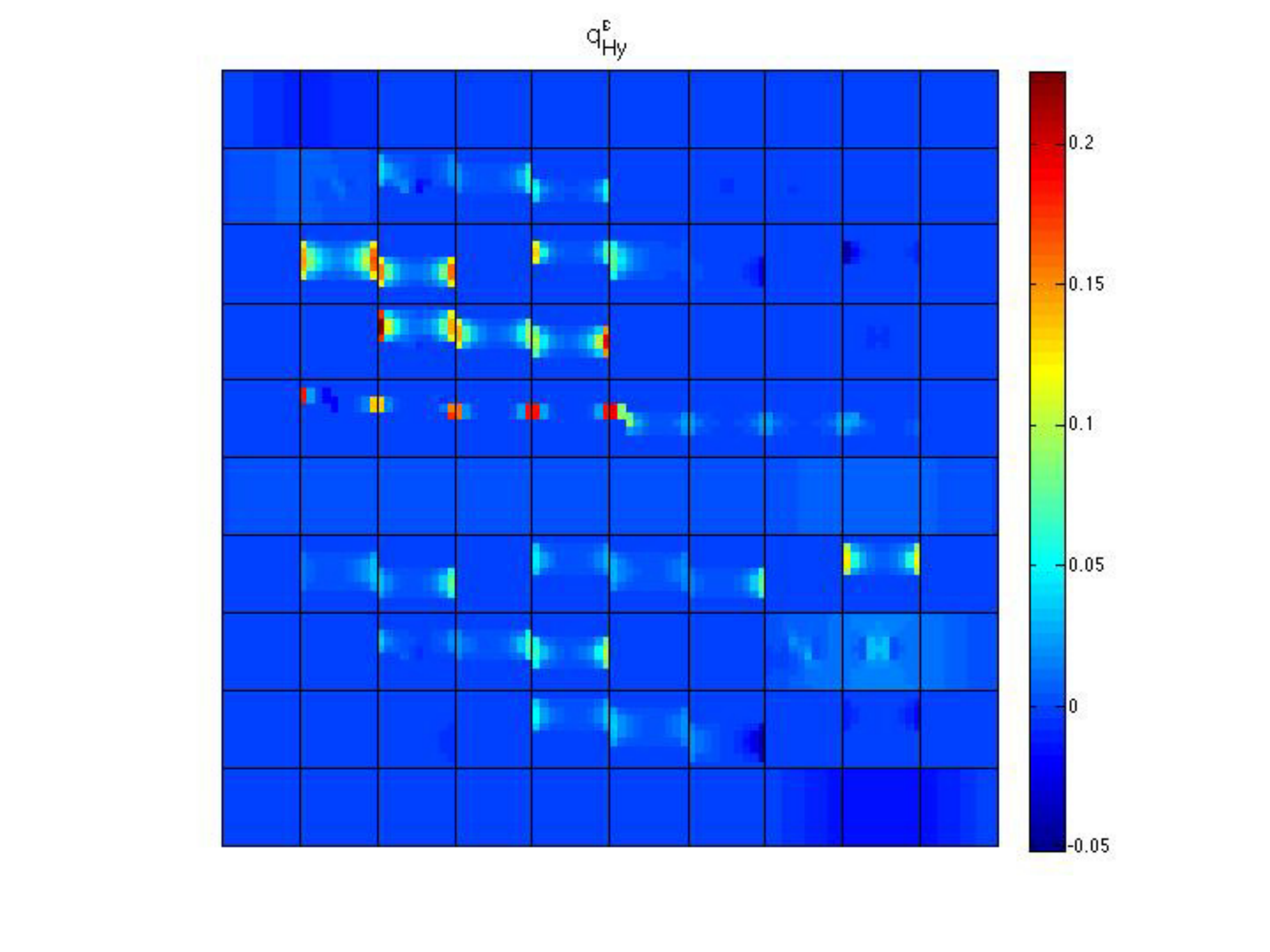}
\end{center}
\caption{Flow problem: exact flux $\bq^\eps$ (left) and MsFEM flux $\bq^\eps_H$ (right) (top row: components $\bq^\eps \cdot \be_1$ and $\bq_H^\eps \cdot \be_1$; bottom row: components $\bq^\eps \cdot \be_2$ and $\bq_H^\eps \cdot \be_2$).}
\label{fig:solMsFEM2Dflowflux}
\end{figure}

\begin{figure}[h]
\begin{center}
\includegraphics[width=75mm]{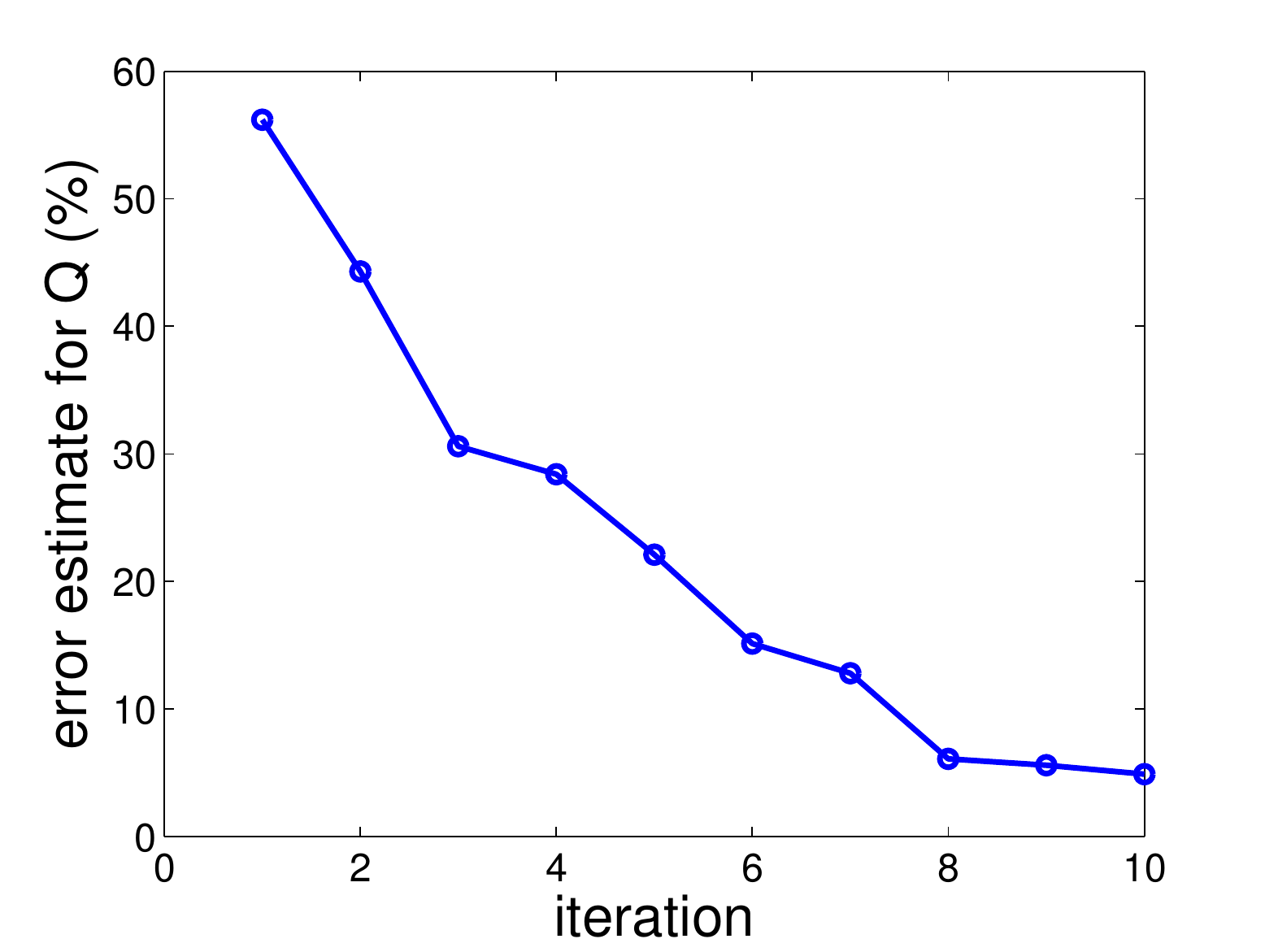} 
\end{center}
\caption{Flow problem: convergence of the error estimate for $Q(u^\eps)$.}
\label{fig:converrorQ2Dflow}
\end{figure}

We take an initial mesh $\mT_H$ made of $10\times 10$ macro elements, with $h_K=H/10$ for any $K$ and no oversampling. We show on Figs~\ref{fig:solMsFEM2Dflowu}, \ref{fig:solMsFEM2Dflowgrad} and~\ref{fig:solMsFEM2Dflowflux} the exact solution $u^\eps$ and its MsFEM approximation $u^\eps_H$, their gradients $\nab u^\eps$ and $\nab u^\eps_H$, and the fluxes $\bq^\eps$ and $\bq^\eps_H$. We can see that the solution $u^\eps$ is essentially supported in the source region $K_1$ and the sink region $K_2$. The initial approximation $u^\eps_H$ is not very accurate, especially in the top half of $\Omega$. As can be seen on Fig.~\ref{fig:converrorQ2Dflow}, the initial error on the quantity of interest is almost 60\%. The discretization is purposely too coarse.  

\begin{figure}[H]
\begin{center}
\includegraphics[width=75mm]{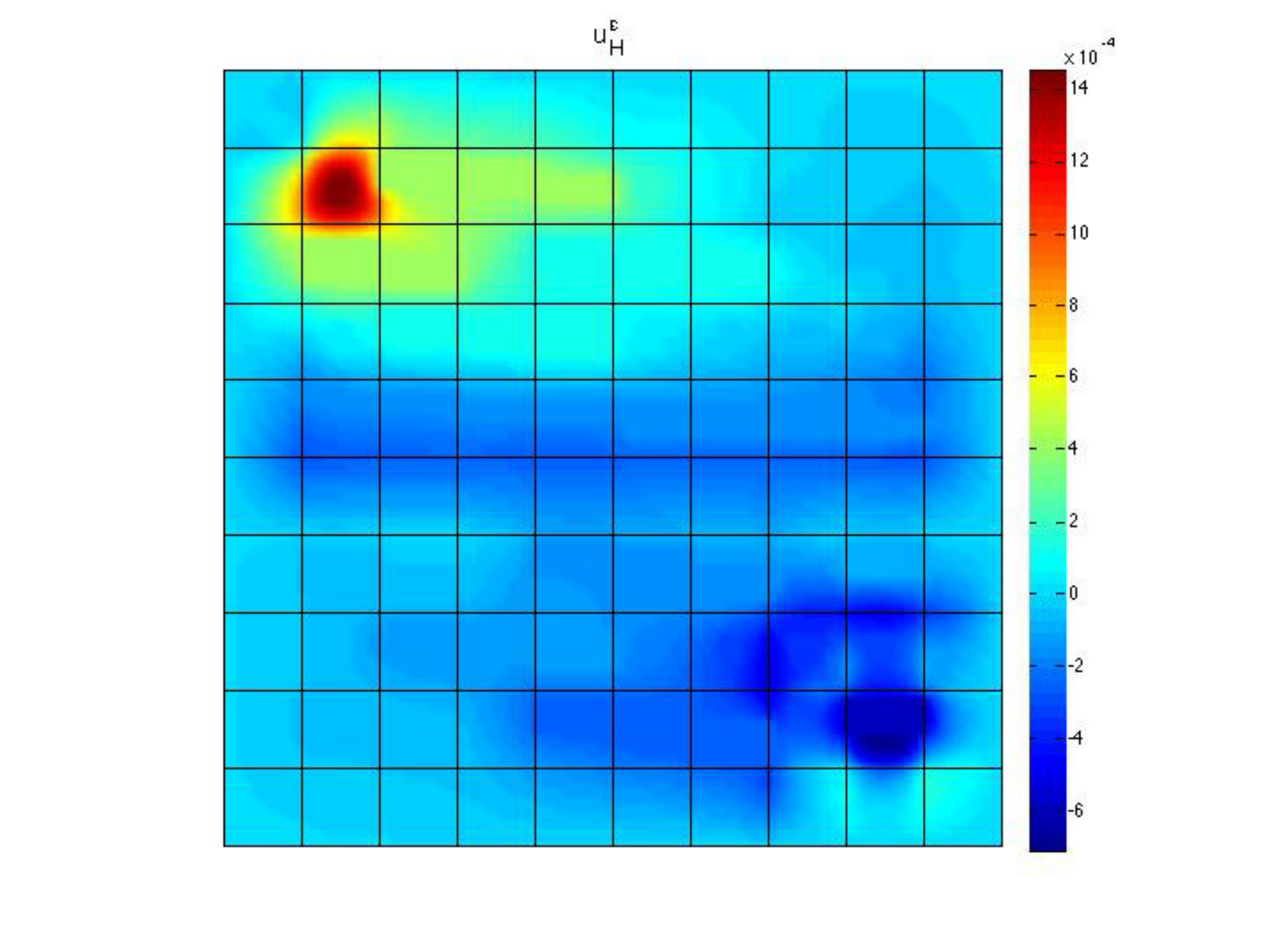}    
\includegraphics[width=75mm]{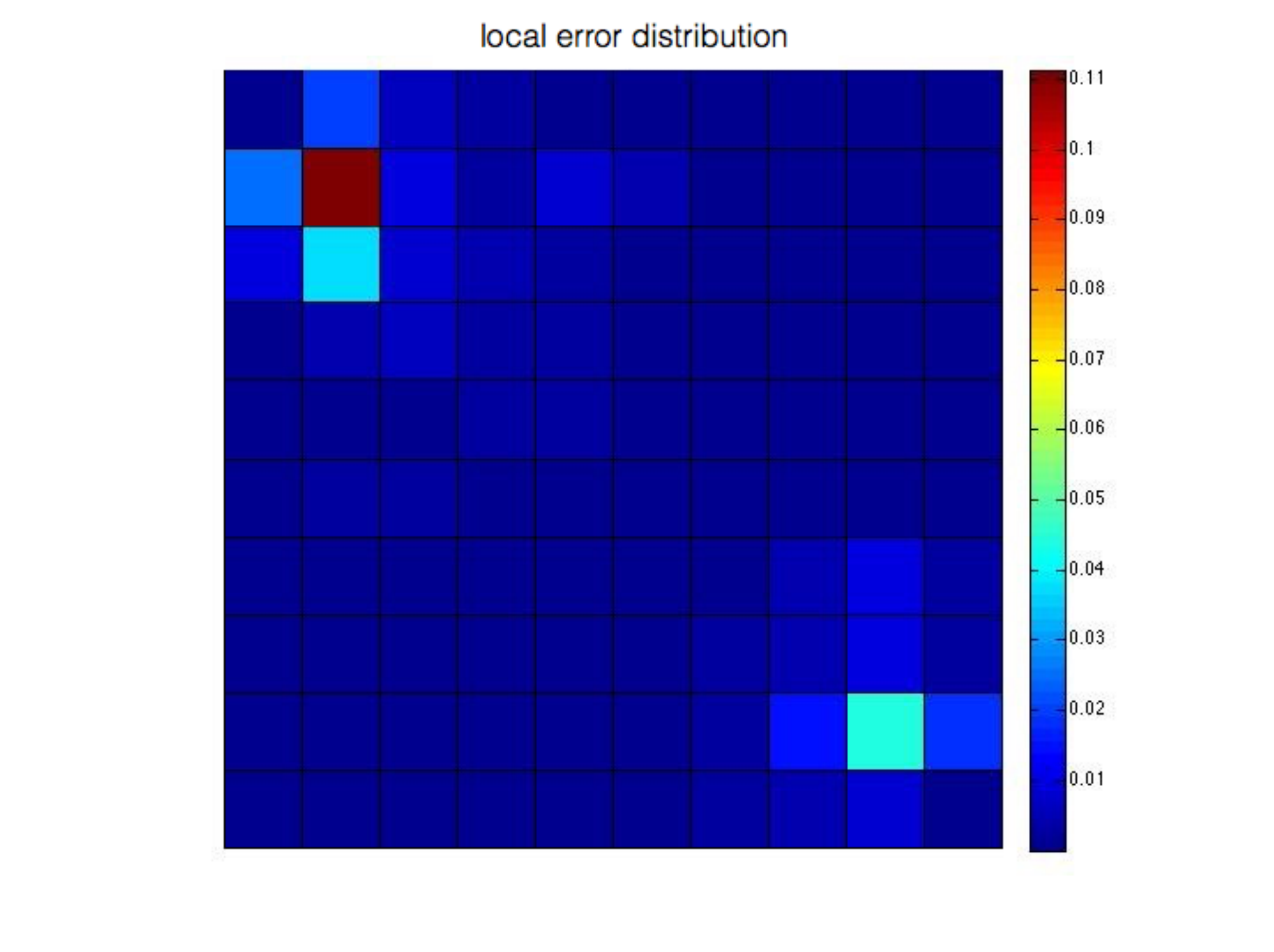}
\end{center}
\caption{Flow problem: distribution of the error. Left: $u^\eps -u^\eps_H$. Right: error in the energy norm.}
\label{fig:ener2Dflow}
\end{figure}

The error between the exact solution and its MsFEM approximation is shown on Fig.~\ref{fig:ener2Dflow}.

\begin{figure}[H]
\begin{center}
\includegraphics[width=50mm]{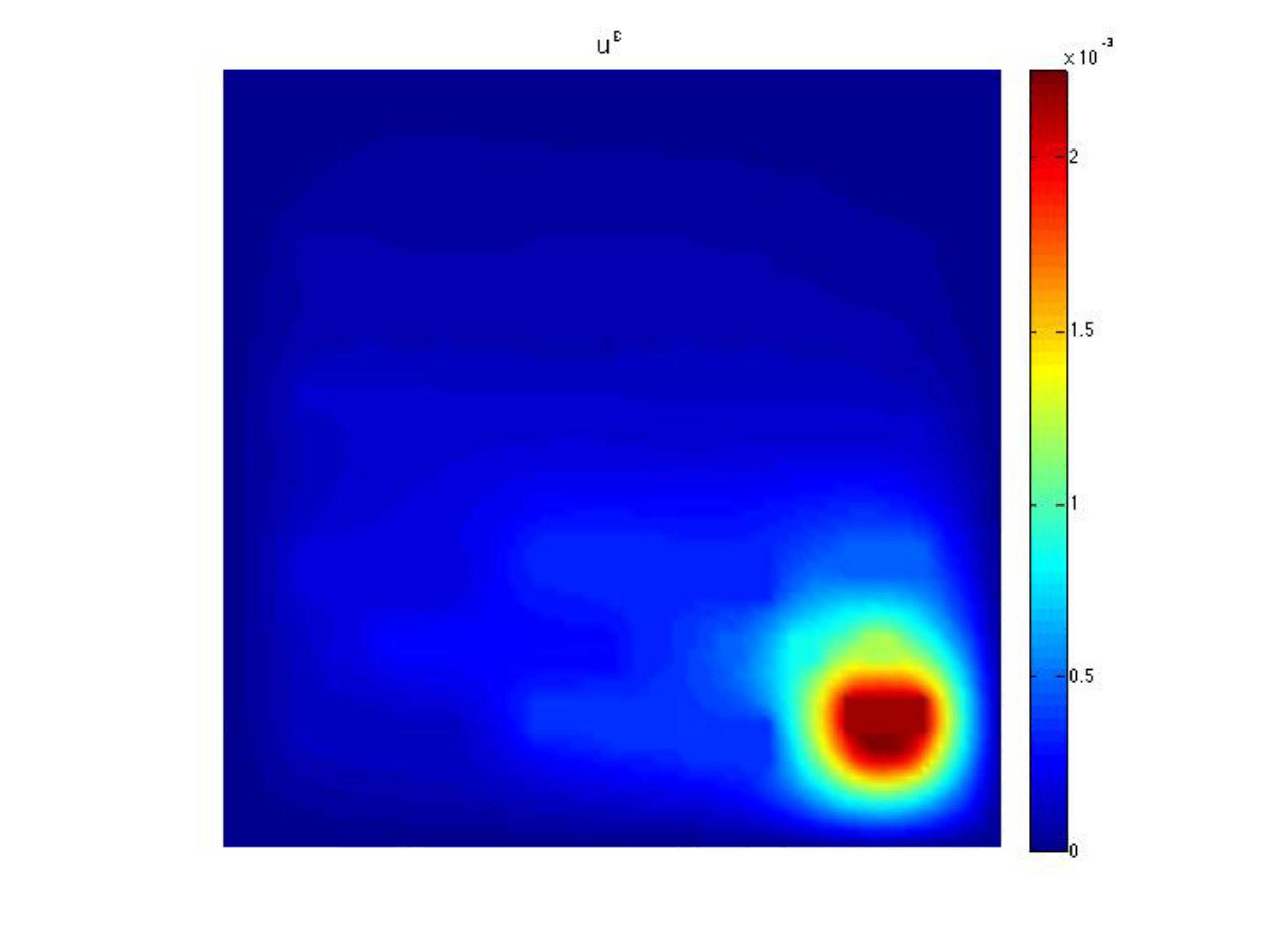}
\includegraphics[width=50mm]{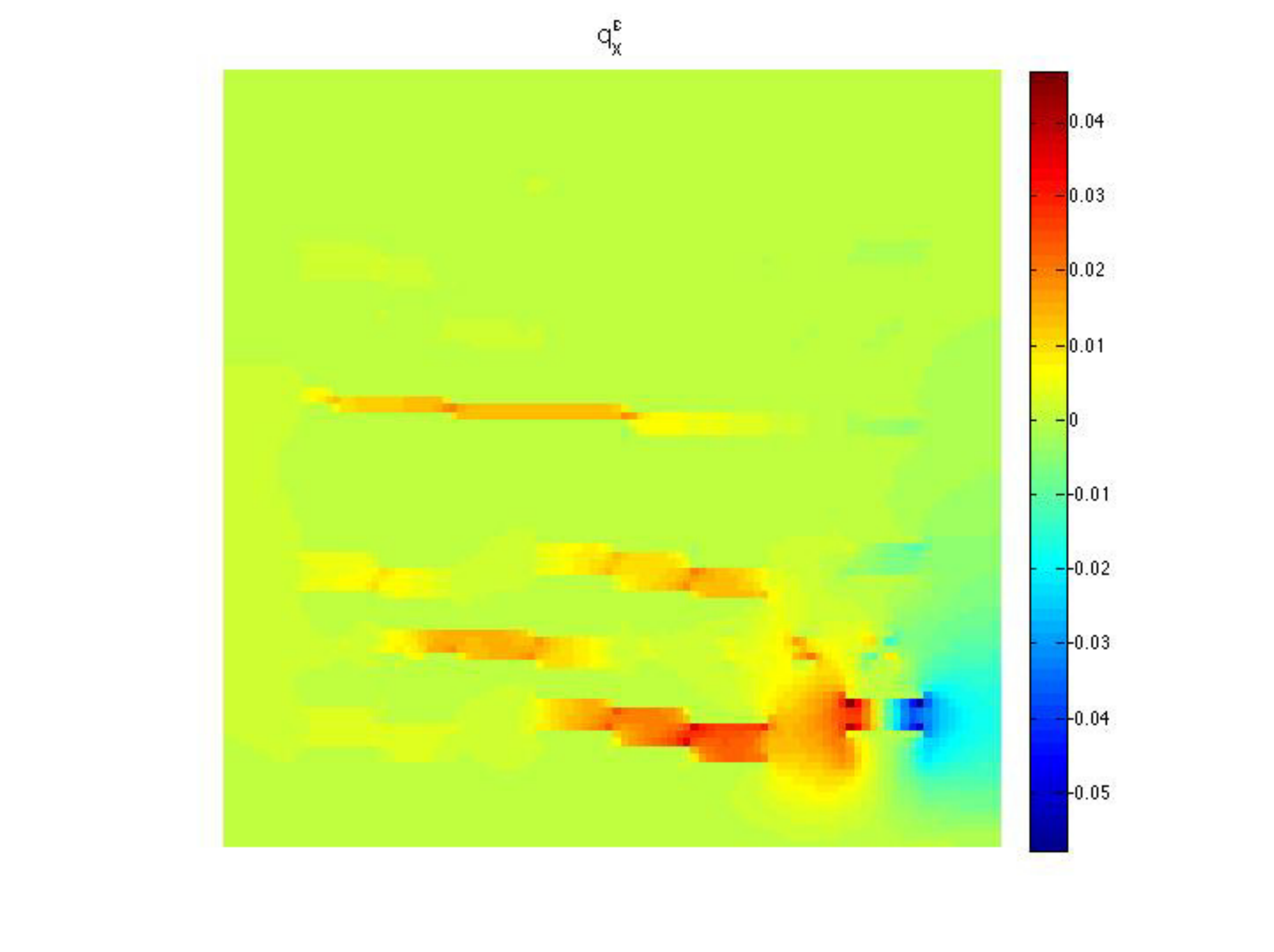}
\includegraphics[width=50mm]{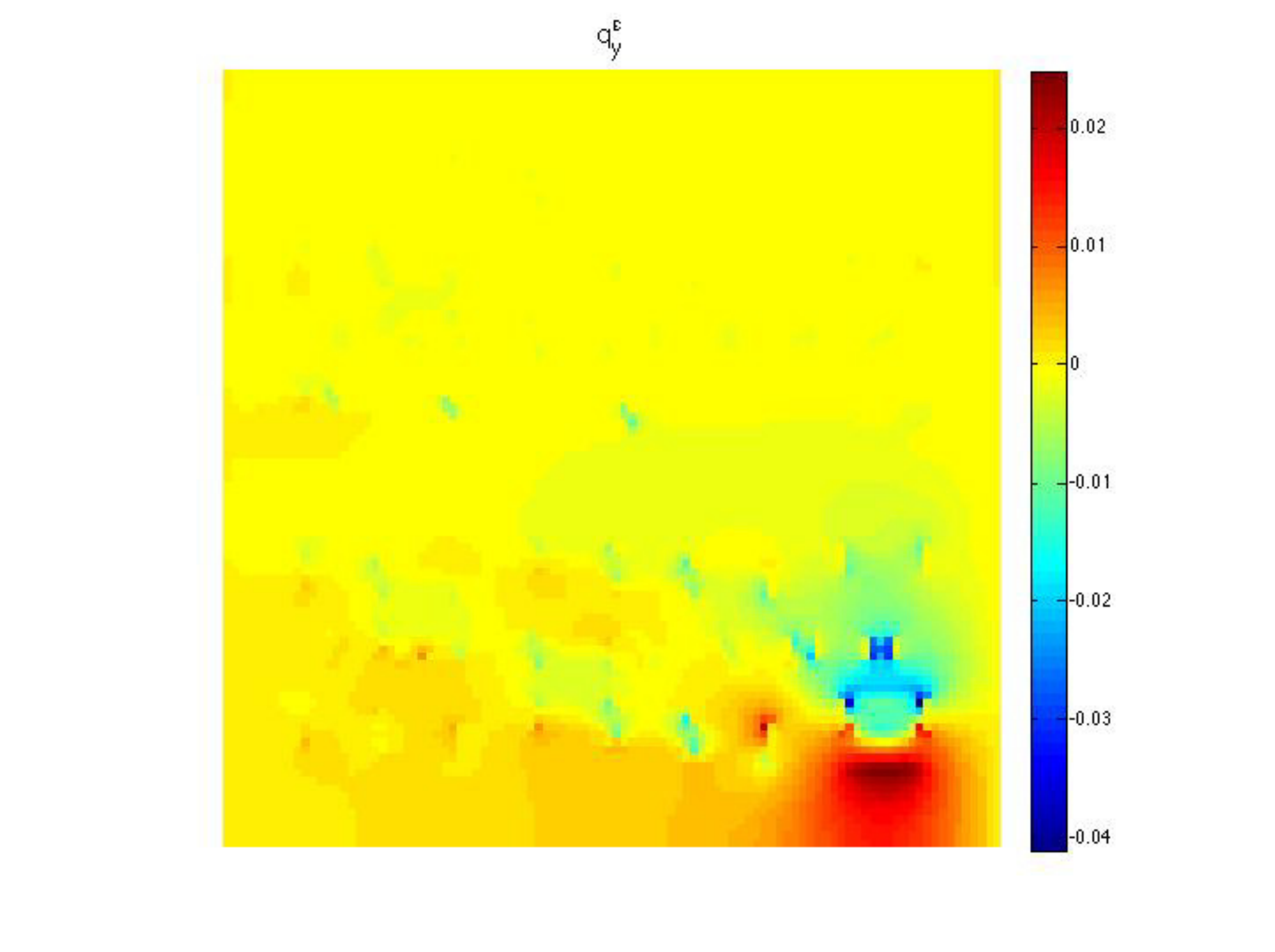}
\end{center}
\caption{Flow problem: exact adjoint solution $\widetilde{u}^\eps$ (left), exact adjoint flux $\widetilde{\bq}^\eps \cdot \be_1$ (center) and $\widetilde{\bq}^\eps \cdot \be_2$ (right).}
\label{fig:soladj2Dflowu}
\end{figure}

The adjoint solution $\widetilde{u}^\eps$ and flux $\widetilde{\bq}^\eps$ corresponding to the quantity of interest $Q$ are shown on Fig.~\ref{fig:soladj2Dflowu}.

\begin{figure}[h]
\begin{center}
\includegraphics[width=150mm]{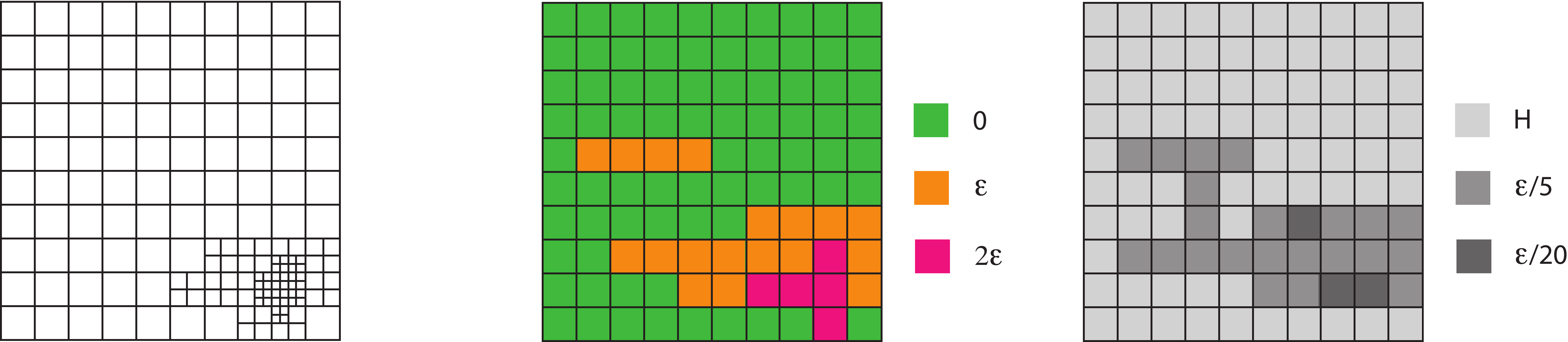} 
\end{center}
\caption{Flow problem: final MsFEM discretization for the quantity of interest $Q$.}
\label{fig:meshconvQ2Dflow}
\end{figure}

\begin{figure}[h]
\begin{center}
\includegraphics[width=150mm]{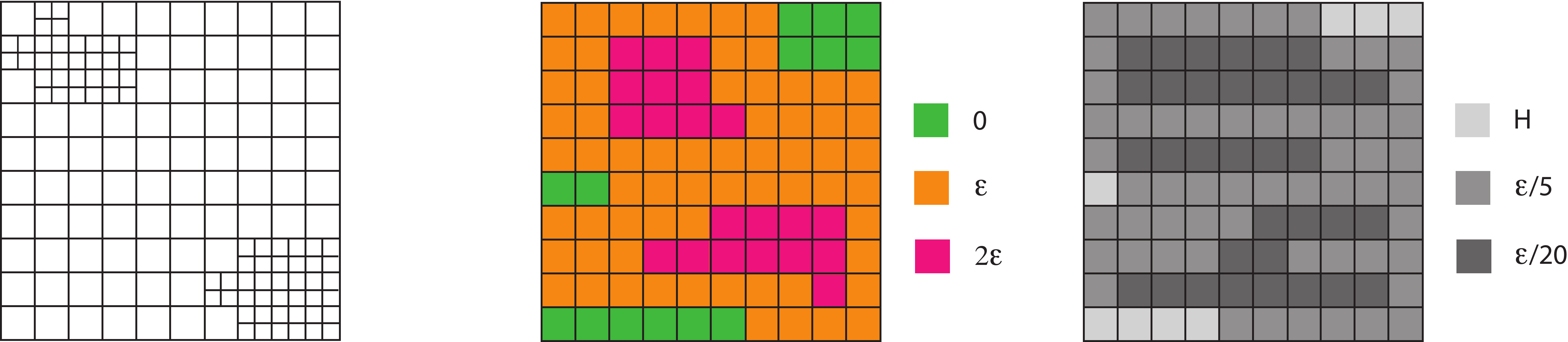} 
\end{center}
\caption{Flow problem: final MsFEM discretization when controlling the global error.}
\label{fig:meshconvglob2Dflow}
\end{figure}

The error threshold is fixed to $1\%$. We show on Fig.~\ref{fig:meshconvQ2Dflow} the final MsFEM discretization that we obtain in order to respect the error threshold for the quantity of interest. As in the previous test cases, this discretization is much coarser than the discretization (shown on Fig.~\ref{fig:meshconvglob2Dflow}) that is obtained when controlling the global error (in energy norm) using tools presented in~\cite{CHA16b}.

The accuracy of the error estimate is again very good in this example. We get $\eta^Q/|Q(u^\eps)-Q(u^\eps_H)| = 1.14$ at the first iteration of the adaptive algorithm, and $\eta^Q/|Q(u^\eps)-Q(u^\eps_H)| = 1.08$ at the final iteration of the adaptive algorithm.

The error estimate on the quantity of interest is shown on Fig.~\ref{fig:converrorQ2Dflow}, as the iterations of the adaptive process proceed. We again observe that this error monotically decreases.

\section{Conclusions and perspectives}\label{section:conclusions}

We have introduced efficient numerical tools for the \textit{a posteriori} error estimation of specific quantities of interest in the context of MsFEM computations. These tools are based on the Constitutive Relation Error concept, that we have previously used in~\cite{CHA16b} to estimate the error in the energy norm, and that we have here successfully extended to the goal-oriented context.

Our approach leads to fully computable and guaranteed upper bounds on the error, and enables to perform adaptive multiscale computations by selecting the relevant MsFEM parameters to reach a prescribed error tolerance. Our procedure is consistent with the {\em offline}/{\em online} paradigm of MsFEM, hence yielding affordable additional computational costs. 

The performances of the approach have been investigated on several multiscale problems. Fine scale features of the solution are adaptively recovered during the iterations, if and where needed. In particular, the fine-scale features of the solution in some parts of the physical domain may not affect the quantity of interest, and our procedure then prevents the discretization approach to capture them, thus achieving a significant computational gain at no loss in accuracy. Our method thus outperforms the standard \textit{a posteriori} error estimation strategy (which assesses and controls the error in the energy norm): for an equal tolerance on the error, a coarser discretization is obtained, leading to a smaller computational cost.

\section*{Acknowledgements}

We thank C. Le Bris for stimulating discussions on the topics addressed in this article. LC thanks Inria for enabling his two-year leave (2014--2016) in the MATHERIALS project team. The work of FL is partially supported by ONR under grant N00014-15-1-2777 and by EOARD under grant FA9550-17-1-0294.


\end{document}